\newif\ifSHOWEXTRA
\SHOWEXTRAtrue
\SHOWEXTRAfalse

\documentclass[12pt]{article}
\usepackage[utf8]{inputenc}
\usepackage[margin=1in]{geometry}
\usepackage{amsmath,amsthm,amssymb,amsfonts, pifont, bbm, xcolor, scrextend}
\definecolor{forestgreen}{RGB}{0,98,51}
\usepackage{appendix}
\usepackage{fancyhdr}
\usepackage{color,colortbl}
\newcounter{CommentCounter}
\newcommand{\daniel}[1]{{\color{purple} \sf $\spadesuit\spadesuit\spadesuit$ Daniel (\stepcounter{CommentCounter}\theCommentCounter): [#1]}}

\definecolor{Gray}{gray}{0.9}
\definecolor{LightGreen}{rgb}{0,0.1,0.1}
\newcommand{\R}{\mathbb R}
\newcommand{\G}{\mathcal G}
\newcommand{\Hg}{\mathcal H}
\newcommand{\Pa}{\mathbb P}
\newcommand{\Ea}{\mathbb E}

\newcommand{\N}{\mathbb N}
\newcommand{\Q}{\mathbb Q}
\newcommand{\Z}{\mathbb Z}
\newcommand{\0}{\vec 0}
\newcommand{\best}{\ast}
\newcommand{\one}{\mathbbm{1}}
\newcommand{\tepsilon}{\varepsilon}
\newcommand{\tdelta}{\delta}
\newcommand{\given}{\hspace{3pt}\vline\hspace{3pt}}
\newcommand\indep{\protect\mathpalette{\protect\independenT}{\perp}}
\def\independenT#1#2{\mathrel{\rlap{$#1#2$}\mkern2mu{#1#2}}}
\usepackage{tikz}
\usepackage{pgfplots}
\usepackage{eqnarray}
\usepackage{outlines}
\usepackage[mathscr]{euscript}
 \let\mathscr\relax
\usepackage[scr]{rsfso}
\usepackage{multicol}
\usepackage{enumerate}
\usepackage[colorlinks=true, pdfstartview=FitV, linkcolor=blue,
citecolor=blue, urlcolor=blue]{hyperref}

\newcommand{\E}{\mathrm{E}}

\newcommand{\BIX}{\overline{\bf X}}
\newcommand{\biX}{\overline{X}}

\newcommand{\totheleft}{-}
\newcommand{\totheright}{+}
\newcommand{\allowable}{\mathcal N}

\newtheorem{thm}{Theorem}[section]
\numberwithin{thm}{section}
\newtheorem{cor}[thm]{Corollary}
\newtheorem{prop}[thm]{Proposition}
\newtheorem{lem}[thm]{Lemma}

\newtheorem{quest}{Question}
\numberwithin{quest}{section}

\theoremstyle{remark}
\newtheorem{rem}{Remark}
\numberwithin{rem}{section}

\theoremstyle{definition}

\newtheorem{exmp}{Example}[section]

\newtheorem{claim}{Claim}
\numberwithin{claim}{exmp}
\newtheorem*{property}{Property}

\newcounter{ExampleSave}[section]
\newcounter{SectionSave}

\pgfplotsset{compat=1.12}

\begin{document}
\author{
Daniel J. Slonim\footnote{Purdue University,
  Department of Mathematics,
150 N. University Street, West Lafayette, IN 47907, dslonim@purdue.edu, 0000-0002-9554-155X}
}
\title{Random Walks in Dirichlet Random Environments on $\Z$ with Bounded Jumps}

\date{\today}

\maketitle

\begin{abstract}
    We examine a class of random walks in random environments on $\Z$ with bounded jumps, a generalization of the classic one-dimensional model. The environments we study have i.i.d. transition probability vectors drawn from Dirichlet distributions. For this model, we characterize recurrence and transience, and in the transient case we characterize ballisticity. For ballisticity, we give two parameters, $\kappa_0$ and $\kappa_1$. The parameter $\kappa_0$ governs finite trapping effects, and $\kappa_1$ governs repeated traversals of arbitrarily large regions of the graph. We show that the walk is right-transient if and only if $\kappa_1>0$, and in that case it is ballistic if and only if $\min(\kappa_0,\kappa_1)>1$. 
    
    \medskip\noindent {\it MSC 2020:}
    60G50 
    60J10 
    60K37 
    \\
    {\it Kewords:}
    random walk,
    random environment,
    Dirichlet environments,
    bounded jumps,
    ballisticity
\end{abstract}

\tableofcontents

\section{Introduction}\label{sec:Introduction}

This paper studies a special model of random walks in random environments (RWRE) on $\Z$. We allow jumps of bounded size, and assume transition probabilities at various sites are drawn in an i.i.d. way according to a Dirichlet distribution. 

RWRE were first treated in depth by F. Solomon \cite{Solomon1975} in 1975 and the study has since grown in many directions. Solomon focused on nearest-neighbor RWRE on $\Z$, characterizing directional transience and calculating limiting speed (under mild conditions) in terms of simple, easily computable expectations involving the environment at a single site. One of the main results from Solomon's paper is that a RWRE can approach infinity almost surely, but with an almost-sure limiting speed of zero, a surprising phenomenon that cannot occur for random walks in homogenous or even periodic environments. 
A random walk is called {\em directionally transient} if it has an almost-sure limiting direction. Under appropriate i.i.d. assumptions on the environment, all directionally transient nearest-neighbor RWRE on $\Z^d$ have an almost-sure limiting velocity $v$ (which can be 0). This is also true of RWRE with bounded jumps on $\Z$; see Appendix \ref{append:Leftoverproofs}. We say a directionally transient random walk is {\em ballistic} if the limiting velocity $v$ is nonzero. The goal of this paper is to characterize directional transience and ballisticity for the model of random walks in Dirichlet environments (RWDE) on $\Z$ with bounded jumps. 

RWRE have proven quite challenging to analyze in settings other than the nearest-neighbor case of $\Z$. 
For instance, although some sufficient conditions for directional transience and ballisticity have been studied for  nearest-neighbor RWRE on $\Z^d$, no general characterizations are known, and known sufficient conditions are often quite difficult to check in practice. Nevertheless, certain special cases have proven to be more tractable. Examples include random environments that almost surely have zero drift at every site (e.g., \cite{Lawler1982}, \cite{Guo&Zeitouni2012}, \cite{Berger&Deuschel2014}), random environments that are small perturbations of simple random walks (e.g., \cite{Bricmont&Kupiainen1991}, \cite{Sznitman2003}, \cite{Sabot2004}, \cite{Bolthausen&Zeitouni2007}, \cite{Baur&Bolthausen2015}), and environments where transition probabilities are deterministic in some directions and random in others (e.g., \cite{Bolthausen&Sznitman&Zeitouni2003}). The case of RWDE is another notable example. Many conjectures that remain open for general nearest-neighbor RWRE on $\Z^d$ have been proven for RWDE, including a characterization of directional transience for all $d$ and a characterization of ballisticity for $d\geq3$.

Dirichlet environments present a helpful case study for general RWRE, giving insight into what sorts of behaviors are possible. For example, a noted conjecture (see \cite{Drewitz&Ramirez2010}) asserts that under an assumption of uniform ellipticity (where all transition probabilities are bounded from below), nearest-neighbor RWRE on $\Z^d$, $d\geq2$, that are directionally transient are necessarily ballistic. Certain RWDE provide a counterexample for the non--uniform elliptic case \cite{Sabot2013}, showing that the uniform ellipticity assumption is necessary for this conjecture. At the same time, in the Dirichlet case, the factors that can cause ballisticity to fail are entirely due to the non--uniform ellipticity property, providing some additional evidence for the conjecture (see \cite{Bouchet2013}, \cite[Remark 5.13]{Sabot&Tournier2016}). 

In the setting of RWRE on $\Z$ with bounded jumps---and of RWRE on a strip, a generalization of the bounded-jump model---directional transience and ballisticity have been given various characterizationsy, most in terms of Lyapunov exponents of infinite products of random matrices (see, for example, \cite{Key1984},\cite{Bolthausen&Goldsheid2000},\cite{Bremont2004},\cite{Bremont2009}, \cite{Roitershtein2008}). These exponents cannot in general be computed exactly, although some of them can be well approximated. Moreover, the characterizations of ballisticity have relied on various strong ellipticity assumptions, which preclude the Dirichlet model.

Although RWDE have proven to be a fruitful model for nearest-neighbor RWRE on $\Z^d$, we know of only one paper studying RWDE on $\Z$ with bounded jumps or on strips. That paper \cite{Keane&Rolles2002} was the first to treat RWDE in {\em any} setting. It uses a connection between RWDE and directed edge reinforced random walks to provide a characterization of directional transience for the latter in the mold of \cite{Key1984}. That paper preceded the development of helpful tools that have since been applied to the analysis of nearest-neighbor RWDE on $\Z^d$. Armed with these tools, we now undertake a more comprehensive study of a special case of RWDE on a strip: namely, RWDE on $\Z$ with bounded jumps. 

We show that as with nearest-neighbor RWDE on $\Z^d$, directional transience is determined by the direction of the annealed expectation of the first step, a fact that does not hold for general RWRE. We also show that in the case where the annealed expectation of the first step is 0, the walk is recurrent.
We next turn to ballisticity. For i.i.d. RWRE on $\Z$ with bounded jumps, there necessarily exists a deterministic $v$ with $\lim_{n\to\infty}\frac{X_n}{n}=v$ almost surely (see Appendix \ref{append:Leftoverproofs}). When the walk is recurrent, necessarily $v=0$. We assume the walk is transient to the right, so that $v\geq0$, and characterize the ballistic regime $v>0$.

Our results on ballisticity comprise a sort of mixture of characterizations for the nearest-neighbor cases on $\Z$ and on $\Z^d$, $d\geq3$. These characterizations are quite different from each other, as they reflect substantially different ways that a walk may ``get stuck." (The case $d=2$ is still open.)
For $d\geq3$, ballisticity has a simple characterization in terms of a parameter $\kappa_0$ (called $\kappa$ in \cite{Bouchet2013}).
The fact that Dirichlet distributions are not uniformly elliptic allows environments to contain arbitrarily severe ``traps" where a walk can get stuck for a long time. In fact, if enough parameters of the Dirichlet distribution are sufficiently small, there are finite subgraphs whose annealed expected exit times are infinite, causing zero limiting speed---see \cite[Proposition 12]{Tournier2009}, \cite[Proposition 2]{Sabot&Tournier2016}. The parameter $\kappa_0$, which represents the minimal amount of weight exiting a finite set, controls finite moments of the quenched expected exit times of finite traps containing the origin. In the case $d\geq3$, a directionally transient walk is ballistic if and only if $\kappa_0>1$, reflecting the idea that finite traps are the only way directional transience with zero speed can occur in the case $d\geq3$. 

The case $d=1$, 
where probabilities $\omega(x,x+1)$ of stepping to the right are given by beta random variables with parameters $(a,b)$, is different. Here, $\kappa_0=a+b$, and this parameter still controls finite traps, but 
it is possible for a walk to have zero speed even if $\kappa_0>1$. In fact, ballisticity of walks transient to the right is controlled by another parameter, $\kappa_1=a-b$, which is the unique positive number (studied in a more general setting by Kesten, Kozlov, and Spitzer in \cite{KestenKozlov&Spitzer1975}, and there called $\kappa$) such that  $E\left[\left(\frac{1-\omega(0,1)}{\omega(0,1)}\right)^{\kappa_1}\right]=1$. In the case $d=1$, a walk that is transient to the right is ballistic if and only if $\kappa_1>1$, by a direct application of the characterization of ballisticity given in \cite{Solomon1975}. Here, as in the case $d\geq3$, $\kappa_0\leq1$ is enough to cause finite trapping that would slow the walk down to zero speed. However, because we always have $\kappa_1<\kappa_0$, the walk is already not ballistic in this case, and thus the value of $\kappa_1$ alone determines ballisticity.

We find that the parameters $\kappa_0$ and $\kappa_1$ can be given definitions that apply to our model as well, that unlike in the nearest-neighbor $\Z$ case, either may be greater than the other, and that both must be greater than 1 in order to achieve ballisticity. As in the nearest-neighbor cases, $\kappa_0$ controls finite traps. We do not give an explicit formula for $\kappa_0$, which is defined as an infimum over an infinite set of sums, but we show that it is in fact a minimum over finitely many sums, and provide an algorithm to compute it directly. The parameter $\kappa_1$ has a simple formula as a weighted sum of Dirichlet parameters, which reduces to $a-b$ in the nearest-neighbor beta case. 
When $\kappa_0\leq1$, the walk has zero speed because of the relatively high likelihood of getting stuck in a region of bounded size for a long time. When $\kappa_1\leq1$, the walk has zero speed because of the relatively high likelihood of repeatedly backtracking over regions of all sizes. When both are greater than 1, the walk is ballistic

The appearance of the dual possibilities of finite trapping and large-scale backtracking seems to be a new phenomenon in RWRE with bounded jumps. Previous characterizations of ballisticity use ellipticity assumptions strong enough to preclude finite trapping, and therefore do not cover cases where walks can get stuck in these two different ways. Part of our characterization for the Dirichlet case involves showing that the two forms of slowing ``act independently" in the sense that if neither finite trapping nor large-scale backtracking is by itself enough to cause zero speed, the two together cannot cause zero speed; the walk is ballistic. We ask whether this is true for general RWRE on $\Z$ with bounded jumps; see Question \ref{question}. 

While our main results are for RWDE, part of the proof requires obtaining some results that apply to general RWRE on $\Z$ with bounded jumps. We present a new abstract criterion for ballisticity, showing that the walk is ballistic if and only if the annealed expected number of returns to the origin is finite. To do this, we define a ``walk from $-\infty$ to $\infty$" in a typical environment where transience to the right holds, and the almost-sure limiting speed of this bi-infinite walk is the reciprocal of the expected amount of time it spends at 0 under an appropriate annealed measure. Thus, the limiting speed is zero precisely when this expectation is infinite. We show that this expectation in turn is infinite if and only if the expected amount of time at 0 for a walk started at 0 is infinite. While for general models we do not know how to check the finiteness of this expectation, we are able to do so for our Dirichlet model, showing that it is finite if and only if $\min(\kappa_0,\kappa_1)>1$.

\subsection{Model}\label{subsec:model}
Because there are many definitions and symbols introduced throughout this paper, we provide an appendix to help the reader keep track of notation that is introduced outside of this subsection or that is particularly unusual. 

Let $L$ and $R$ be positive integers, and let $(\alpha_{i})_{i=-L}^R$ be non-negative real numbers, with $\alpha_L,\alpha_R>0$. 
Assume that the $\gcd$ of all $i$ with $\alpha_i>0$ is 1. We are interested in a random walks in Dirichlet random environments on $\Z$ with jumps to the left up to $L$ steps and to the right up to $R$ steps, with transition probability vectors given by i.i.d. Dirichlet random vectors with parameters $(\alpha_{i})_{i=-L}^R$. To define these concepts explicitly, we will first define the more general notion of random walks in random environments.

{\bf Random Walks in Random Environments}

Let $V$ be a finite or countable set, and let $\Omega=\prod_{x\in V}\mathcal{M}_1(V)$, where $\mathcal{M}_1(V)$ is the set of probability measures on $V$, endowed with the topology of weak convergence. An {\em environment on $V$} is an element $\omega\in\Omega$, which can be thought of as a function from $V\times V$ to $[0,1]$, with $\sum_{y\in V}\omega(x,y)=1$ for all $x$. 
For a given environment $\omega$ and $x\in V$, we can define $P_{\omega}^x$ to be the measure on $V^{\N_0}$ giving the law of a Markov chain $(X_n)_{n=0}^{\infty}$ started at $x$ with transition probabilities given by $\omega$. That is, $P_{\omega}^x(X_0=x)=1$, and for $n\geq1$, $y\in V$, $P_{\omega}^x(X_{n+1}=y|X_0,\ldots,X_n)=\omega(X_n,y)$.

Let $\mathcal{F}_V$ be the Borel sigma field on $\Omega_V$ (with respect to the product topology), and let $P$ be a probability measure on $(\Omega_V,\mathcal{F}_V)$ (we often leave the $\mathcal{F}_V$ implicit and say $P$ is a probability measure on $\Omega_V$). For a given $x\in V$, we let $\Pa^x=P\times P_{\omega}^x$ be the measure on $\Omega_V\times V^{\N_0}$ induced by both $P$ and $P_{\omega}$. That is, for measurable events $A\subset\Omega_V,B\subset V^{\N_0}$,
\begin{equation*}
\Pa^x(A\times B)=\int_{A}P_{\omega}^x(B)P(d\omega)
\end{equation*}
In particular, $\Pa^x(\Omega_V\times B)=E\left[P_{\omega}^x(B)\right]$. 
For convenience, we commit a small abuse of notation by using $\Pa^x$ to refer both to the measure we've described on $\Omega_V \times V^{\N_0}$ and also to its marginal $\Pa^x(\Omega_V \times \cdot)$ on $V^{\N_0}$. 
We call a measure $P_{\omega}^x$ on $V^{\N_0}$ a {\em quenched measure} of a random walk in random environment on $V$ started at $x$, and we call the measure $\Pa^x$ 
the {\em annealed measure}. 

In this paper, $V$ is nearly always $\Z$ or a subset of $\Z$. In the case where $V\subset\Z$, we can consider an environment $\omega$ on $V$ to be the projection of an environment on $\Z$ by assigning arbitrary transition probability vectors at states $x\notin V$. Likewise, a measure $P$ on $\Omega_V$ may be thought of as, for example, a measure on $\Omega_{\Z}$ with $\omega(x,y)=\one_{\{y=x\}}$ for $x\notin V$, $P$--a.s. Though we do not emphasize it in the body of this paper, all environments considered are assumed to live in the space $\Omega_{\Z}$, and all measures on environment spaces are measures on $\Omega_{\Z}$. This allows us to intelligibly compare probabilities of the same event under different annealed measures, even if one of the measures is concentrated on walks that remain in a finite part of $\Z$. 

As another notational convenience, we will use interval notation to denote sets of consecutive integers in the state space $\Z$, rather than subsets of $\R$. Thus, for example, we will use $[1,\infty)$ to denote the set of integers to the right of 0. However, we make one exception, using $[0,1]$ to denote the set of all real numbers from 0 to 1.  

In the case where $V=\Z$, for each $x\in\Z$ we let $\omega^x\in\mathcal{M}_1(\Z)$ be the measure on $\Z$ given by $\omega^x(y)=\omega(x,x+y)$. Thus, $\omega=(\omega^x)_{x\in\Z}$ is the generic element of $\Omega$. For a subset $S\subseteq\Z$, let $\omega^S=(\omega^x)_{x\in S}$. In the case where $S$ is a half-infinite interval, we simplify our notation by using $\omega^{\leq x}$ to denote $\omega^{(-\infty,x]}$, and similarly with $\omega^{<x}$, $\omega^{\geq x}$, and $\omega^{>x}$. 
A common assumption, which our model will satisfy, is that the $\{\omega^x\}_{x\in\Z}$ are i.i.d. under $P$. 

\newpage
{\bf Random Walks in Dirichlet Environments}

Let $\Hg=(V,E,w)$ be a weighted directed graph with vertex set $V$, edge set\footnotemark
$~E\subseteq V\times V$, and a weight function $w:E\to\R^{>0}$. If $e=(x,y)\in E$, we say that $e$ is an edge from $x$ to $y$, and we say the {\em head} of $e$ is $\overline{e}=y$ and the {\em tail} of $e$ is $\underline{e}=x$. We say a set $S\subset V$ is {\em strongly connected} if for all $x,y\in S$, there is a path from $x$ to $y$ in $\Hg$ using only vertices in $S$.
To the weighted directed graph $\Hg$, we can associate the Dirichlet measure $P_{\Hg}$ on $(\Omega_{V},\mathcal{F}_{V})$, which we now describe. 

\footnotetext{We define weighted directed graphs in a way that precludes multiple edges from sharing the same head and tail. However, we could expand our definition to include weighted directed multigraphs, and natural generalizations of the results that are true for graphs as we define them would still hold. Describing these generalizations would cause some notational inconvenience that is unnecessary for our purposes, such as defining random walks that keep track of edges taken as well as vertices visited. Nevertheless, we use multigraphs in our illustrations as a visual aid. Multiple edges from one vertex to another in our illustrations can be interpreted as a single edge whose weight is the sum of the weights of the edges depicted. By the amalgamation property reviewed in Section \ref{sec:Background}, identifying or splitting these edges that share the same head and tail does not affect the distributions of transition probability vectors between sites.}

Recall the definition of the Dirichlet distribution: for a finite set $I$, take parameters ${\bf a}=(a_i)_{i\in I}$, with $a_i>0$ for all $i$. The Dirichlet distribution with these parameters is a probability distribution on the simplex $\Delta_I:=\{(x_i)_{i\in I}:\sum_{i\in I}x_i=1\}$ with density
\begin{equation*}D\left((x_i)_{i\in I}\right)=C({\bf a})\prod_{i\in I}x_i^{a_i-1},\end{equation*}
where $C({\bf a})$ is a normalizing constant. 

Define $P_{\Hg}$ to be the measure on $\Omega_V$ under which transition probabilities at the various vertices $x\in V$ are independent, and for each vertex $x\in V$, $(\omega(x,\overline{e}))_{\underline{e}=x}$ is distributed according to a Dirichlet distribution with parameters $(w(e))_{\underline{e}=x}$. (Or, if $V\subset\Z$, let $P_{\Hg}$ be any measure on $\Omega_{\Z}$ satisfying this description.\footnote{The measure $P_{\Hg}$ is then technically not unique, but its marginal on the set $\{\omega^V:\omega\in\Omega_{\Z}\}$ is unique, and that is all that will matter for our purposes.}) With $P_{\Hg}$-probability 1, $\omega(x,y)>0$ if and only if $(x,y)\in E$ for all $x,y\in V$. We will call a random environment chosen according to $P_{\Hg}$ a {\em Dirichlet environment on $\Hg$}.  We will use $E_{\Hg}$ to denote the associated expectation, and $\Pa_{\Hg}^x$ and $\Ea_{\Hg}^x$ to denote the annealed measure and expectations.

\begin{figure}
    \centering
    \includegraphics[width=6.5in]{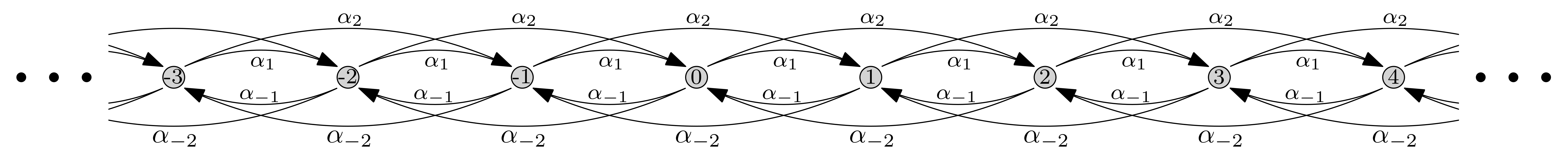}
    \caption{A portion of the graph $\G$ with $L=R=2$}
    \label{fig:graphG}
\end{figure}

We let $\G$ be the graph with vertex set $\Z$, edge set $\{(x,y)\in \Z\times\Z:-L\leq y-x\leq R,\alpha_{y-x}>0\}$, and weight function $(x,y)\mapsto\alpha_{y-x}$. An example with $L=R=2$ is represented in Figure \ref{fig:graphG} (here, and in other illustrations of graphs, we depict the case where $\alpha_0=0$, but our model does allow for $\alpha_0>0$). 

Our main concern in this paper is to characterize almost-sure recurrence, directional transience, and ballisticity of a random walk in a Dirichlet environment on $\G$ started at 0 in terms of the $\alpha_i$. 

Our assumption that $\alpha_L,\alpha_R>0$ is without loss of generality, since $L$ can be defined as the minimum $i$ with $\alpha_i>0$, and $R$ as the maximum of these $i$. Our assumption that $L$ and $R$ are positive means it is possible for the walk to move in both directions. Otherwise, directional transience and ballisticity are trivial.\footnote{Indeed, the walk has $|v|\geq1$ if $\alpha_0=0$, and if $\alpha_0>0$, one can easily check that $v>0$ if and only if $\sum_{i\neq0}\alpha_i>1$, which is what is needed to guarantee finite expected time to exit a self-loop.} The assumption that the $\gcd$ of all the $i$ such that $\alpha_i>0$ is 1 is needed to ensure that the Markov chain under $P_{\omega}^0$ is almost surely irreducible. It is without loss of generality because dividing all the $i$ by their $\gcd$ $g$ yields a graph $\G'$ where the law of $(\omega(gx,gy))_{x,y\in\Z}$ under $P_{\G'}$ is the same as the law of $(\omega(x,y))_{x,y\in\Z}$ under $P_{\G}$, and thus the law of $(gX_n)_{n=0}^{\infty}$ under $\Pa_{\G'}$ is the law of $(X_n)_{n=0}^{\infty}$ under $\Pa_{\G}$. From this assumption, it follows that there is an $m$ large enough that every interval of length $m$ is strongly connected in $\G$. Let $m_0$ be such an integer, chosen large enough that also $m_0\geq\max(L,R)$. We will use this $m_0$ in several proofs throughout this paper.

\subsection{Results}\label{subsec:results}
It has been shown \cite{Tournier2015}
that for nearest-neighbor random walks in Dirichlet environments in $\Z^d$, $d\geq1$, if the annealed expected position of the first step is nonzero, then its direction is the direction of almost-sure transience. A remark in that paper points out that the proof also works for the bounded-jump case. 

Let $d^{\totheright}=\sum_{i=1}^Ri\alpha_i$ and $d^{\totheleft}=\sum_{i=-L}^{-1}|i|\alpha_i$. Then $\kappa_1:=d^{\totheright}-d^{\totheleft}$ is the weighted sum of the weights $\alpha_i$, and its sign is the sign of $\Ea_{\G}^0[X_1]$. In fact, let $c^{\totheright}$ and $c^{\totheleft}$ be the unweighted sums $\sum_{i=1}^{R}\alpha_i$ and $\sum_{i=-L}^{-1}\alpha_i$, respectively. Then one can check using the well known Lemma \ref{DERRW}, which we recall below, or simply using the expectation of a beta random variable, that $\Ea_{\G}^0[X_1]=\frac{\kappa_1}{c^{\totheleft}+\alpha_0+c^{\totheright}}$. 

Applying the results from \cite{Tournier2015} with $d=1$ yields a characterization of directional transience for our model in the case where $\kappa_1\neq0$. To complete the characterization, we prove in Section \ref{sec:RecurrenceAndTransience} that the walk is recurrent in the case $\kappa_1=0$. We also outline a version of the argument from \cite{Tournier2015}, because it is a bit simpler in one dimension and because similar ideas are used in the proof for the case $\kappa_1=0$, and to a lesser extent for later proofs in this paper.

\begin{thm}\label{thm:transiencecriterion}
The directional transience of the walk is determined by the sign of $\kappa_1$ in the following way:

\begin{enumerate}
\item If $\kappa_1>0$, then $\lim_{n\to\infty}X_n=\infty,~\Pa_{\G}^0\text{--a.s.}$
\item If $\kappa_1<0$, then $\lim_{n\to\infty}X_n=-\infty,~\Pa_{\G}^0\text{--a.s.}$
\item If $\kappa_1=0$, then $-\infty=\liminf_{n\to\infty}X_n<\limsup_{n\to\infty}X_n=\infty,~\Pa_{\G}^0\text{--a.s.}$
\end{enumerate}
\end{thm}

Our goal in Sections \ref{sec:abstract} and \ref{sec:ballisticity} is to investigate the limit $v:=\lim_{n\to\infty}\frac{X_n}{n}$, showing that it exists and characterizing when it is positive. 
Section \ref{sec:abstract} provides characterizations of ballisticity for a more general setting than the Dirichlet model. We consider the following conditions for a probability measure $P$ on $\Omega_{\Z}$:
\begin{enumerate}[(C1)]
    \item\hypertarget{cond:C1} The $\{\omega^x\}_{x\in\Z}$ are i.i.d. under $P$.
    \item\hypertarget{cond:C2} For $P$--a.e. environment $\omega$, the Markov chain induced by $P_{\omega}^0$ is irreducible.
    \item\hypertarget{cond:C3} For $P$--a.e. environment $\omega$, $\omega(a,b)=0$ whenever $b$ is outside $[a-L,a+R]$. 
\end{enumerate}
It was shown in \cite{Key1984} that under these assumptions, a 0-1 law holds for directional transience. That is, the walk is either almost surely transient to the right, almost surely transient to the left, or almost surely recurrent. We show that under these assumptions, the limit $v$ necessarily exists. In the recurrent case, it necessarily holds that $v=0$. We then provide two abstract characterizations of ballisticity under the following additional assumption.
\begin{enumerate}[(C4)]
    \item\hypertarget{cond:C4} For $P$--a.e. environment $\omega$, $\lim_{n\to\infty}X_n=\infty$, $P_{\omega}^0$--a.s.
\end{enumerate}
By symmetry, our characterizations also handle the case where the walk is transient to the left, and thus by the 0-1 law for directional transience, completely characterize the regime $v\neq0$ for all measures $P$ satisfying (\hyperlink{cond:C1}{C1}), (\hyperlink{cond:C2}{C2}), and (\hyperlink{cond:C3}{C3}).

The first characterization strengthens one given by Brémont, who showed (see \cite[Theorem 3.7]{Bremont2002}, \cite[Proposition 9.1]{Bremont2004}) that for a walk that is transient to the right, $v>0$ if and only if the annealed expected time to reach $[1,\infty)$ is finite. Brémont's works used an ellipticity assumption that is too strong to apply to our Dirichlet environments. We therefore prove the lemma without the assumption.

To formally state it, we must establish notation for hitting times. For a given walk ${\bf X}=(X_n)_{n=0}^{\infty}$, we define $H_x({\bf X})$ to be the first time the walk hits $x\in\Z$. That is,
\begin{equation*}
    H_x({\bf X})=\inf\{n\in\N_0:X_n=x\}
\end{equation*}
We usually write it as $H_x$ when we can do so without ambiguity. For a subset $S\subset\Z$, let $H_S=\min_{x\in S}H_x$. First positive hitting times are denoted as $\tilde{H}_x$ or $\tilde{H}_S$. That is,
\begin{equation*}
    \tilde{H}_x({\bf X})=\inf\{n\in\N:X_n=x\},
\end{equation*}
and $\tilde{H}_S=\min_{x\in S}\tilde{H}_x$.
If the set is the half-infinite interval $[x,\infty)$, we use $H_{\geq x}$ to denote its hitting time, and similarly with $H_{>x}$, $H_{\leq x}$, and $H_{>x}$. 

\begin{lem}\label{lem:BremontLemma}
Let $P$ be a probability measure on $\Omega_{\Z}$ satisfying (\hyperlink{cond:C1}{C1}), (\hyperlink{cond:C2}{C2}), (\hyperlink{cond:C3}{C3}), and (\hyperlink{cond:C4}{C4}). Then $v>0$ if and only if $\Ea^0[H_{\geq1}]<\infty$, where $H_{\geq 1}$ is the first time the walk hits $[1,\infty)$.
\end{lem}

This characterization is quite natural, given that in the nearest-neighbor case we in fact have the identity $v=1/\Ea^0[H_{\geq1}]$, where the fraction is understood to be 0 if the denominator is infinite. 
However, although it is natural, we do not know a way to check it directly in the $(L,R)$ case, even for Dirichlet environments. We therefore use Lemma \ref{lem:BremontLemma} as well as a construction of a ``walk from $-\infty$ to $\infty$", to provide another criterion for ballisticity, which is based on the expected number of visits to a specific site.

For a walk ${\bf X} = (X_n)_{n=0}^{\infty}$ on a vertex set $V$ with $x\in V$, $N_x({\bf X})=\#\{n\in\N_0:X_n=x\}$ is the number of times the walk is at site $x$. We usually write it as $N_x$ if we are able to do so without ambiguity. For a subset $S\subset\Z$, let $N_S=\sum_{x\in S}N_x$. We prove the following lemma in Section \ref{sec:abstract}.

\begin{lem}\label{lem:ballistic}
Let $P$ be a probability measure on $\Omega_{\Z}$ satisfying (\hyperlink{cond:C1}{C1}), (\hyperlink{cond:C2}{C2}), (\hyperlink{cond:C3}{C3}), and (\hyperlink{cond:C4}{C4}). Then $v>0$ if and only if $\Ea^{0}[N_{0}]=E[E_{\omega}^0[N_0]]<\infty$.
\end{lem}

Thus, the question of ballisticity is reduced to the integrability of the ``Green function" $E_{\omega}^0[N_0]$ under the measure on environments. We devote Section \ref{sec:ballisticity} to answering this question in the case of our Dirichlet measure $P_{\G}$. In fact, we go further and characterize integribility of $E_{\omega}^0[N_0]^s$ for any $s>0$. This is done in terms of two parameters; one we call $\kappa_0$, and the other is $\kappa_1=d^{\totheright}-d^{\totheleft}$. Although we always have $\kappa_0>\kappa_1$ in the case $L=R=1$, the ordered pair $(\kappa_0,\kappa_1)$ can take on any value in the first quadrant of $\R^2$ in the general case (see Proposition \ref{prop:orderedpair}).

For a weighted, directed graph $\Hg=(V,E,w)$ and vertex $x\in V$, one may define $\kappa_0=\kappa_0(\Hg,x)$ as the minimal total weight of edges exiting a finite, strongly connected set of vertices containing $x$ in $\Hg$. We give a precise definition for our graph $\G$ in Section \ref{subsec:FiniteTraps}; see \eqref{eqn:kappa0def} (the definition is not vertex-dependent due to the translation-invariance of $\G$).
The smaller $\kappa_0(\Hg,x)$ is, the greater is the propensity of a walk drawn according to $\Pa_{\Hg}^x$ to get stuck for a long time in a finite trap containing $x$ \cite{Tournier2009}.
This parameter is the ``$\kappa$" defined for $\Z^{d}$ in \cite{Bouchet2013}. In that model, the walks are nearest-neighbor, and so the underlying directed graph has an edge from $x$ to $y$ precisely when $x$ and $y$ are adjacent. There, it can easily be shown that the worst traps are just pairs of vertices, and so $\kappa_0$ has an explicit formula as a minimum of $d$ different sums of edge weights. By contrast, our model encompasses many underlying directed graphs (even before assignment of weights). For each underlying directed graph there is a different formula for $\kappa_0$ as a minimum of finitely many sums, but we do not have a general method to find the formula given a particular underlying directed graph. This is because we have no simple general way to know what the worst finite traps look like. However, we find the formula in several examples in Appendix \ref{append:kappa}, and show in Proposition \ref{prop:GammaAttained} that $\kappa_0$ can be calculated directly from $L$, $R$, and the specific values of the $\alpha_i$, even without a general formula in terms of the $\alpha_i$.

This parameter $\kappa_0$ plays an important role in the integrability of $E_{\omega}^0[N_0]$. 
For a set $S\subseteq\Z$, and for $x \in S$, define $N_x^S$ to be the amount of time a walk spends at $x$ before leaving $S$ for the first time (we always have $N_x^S\leq N_x$). We define $\kappa_0$ as an infimum of sums of edge weights. This infimum is over an infinite set, but once we can show that it is actually a minimum, the following theorem follows almost immediately from \cite[Theorem 1]{Tournier2009}.

\begin{thm}\label{thm:gamma}
For $s>0$, the following are equivalent:
\begin{enumerate}[(a)]
    \item $\kappa_0\leq s$.
    \item For all sufficiently large $M$, $E_{\G}\left[E_{\omega}^0\left[N_0^{[-M,0]}\right]^s\right]=\infty$.
    \item For some $M\geq0$, $E_{\G}\left[E_{\omega}^0\left[N_0^{[-M,M]}\right]^s\right]=\infty$.
\end{enumerate}
\end{thm}

Letting $s=1$, the implication (a)$\Rightarrow$(b) or (a)$\Rightarrow$(c) shows that if $\kappa_0\leq1$, then $\Ea_{\G}^0[N_0]=\infty$, which by Lemma \ref{lem:ballistic} implies $v=0$.
We include condition (b) because the implication (a)$\Rightarrow$(b) allows one to arrive at the same conclusion using Lemma \ref{lem:BremontLemma} (by showing that $\kappa_0\leq1$ implies $\Ea_{\G}^0[H_{\geq1}]=\infty$).

While $\kappa_0$ controls the moments of the quenched expected amount of time the walk spends at 0 before exiting a finite region of the graph, the parameter $\kappa_1=d^{\totheright}-d^{\totheleft}$ controls, in the same way, the moments of the quenched expected number of times the walk traverses arbitrarily large regions of the graph. For $x<y\in\Z$, we define the following functions of a walk ${\bf X}$:

\begin{itemize}
    \item $N_{x,y}({\bf X})=\#\big\{n\in\N_0:X_n=x,\sup\{j<n:X_j=y\}>\sup\{j<n:X_j=x\}\big\}$ is the number of times the walk hits $x$ after more recently having hit $y$, or the number of ``trips from $y$ to $x$".
    \item $N_{x,y}'({\bf X}):=\#\big\{n\in\N_0:X_n\leq x,\sup\{j<n:X_j\geq y\}>\sup\{j<n:X_j \leq x\}\big\}$ is the number of trips leftward across $[x,y]$.
\end{itemize} 
Again, we write these as $N_{x,y}$ and $N_{x,y}'$ if we can do so without ambiguity. Note that $N_x\geq N_{x,y}$, and also $N_{x,y}'\geq N_{x,y}$.
We prove the following theorem.

\begin{thm}\label{thm:TFAE2}
Let $\kappa_1>0$, so that the walk is transient to the right. Then, if $s>0$, the following are equivalent:
\begin{enumerate}[(a)]
    \item $\kappa_1>s$.
    \item There is an $M\geq0$ such that for all $x,y\in\Z$ with $y-x\geq M$, $E_{\G}\left[E_{\omega}^0[N_{x,y}']^s\right]<\infty$.
    \item There exist $x<y\in\Z$ such that $E_{\G}\left[E_{\omega}^0[N_{x,y}]^s\right]<\infty$.
\end{enumerate}
\end{thm}
The proof of Theorem \ref{thm:TFAE2} is long and naturally divides into two parts, so we prove the parts separately as Proposition \ref{prop:InfiniteTraps} and Proposition \ref{prop:FewOscillations}. 
Letting $s=1$, the contrapositive of the implication (c)$\Rightarrow$(a) tells us that if $\kappa_1\leq1$, then $\Ea_{\G}^0[N_0]=\infty$, which by Lemma \ref{lem:ballistic} implies $v=0$.

Combining the theorems stated so far, we can see that if $\kappa_0\leq1$, then the walk is not ballistic due to finite trapping, and that if $\kappa_1\leq1$, then the walk is not ballistic due to large-scale backtracking. We would like to show that if both parameters are greater than 1, then the walk is ballistic. For every environment $\omega$ on $\Z$ and every $M>0$, we have
\begin{equation}\label{eqn:447}
E_{\omega}^0[N_0]=
E_{\omega}^0\left[N_0^{[-M,M]}\right]
E_{\omega}^0\left[\#\left\{\parbox{1.75in}{Times exiting $[-M,M]$ and then returning to 0}\right\}\right].
\end{equation}
The first expectation on the right relates to finite trapping, and the second to large-scale backtracking. By Theorem \ref{thm:gamma}, the term $E_{\omega}^0\left[N_0^{[-M,M]}\right]$ has finite moments up to $\kappa_0$ under $P_{\G}$ for $M$ sufficiently large. 
And the number of times exiting $[-M,M]$ and returning to $0$ is between $N_{-M-1,0}+N_{0,M+1}$ and $N_{-M-1,0}'+N_{0,M+1}'$, so the term $E_{\omega}^0\left[\#\left\{\parbox{1.75in}{Times exiting $[-M,M]$ and then returning to 0}\right\}\right]$ has finite moments up to $\kappa_1$ by Theorem \ref{thm:TFAE2}. If the two terms on the right side of \eqref{eqn:447} were independent under $P_{\G}$, we could conclude that $E_{\G}[E_{\omega}^0[N_0]^s]<\infty$ if and only if $s<\min(\kappa_0,\kappa_1)$. However, they are not independent. We therefore ask whether it is possible that the phenomena of finite trapping and large-scale backtracking may conspire together to prevent ballisticity, even if neither is strong enough to do it on its own. This question is of interest for general RWRE on $\Z$ with bounded jumps. 

\begin{quest}\label{question}
Let $P$ be a measure satisfying (\hyperlink{cond:C1}{C1}), (\hyperlink{cond:C2}{C2}), (\hyperlink{cond:C3}{C3}), and (\hyperlink{cond:C4}{C4}), under which both terms on the right of \eqref{eqn:447} have finite expectation for all $M$; that is, $E\left[E_{\omega}^0\left[N_0^{[-M,M]}\right]\right]<\infty$ and $E\left[E_{\omega}^0\left[\#\left\{\parbox{1.75in}{Times exiting $[-M,M]$ and then returning to 0}\right\}\right]\right]<\infty$. Does it necessarily follow that $\Ea^0[N_0]=E[E_{\omega}^0[N_0]]<\infty$?
\end{quest}

We are able to answer this question in the affirmative for our Dirichlet model. In fact, we characterize all finite moments of $E_{\omega}^0[N_0]$ under $P_{\G}$.

\begin{thm}\label{thm:MainMomentTheorem}
Assume $\kappa_1>0$. Then $E_{\G}\left[\left(E_{\omega}^0[N_0]\right)^s\right]<\infty$ if and only if $s<\min(\kappa_0,\kappa_1)$.
\end{thm}

Combining this with Lemma \ref{lem:ballistic}, we get a complete characterization of ballisticity.

\begin{thm}\label{thm:MainBallisticTheorem}
Assume $\kappa_1>0$. Then the walk is ballistic if and only if $\min(\kappa_0,\kappa_1)>1$. 
\end{thm}

\subsection{Acknowledgments}

The author thanks his advisor, Jonathon Peterson, for suggesting the problem and for his mentorship.

\section{Background on Dirichlet Environments}\label{sec:Background}

For a more comprehensive overview of random walks in Dirichlet environments, their properties, known results, and techniques used to achieve those results, see \cite{Sabot&Tournier2016}. Here, we review specific results that will be useful to us. To begin with, we recall some important properties of Dirichlet distributions.

\begin{property}[Amalgamation]\label{amalgamation}
Assume $(U_i)_{i\in I}$ has Dirichlet distribution on $\Delta_I$ with parameters $(a_i)_{i\in I}$. Let $I_1,\ldots,I_r$ be a partition of $I$. The random vector $\left(\sum_{i\in I_k}U_i\right)_{1\leq k\leq r}$ on the simplex $\{(x_i)_{i=1}^r:\sum_{i=1}^rx_i=1\}$ follows the Dirichlet distribution with parameters $\left(\sum_{i\in I_k}a_i\right)_{1\leq k\leq r}$.
\end{property}

\begin{property}[Restriction]\label{restriction}
Assume $(U_i)_{i\in I}$ has Dirichlet distribution on $\Delta_I$ with parameters $(a_i)_{i\in I}$. Let $J$ be a nonempty subset of $I$. The random vector $\left(\frac{U_i}{\sum_{j\in J}U_j}\right)_{i\in J}$ , which takes values on the simplex $\Delta_J$, follows the Dirichlet distribution with parameters $\left(a_i\right)_{i\in J}$ and is independent of $\sum_{j\in J}U_j$. It is also independent of $\left(U_k\right)_{k\notin J}$.
\end{property}

By the amalgamation property, the marginal distribution of each coordinate of a Dirichlet random vector is a beta distribution. The next property bounds the probability that a beta random variable is small.

\begin{property}[Moments]\label{moments}
Assume $X$ is a beta random variable with parameters $(a,b)$. Then there exists constants $0<c<C$ such that for all $\tepsilon\in[0,1]$,
\begin{equation}\label{eqn:moments}
c\tepsilon^a \leq P(X<\tepsilon)\leq C\tepsilon^a.    
\end{equation}
In particular, $E\left[\frac1{X^s}\right]<\infty$ if and only if $s<a$. 
\end{property}

Dirichlet environments were first studied for their connection to a stochastic process called a {\em directed edge reinforced random walks} (DERRW), which we now define. For a weighted directed graph $\Hg = (V,E,w)$, and for an initial vertex $x_0\in V$, we define the stochastic processes $(X_n)_{n=0}^{\infty}$ and $\left(r(e,n)\right)_{e\in E,n\geq0}$ as follows: with probability 1, $X_0=x_0$ and $r(e,0)=w(e)$ for all $e$. If $X_n=x$ and $e_1$ is an edge with $\underline{e}=x$ (i.e., the edge is rooted at $x$), then the walk takes edge $e_1$, so that $X_{n+1}=\overline{e_1}$, with probability $\frac{r(e_1,n)}{\sum_{\underline{e}=x}r(e,n)}$. Each time an edge is taken, its weight $r$ is increased by 1; otherwise, weights do not change. That is, $r(e,n+1)=\begin{cases}r(e,n)+1~~\text{if }(X_n,X_{n+1})=e\\ r(e,n)~~~\text{otherwise} \end{cases}$. The process $(X_n)_{n=0}^{\infty}$ is the DERRW\label{page:619} on $\Hg$ started at $x_0$. The following lemma was first shown in \cite{Enriquez&Sabot2002} and \cite{Keane&Rolles2002}.

\begin{lem}\label{DERRW}
 Let $V$ be a set, and let $\Hg=(V,E,w)$ be a weighted directed graph with vertex set $V$. Then the law of a DERRW on $\Hg=(V,E,w)$ started at vertex $x$ is the annealed law $\Pa_{\Hg}^x$ of the RWDE on $\Hg$.
\end{lem}

We now describe an important time-reversal lemma for Dirichlet environments. Let $\Hg=(V,E,w)$ be a weighted directed graph, and let $P_{\Hg}$ be the associated product Dirichlet measure on the set $\Omega_V$ of environments on $V$.

For a vertex $x\in V$, the {\em divergence} of $x$ in $\Hg$ is $\text{div}(x)=\sum_{\overline{e}=x}w(e)-\sum_{\underline{e}=x}w(e)$. If the divergence is zero for all $x$, we say the graph $\Hg$ has zero divergence.

Now if $\Hg$ is a finite, strongly connected graph, then for $P_{\Hg}$--a.e. $\omega\in\Omega_V$, there exists an invariant probability for the Markov chain corresponding to $\omega$. Call this invariant measure $\pi^{\omega}$. Define the {\em time reversed} environment $\check{\omega}$ by
\begin{equation*}\check{\omega}(x,y):=\frac{\pi^{\omega}(y)}{\pi^{\omega}(x)}\omega(y,x)\end{equation*}
One can check that $\check{\omega}$ is an 
\ifSHOWEXTRA
\else
environment and that the probability, under $\omega$, of taking any loop is equal to the probability, under $\check{\omega}$, of taking the reversed loop; that is, for a path $(v_0,v_1,\ldots,v_n=v_0)$, we have
\begin{equation*}
    P_{\omega}^{v_0}(X_1=v_1,X_2=v_2,\ldots,X_n=v_n)
    =
    P_{\check{\omega}}^{v_0}(X_1=v_{n-1},X_2=v_{n-2},\ldots,X_n=v_0).
\end{equation*}
See \cite{Sabot&Tournier2011} for details.
\fi
The following lemma can be proven using Lemma \ref{DERRW}. It was first proven analytically in \cite{Sabot2011}, and its probabilistic proof was first given in \cite{Sabot&Tournier2011}. Let $\check{\Hg}$ be the graph made by reversing all edges of $\Hg$ and keeping the same weights.\footnote{Formally, if $e=(a,b)$, then let $\check{e}=(b,a)$. Then define $\check{w}(\check{e})=w(e)$ and $\check{E}=\{\check{e}:e\in E\}$. Now let $\check{\Hg}=(V,\check{E},\check{w})$.}, and let $P_{\check{\Hg}}$ be the associated measure on $\Omega_V$. 

\begin{lem}[\cite{Sabot&Tournier2016}, Lemma 3.1]\label{timereversal}
If the graph $\Hg$ has zero divergence, then the law of $\check{\omega}$ is $P_{\check{\Hg}}$.
\end{lem}

In other words, drawing an environment $\omega$ according to $P_{\Hg}$ and then time-reversing it is the same as reversing the edges of $\Hg$ to get $\check{\Hg}$ and then drawing an environment according to $P_{\check{\Hg}}$. This lemma implies that the probability, under $\Pa_{\Hg}$, of taking any loop is equal to the probability, under $\Pa_{\check{\Hg}}$, of taking the reversed loop. Indeed, for our purposes, the use of Lemma \ref{timereversal} comes from the following corollary.

\begin{cor}\label{cor:loopreversal}
Let $\Hg$ be as described above, and let $x,y\in V$ such that there is an edge $e$ from $y$ to $x$ in $\Hg$. Then, letting $\tilde{H}_x$ denote the first positive hitting time of $x$,
\begin{enumerate}
    \item The law of $P_{\omega}^x(X_{\tilde{H}_x-1}=y)$ under $P_{\Hg}$ is the law of $P_{\omega}^x(X_1=y)=\omega(x,y)$ under $P_{\check{\Hg}}$.
    \item $\Pa_{\Hg}^x(X_{\tilde{H_x}-1}=y)=\Pa_{\check{\Hg}}^x(X_1=y)=\frac{w(y,x)}{\sum_{v\in V}w(v,x)}$.
\end{enumerate}
\end{cor}
The formula for the probability as a fraction comes from either Lemma \ref{DERRW} or the fact that the expectation of a beta random variable with parameters $(a,b)$ is $\frac{a}{a+b}$.

The next lemma we recall was proven by Tournier \cite{Tournier2009}. We will refer to it as Tournier's lemma. For a set $S\subseteq V$, define
\begin{equation}\label{betaS}
    \beta_S:=\sum_{\underline{e}\in S,~\overline{e}\notin S}w(e).
\end{equation}
This parameter $\beta_S$ is the sum of the weights of all edges exiting the set $S$.

\begin{lem}[see \cite{Tournier2009}, Theorems 1 and 2]\label{lem:TourniersLemma}
Let $\Hg=(V\cup\{\partial\},E,w)$ be a finite weighted directed graph with $\partial$ a unique sink reachable from every other site. We denote by $P_{\Hg}$ the corresponding Dirichlet distribution on environments.

For every $s>0$, the following statements are equivalent:
\begin{enumerate}
    \item $E_{\Hg}[E_{\omega}^x[N_x]^s]<\infty$.
    \item For every strongly connected subset $S$ of $V$ with $x\in S$, $\beta_S>s$.
\end{enumerate}
\end{lem}

In particular, by letting $s=1$, we see that $\Ea_{\Hg}^x[N_x]<\infty$ if and only if for every strongly connected subset $S$ of $V$ containing $x$, $\beta_S>1$. The formulation given in Theorem 1 of \cite{Tournier2009} is in terms of strongly connected sets of edges rather than vertices, but implies ours. 
\ifSHOWEXTRA
\else{Tournier defines a set $A$ of edges to be strongly connected if all heads and tails of edges in $A$ can communicate only through $A$, and for such a set defines $\beta_A$ to be the sum of weights of edges that are not in $A$ but share tails with edges in $A$. }
\fi
Every $\beta_S$ for a set $S$ of vertices is $\beta_A$ for the set $A$ of edges between vertices in $S$. On the other hand, for any strongly connected set $A$ of edges, one can take $S$ to be the set of heads or tails of edges in $A$ and take $A'$ to be the set of edges between vertices in $S$ (so that $A\subset A'$). Then $\beta_S=\beta_{A'}<\beta_A$. From this, one can check that Tournier's formulation implies ours.

\section{Recurrence and Transience}\label{sec:RecurrenceAndTransience}

In this section, we are to prove the characterization of recurrence and directional transience given in Theorem \ref{thm:transiencecriterion}.

\setcounter{ExampleSave}{\value{exmp}}
\setcounter{SectionSave}{\value{section}}
\setcounter{section}{1}
\setcounter{exmp}{1}
\setcounter{claim}{0}

\begin{proof}[Proof of Theorem \ref{thm:transiencecriterion}]

{\em Case 1:} $\kappa_1>0$. 
\hypertarget{proof:transiencecriterion}
This follows from \cite[Corollary 1]{Tournier2015}, but we outline the proof for the one-dimensional case. We are to show that $\lim_{n\to\infty}\frac{X_n}{n}=\infty$, $\Pa_{\G}^0$--a.s. The steps are as follows:

\begin{enumerate}[(a)]
    \item Show that $\Pa_{\G}^0(H_{\geq M}<\tilde{H}_{\leq0})$ is bounded away from 0 as $M$ approaches $\infty$. (This is shown in \cite[Theorem 1]{Tournier2015}). 
    \item Taking limits, conclude that $\Pa_{\G}^0(\tilde{H}_{\leq0}=\infty)>0$, implying the walk is transient to the right with positive probability.
    \item Use the 0-1 law \cite[Theorem 11]{Key1984} to conclude that if the probability of transience to the right is positive, it is 1. 
\end{enumerate}

We outline the proof of (a) for convenient reference. For $M>R+L$, we consider a weighted directed graph $\G_M$ with vertex set $[0,M]$. Each vertex in $[1,M-1]$ has the same edges with the same weights as those on $\G$, except that edges that would terminate at points less than zero are simply edges to the point $0$, and edges that would terminate at points greater than $M$ are simply edges to the point $M$. If this would result in multi-edges, each multi-edge is replaced with a single edge whose weight is the sum of the weights in the multi-edge; however, we leave the multi-edges in our illustrations in order to show more clearly where this occurs. Based on the edges and weights we've described so far, zero divergence already holds at points from $R+1$ to $M-L-1$, but points to the left of $R$ and to the right of $M-L$ are ``missing" incoming weights from vertices to the left of $0$ and to the right of $M$, respectively. Therefore, to each vertex $1\leq j\leq R$, we add an edge from 0 with weight $\sum_{i=j}^R\alpha_i$ (depicted in Figure \ref{fig:graphGM} as a multi-edge) in order to achieve zero divergence at $j$. Likewise, to each edge $M-L\leq j\leq M$, we add an edge from $M$ with weight $\sum_{i=M-L}^{M-j}\alpha_{i-M}$. Based on these weights, the site $0$ has incoming weight $d^{\totheleft}$ and outgoing weight $d^{\totheright}$, and the site $M$ has incoming weight $d^{\totheright}$ and outgoing weight $d^{\totheleft}$. To adjust for this, we add a special edge from $M$ to $0$ with weight $d^{\totheright}-d^{\totheleft}=\kappa_1$.

\begin{figure}
    \centering
    \includegraphics[width=6.5in]{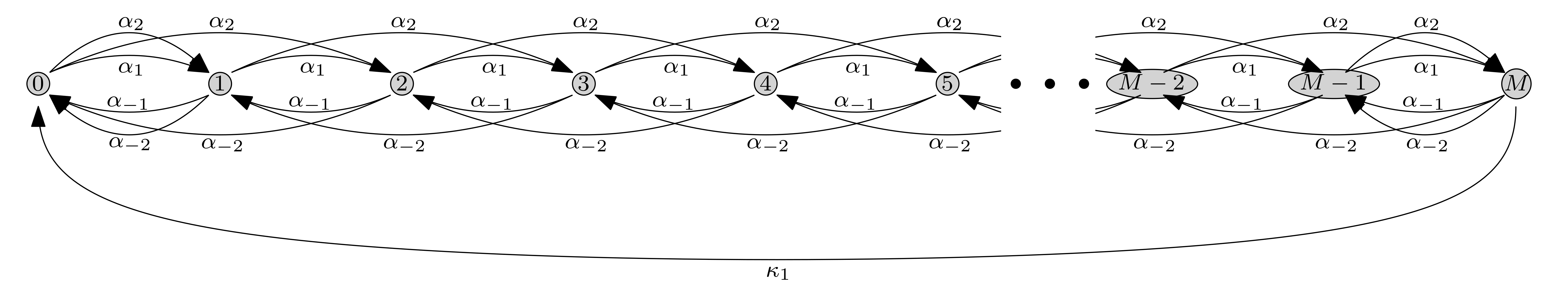}
    \caption{The graph $\G_M$}
    \label{fig:graphGM}
\end{figure}

Now the graph satisfies the zero-divergence property. This is the graph $\G_M$, pictured in Figure \ref{fig:graphGM}. 
For now, the usefulness of $\G_M$ comes from the following claim, which we prove in its entirety because the ideas in its proof will be referenced several times throughout this paper.
\begin{claim}\label{claim:ComparingGraphs}
\begin{equation}\label{eqn:613}
\Pa_{\G}^x(H_{\geq M}<H_{\leq0}) = \Pa_{\G_M}^x(H_M<H_0),\quad \quad 0<x<M.
\end{equation}
\end{claim}

To prove this, consider for each environment $\omega$ on $\Z$ a modified environment $\omega'$, where transition probabilities between sites in $[1,M-1]$ are the same as in $\omega$, but for each $i\in[M-R,M-1]$, $\omega'(i,M)=\sum_{j\geq M}\omega(i,j)$, and for each $i\in[1,R]$, $\omega'(i,0)=\sum_{j\leq0}\omega(i,j)$. Then by construction, a walk drawn according to $P_{\omega'}^x$ for any $x$ strictly between $0$ and $M$ and stopped when it hits $0$ or $M$ follows the same law (except possibly for the terminating site) as the law of a walk drawn according to $P_{\omega}^x$ and stopped when it reaches $(-\infty,0]$ or $[M,\infty)$. In particular,
\begin{equation*}
    P_{\omega}^x(H_{\geq M}<H_{\leq0}) = P_{\omega'}^x(H_M<H_0),\quad \quad 0<x<M.
\end{equation*}
On the other hand, by the amalgamation property of Dirichlet random vectors, we also see that for every $y\in[1,M-1]$, the law of $(\omega')^x$ under $P_{\G}$ is a Dirichlet distribution, and in fact is the same as the law of $\omega^x$ under $P_{\G_M}$. Hence, for each $0<x<M$, we have
\begin{align*}
    \Pa_{\G}^x(H_{\geq M}<H_{\leq0}) &=
    E_{\G}\left[P_{\omega}^x(H_{\geq M}<H_{\leq0})\right]
    \\
    &=E_{\G}\left[P_{\omega'}^x(H_M<H_0)\right]
    \\
    &=E_{\G_M}\left[P_{\omega}^x(H_M<H_0)\right]
    \\
    &=\Pa_{\G_M}^x(H_M<H_0).
\end{align*}
This proves the claim.

From Corollary \ref{cor:loopreversal} (2), we can get $\Pa_{\G_M}^0(H_M<\tilde{H}_0) \geq\frac{\kappa_1}{d^{\totheright}}$. We can use Claim \ref{claim:ComparingGraphs} to show that $\Pa_{\G_M}^0(H_M<\tilde{H}_0)\leq \Pa_{\G}^R(H_{\geq M}<H_{\leq 0})$.
Putting the two together, we get 
\begin{equation*}
\Pa_{\G}^R(H_{\geq M}<H_{\leq0})\geq \frac{\kappa_1}{d^{\totheright}}
\end{equation*} for all $M$. By independence of sites, we then have 
\begin{align*}
\Pa_{\G}^0(H_{\geq M}<\tilde{H}_{\leq0}) &\geq \Pa_{\G}(X_1=R)\Pa_{\G}^R(H_{\geq M}<H_{\leq0})
\\
&\geq \Pa_{\G}(X_1=R)\frac{\kappa_1}{d^{\totheright}}.
\end{align*}
This is the bound we needed to finish case 1.

{\em Case 2:} $\kappa_1<0$. This follows from Case 1 by symmetry.

{\em Case 3:} $\kappa_1=0$.

Again, let $M>R+L$. We will define a graph $\Hg_M$ with vertex set $[-M,M]$. The edges in $\Hg_M$ are defined in much the same way as those in $\G_M$: for vertices other than endpoints, edges are the same, except that edges with heads to the left of $-M$ and the right of $M$ become edges to $-M$ and $M$, respectively (again, multi-edges can be amalgamated into single edges to fit our definitions, but we leave multi-edges in our illustration for clarity). 
As in the case of $\G_M$, we add edges from the left and right endpoints to vertices nearby as needed to achieve zero divergence on those vertices, but since $d^{\totheright}=d^{\totheleft}$, no edge from $M$ to $-M$ is then needed to achieve zero divergence at the endpoints.
Finally, we add edges of weight 1 from 0 to $-M$ and from $-M$ to $0$, which maintains the zero-divergence condition (the weight 1 is arbitrary; any positive weight would work). An example of the graph $\Hg_M$ is depicted in Figure \ref{fig:graphHM}.

\begin{figure}
    \centering
    \includegraphics[width=6.5in]{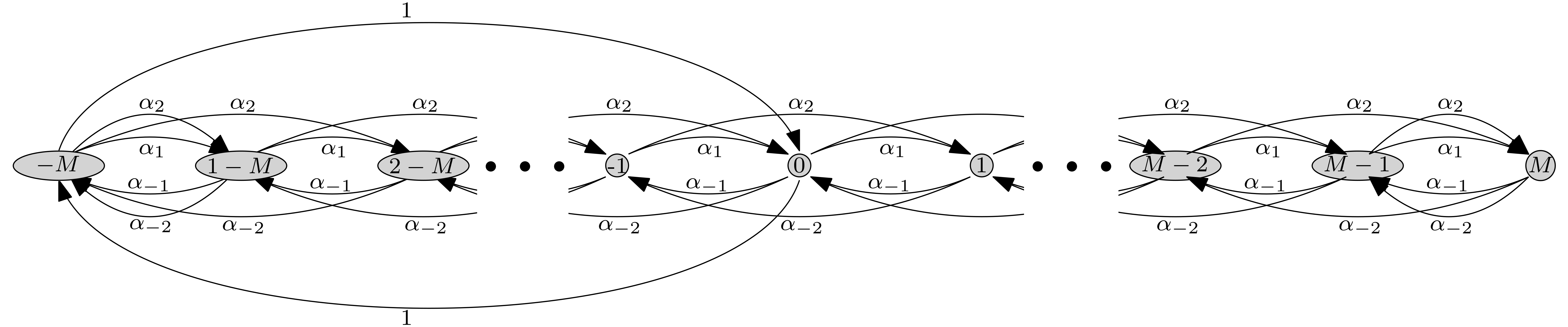}
    \caption{The graph $\Hg_M$, with $L=R=2$}
    \label{fig:graphHM}
\end{figure}

Assume for a contradiction that the walk on $\Z$ is transient, $\Pa_{\G}$--a.s. By the 0-1 law of \cite{Key1984} it is directionally transient; and without loss of generality assume it is transient to the right. We define $T_{-M\to0}=\inf\{n\in\N:X_n=0, X_{n-1}=-M\}$. Note that $\{\tilde{H}_0=T_{-M\to0}\}=\{X_{\tilde{H}_0-1}=-M\}$ is the event that the first return to zero is by the special edge. We make the following claim.

\begin{claim}\label{claim:recurrenceclaim}
There exists $c<1$ such that for all large enough $M$, $\Pa_{\Hg_M}^{-M}(H_0=T_{-M\to0})<c$. That is, the probability, starting from $-M$, that the first visit to 0 is by the special edge from $M$ to $0$, is bounded away from 1.
\end{claim}

To prove the claim, it suffices to show that $\Pa_{\Hg_M}^{-M}(H_0\neq T_{-M\to0})$ is bounded away from 0. Recall that $m_0\geq\max(L,R)$ is large enough that every interval of length $m_0$ is strongly connected in $\G$. Let $M\geq m$. We then have
\begin{equation}\label{eqn:821}
\min_{0\leq i \leq R-1}\Pa_{\Hg_M}^i(H_0<H_{[0,m_0-1]^c})=b>0,
\end{equation} and $b$ does not depend on $M$. 
Now by independence of sites and the strong Markov property,
\begin{align}
\notag
\Pa_{\Hg_M}^{-M}(\tilde{H}_0\neq T_{-M\to0})
&\geq\Pa_{\Hg_M}^{-M}(X_1=R-M)\Pa_{\Hg_M}^{-M+R}(H_{\geq0}<H_{<-M+R})\min_{0\leq i \leq R-1}\Pa_{\Hg_M}^i(H_0<H_{[0,m_0-1]^c})
\\
\label{eqn:701}
&=\frac{\alpha_R}{d^{\totheright}}\Pa_{\Hg_M}^{-M+R}(H_{\geq0}<H_{<-M+R})b.
\end{align}
Again, the equality comes from either Lemma \ref{DERRW} or the expectation of a beta random variable, along with \eqref{eqn:821}. But arguing along the same lines as in Claim \ref{claim:ComparingGraphs}, we can see that paths that do not exit $(-M,M)$ have the same probability under $\Pa_{\G}^{-M+R}$ and $\Pa_{\Hg_M}^{-M+R}$. Therefore, we have
\begin{align}
 \notag   
 \Pa_{\Hg_M}^{-M+R}(H_{\geq0}<H_{<-M+R})
 &=\Pa_{\G}^{-M+R}(H_{\geq0}<H_{<-M+R})
\\
\notag
&\geq \Pa_{\G}^{-M+R}(H_{<-M+R}=\infty)
\\
\label{eqn:715}
&=\Pa_{\G}^0(H_{<0}=\infty)
>0.
\end{align}

Now \eqref{eqn:701} and \eqref{eqn:715} together  give us

\begin{equation*}\Pa_{\Hg_M}^{-M}(\tilde{H}_0\neq T_{-M\to0})
\geq\frac{\alpha_R}{d^{\totheright}}\Pa_{\G}^0(H_{<0}=\infty)b.
\end{equation*}
Because the last term is positive and does not depend on $M$, our claim is proven.

We now note, by Corollary \ref{cor:loopreversal} (2),
\begin{equation*}\Pa_{\Hg_M}^0(\tilde{H}_0= T_{-M\to0})=\Pa_{\check{\Hg}_M}^0(X_1=-M)=\frac{1}{1+c^{\totheleft}+c^{\totheright}}\end{equation*}
On the other hand, we also have $\Pa_{\Hg_M}^0(X_1=-M)=\frac{1}{1+c^{\totheleft}+c^{\totheright}}$.
Now we have

\begin{equation*}\Pa_{\Hg_M}^0(\tilde{H}_0= T_{-M\to0})\leq \Pa_{\Hg_M}^0(X_1=-M)\Pa_{\Hg_M}^{-M}(\tilde{H}_0= T_{-M\to0})+\Pa_{\Hg_M}^0(X_1\neq -M,H_{-M}<\tilde{H_0}),\end{equation*}
which can be rewritten as
\begin{equation}\label{eqn:804}
\frac{1}{1+c^{\totheleft}+c^{\totheright}}\leq \frac{1}{1+c^{\totheleft}+c^{\totheright}}\Pa_{\Hg_M}^{-M}(\tilde{H}_0= T_{-M\to0})+\Pa_{\Hg_M}^0(X_1\neq -M,H_{-M}<\tilde{H_0}).\end{equation}
For now, assume this claim.
\begin{claim}\label{claim:867}
The second term of the right hand side of \eqref{eqn:804} approaches 0 as $M\to\infty$.
\end{claim}

By Claim \ref{claim:recurrenceclaim}, $\Pa_{\Hg_M}^{-M}(\tilde{H}_0= T_{-M\to0})\leq c<1$ for all $M$. Taking the limsup as $M$ approaches infinity in $\eqref{eqn:804}$ therefore yields the contradiction
\begin{equation*}\frac{1}{1+c^{\totheleft}+c^{\totheright}}\leq \frac{1}{1+c^{\totheleft}+c^{\totheright}}c<\frac{1}{1+c^{\totheleft}+c^{\totheright}}.\end{equation*}
This completes the proof, pending Claim \ref{claim:867}, which we now prove. We must show that the expectation, under $P_{\Hg}$, of $P_{\omega}^0(X_1\neq-M,H_{-M}<\tilde{H}_0)$, approaches 0. Because $P_{\omega}^0(X_1\neq-M,H_{-M}<\tilde{H}_0)$ is bounded between 0 and 1, this is equivalent to showing that it approaches 0 in probability. 

First note that for $P_{\Hg}$--a.e. environment $\omega$, we have
\begin{equation*}
    P_{\omega}^0(X_1\neq-M,H_{-M}<\tilde{H}_0)\leq\max_{-L\leq x\leq -1}P_{\omega}^x(H_{-M}<H_0).
\end{equation*}
This is because a walk from 0 that does not immediately step to $-M$ but hits $-M$ before returning to 0 must first visit one of the sites in $[-L,-1]$. Now for almost every environment, we also have for $x\in[-L,-1]$,
\begin{align}
\notag
    P_{\omega}^x(H_{-M}<H_0)
    &\leq P_{\omega}^x(H_{-M}<H_{[0,R-1]})+P_{\omega}^x(H_{[0,R-1]}<H_{-M})\max_{0\leq i\leq R-1}P_{\omega}^i(H_{-M}<H_0)
    \\
\label{eqn:885}    
    &\leq \max_{ -L\leq y\leq -1}P_{\omega}^y(H_{-M}<H_{[0,R-1]})+\max_{0\leq i\leq R-1}P_{\omega}^i(H_{-M}<H_0).
\end{align}
But 
\begin{equation}\label{eqn:889}
    \max_{0\leq i\leq R-1}P_{\omega}^i(H_{-M}<H_0)
    \leq 
    \max_{0\leq i\leq R-1}P_{\omega}^i(H_{[-L,-1]}<H_0)\max_{-L\leq x\leq -1}P_{\omega}^x(H_{-M}<H_0),
\end{equation}
again because the walk must first visit $[-L,-1]$ if it is to hit $-M$ before 0. Combining \eqref{eqn:885} and \eqref{eqn:889}, we get
\begin{equation*}
    \hspace{-.6in}\max_{-L\leq x\leq -1}P_{\omega}^x(H_{-M}<H_0)
    \leq 
    \max_{ -L\leq y\leq -1}P_{\omega}^y(H_{-M}<H_{[0,R-1]})
    +
    \max_{0\leq i\leq R-1}P_{\omega}^i(H_{[-L,-1]}<H_0)
    \max_{-L\leq x\leq -1}P_{\omega}^x(H_{-M}<H_0).
\end{equation*}
Rearranging terms then gives us
\begin{equation*}
\max_{-L\leq x\leq -1}P_{\omega}^x(H_{-M}<H_0)
\leq
\frac{\max_{ -L\leq y\leq -1}P_{\omega}^y(H_{-M}<H_{[0,R-1]})}{1-\max_{0\leq i\leq R-1}P_{\omega}^i(H_{[-L,-1]}<H_0)}.
\end{equation*}
The numerator and denominator are independent under $P_{\Hg_M}$. Since the denominator is almost surely positive, it suffices to show that the numerator approaches 0 in probability. But again arguing as in Claim \ref{claim:ComparingGraphs}, we see that the distribution of $\max_{ -L\leq y\leq -1}P_{\omega}^y(H_{-M}<H_{[0,R-1]})$ under $P_{\Hg_M}$ is the distribution of $\max_{ -L\leq y\leq -1}P_{\omega}^y(H_{\leq-M}<H_{[0,R-1]})$ under $P_{\G}$. By non-transience to the left, the latter approaches 0 $P_{\G}$--a.s. as $M$ increases. This proves Claim \ref{claim:867}, and with it Theorem \ref{thm:transiencecriterion}.
\end{proof}

\setcounter{section}{\value{SectionSave}}
\setcounter{exmp}{\value{ExampleSave}}

\section{An abstract ballisticity criterion}\label{sec:abstract}

The main goal of this section is to prove Lemma \ref{lem:ballistic}, which says that a walk is ballistic if and only if $\Ea^0[N_0]<\infty$.





Before we can characterize when the almost-sure limiting velocity $v$ is positive, we must first note that it exists.
This has been shown under an ellipticity assumption too strong for our model \cite{Bremont2004}, but it can be proven in the more general case with standard techniques. 
The proof for the recurrent case (where, necessarily, $v=0$) can be done by a slight modification of arguments in \cite{Zerner2002}. The proof for the directionally transient case follows \cite{Kesten1977} in defining {\em regeneration times} $(\tau_k)_{k=0}^{\infty}$. Let $\tau_0:=0$, and for $k\geq1$, define 
\begin{equation}\label{eqn:RegenerationTimes}
    \tau_k:=\min\{n>\tau_{k-1}:X_n>X_j\text{ for all }j<n,X_n\leq X_j\text{ for all }j>n\}.
\end{equation}
A crucial fact is that the sequences ($X_{\tau_n}-X_{\tau_{n-1}})_{n=2}^{\infty}$ and $(\tau_n-\tau_{n-1})_{n=2}^{\infty}$ are i.i.d. Using these regeneration times, we are able to derive a formula for $v$, as well as a characterization in terms of hitting times.

\begin{prop}\label{prop:LimitingVelocity} Let $P$ be a probability measure on $\Omega_{\Z}$ satisfying (\hyperlink{cond:C1}{C1}), (\hyperlink{cond:C2}{C2}), and (\hyperlink{cond:C3}{C3}). Then the following hold:
\begin{enumerate}
    \item There is a $\Pa^0$--almost sure limiting velocity 
\begin{equation}\label{eqn:856}
    v:=\lim_{n\to\infty}\frac{X_n}{n}=\frac{\Ea^0[X_{\tau_{2}}-X_{\tau_1}]}{\Ea[\tau_{2}-\tau_{1}]},
\end{equation}
where the numerator is always finite, and the fraction is understood to be 0 if the denominator is infinite.
    \item $\lim_{x\to\infty}\frac{H_{\geq x}}{x} = \frac{1}{v}$, where $\frac{1}{v}$ is understood to be $\infty$ if $v=0$.
\end{enumerate}
\end{prop}
We outline some details of the argument for Proposition \ref{prop:LimitingVelocity} in Appendix \ref{append:Leftoverproofs}.

For the rest of this section, assume $P$ satisfies (\hyperlink{cond:C1}{C1}), (\hyperlink{cond:C2}{C2}), (\hyperlink{cond:C3}{C3}), and (\hyperlink{cond:C4}{C4}).
We also use regeneration times to derive the following lemma.

\begin{lem}\label{known}
For any $a,c\in\Z$,
\begin{equation*}\lim_{x\to\infty}\frac{1}{x}\sum_{k=c}^xN_k=\frac1{v},~\Pa^a\text{--a.s.}\end{equation*}
If $v=0$, then the limit is infinity.
\end{lem}

\begin{proof}
Fix $a$. Recall that $N_{k}^{(-\infty,x)}$ is the amount of time the walk spends at $k$ before $H_{\geq x}$.
Then for $x>c$,


\begin{equation*}\frac{H_{\geq x}}{x}=\frac1x\sum_{k=-\infty}^{c-1}N_k^{(-\infty,x)}+\frac1x\sum_{k=c}^{x-1}N_k^{(-\infty,x)}.
\end{equation*}
The first term approaches 0 almost surely by assumption (\hyperlink{cond:C4}{C4}); hence, by Proposition \ref{prop:LimitingVelocity} (2),

\begin{equation}\label{eqn:890}
\lim_{x\to\infty}\frac1x\sum_{k=c}^{x-1}N_k^{(-\infty,x)}=\frac1{v}, \Pa^a\text{--a.s.}
\end{equation}

We note that $N_k$ and $N_k^{(-\infty,x)}$ differ only if the walk backtracks and visits $k$ after reaching $[x,\infty)$. The sum, over all $k<x$, of these differences, is the total amount of time the walk spends to the left of $x$ after $H_{\geq x}$, and it is bounded above by the time from $H_{\geq x}$ to the next regeneration time (defined as in \eqref{eqn:RegenerationTimes}), which is in turn bounded above by $\tau_{J(x)}-\tau_{J(x)-1}$, where $J(x)$ is the (random) $j$ such that $\tau_{j-1}\leq H_{\geq x}<\tau_j$. Hence

\begin{equation}\label{eqn:896}
\frac1x\sum_{k=c}^{x-1}N_k^{(-\infty,x)}
\leq \frac1x\sum_{k=c}^{x-1}N_k
\leq \frac1x\sum_{k=c}^{x-1}N_k^{(-\infty,x)} + \frac1x[\tau_{J(x)}-\tau_{J(x)-1}]
\end{equation}

Assume $v=0$. Then by \eqref{eqn:890}, the left side of \eqref{eqn:896} approaches $\infty$ as $x$ approaches $\infty$, and therefore so does the middle. 
On the other hand, suppose $v>0$. By \eqref{eqn:856}, $\Ea[\tau_{2}-\tau_{1}]<\infty$.
Then by the strong law of large numbers, $\frac{\tau_n}{n}\to \Ea[\tau_{2}-\tau_{1}]<\infty$, which implies that $\frac{\tau_n-\tau_{n-1}}{n}$ approaches 0. Since $J(x)\leq x+1$, the term $\frac1x[\tau_{J(x)}-\tau_{J(x)-1}]$ approaches zero almost surely; hence the Squeeze Theorem yields the desired result.
\end{proof}

Suppose for now that $R=1$. Then, for almost every $\omega$, it is possible to define a bi-infinite walk $\BIX=(\biX_n)_{n\in\Z}$ whose ``right halves" are distributed like random walks under $\omega$. From each site $a$, run a walk according to the transition probabilities given by $\omega$ until it reaches $a+1$ (which occurs in finite time $P_{\omega}^a$--a.s. for $P$--a.e. $\omega$). Concatenating all of these walks then gives, up to a time shift\footnote{Choose, for example, the time shift where $\biX_0=0$ and where $\biX_n<0$ whenever $n<0$.}, a unique walk $\BIX=(\biX_n)_{n\in\Z}$ such that for any $x\in\Z$, the distribution of $(\biX_k)_{k=n}^{\infty}$, conditioned on $\biX_n=x$, is $P_{\omega}^x$. We may think of $\BIX$ as a walk from $-\infty$ to $\infty$ in the environment $\omega$.

\label{page:cascade}With a bit more work, we can define a similar bi-infinite walk in the general case $R>0$. Call the set of vertices $((k-1)R,kR]$ the $k$th {\em level} of $\Z$, and for $x\in\Z$, let $[[x]]_R$ denote the level containing $x$. Let $\omega$ be a given environment. From each point $a\in\Z$, run a walk according to the transition probabilities given by $\omega$ until it reaches the next level (i.e., $[[a+R]]_R$). This will happen $P_{\omega}^a$--a.s. for $P$--a.e. $\omega$, by transience to the right and because it is not possible to jump over a set of length $R$. Do this independently at every point for every level. This gives what we'll call a {\em cascade}: a set of (almost surely finite) walks indexed by $\Z$, where the walk indexed by $a\in\Z$ starts at $a$ and ends upon reaching level $[[a+R]]_R$. Then for almost every cascade, concatenating these finite walks gives, for each point $a$, a right-infinite walk ${\bf X}^{a}=(X_n^a)_{n=0}^{\infty}$. Let $P_{\omega}$ be the probability measure we have just described on the space of cascades, and let $\Pa=P\times P_{\omega}$. 

It is crucial to note that by the strong Markov property, the law of ${\bf X}^{a}$ under $P_{\omega}$ is the same as the law of ${\bf X}$ under $P_{\omega}^{a}$, which also implies that the law of ${\bf X}^{a}$ under $\Pa$ is the same as the law of ${\bf X}$ under $\Pa^{a}$.


For each $x\in\Z$, let the ``coalescence event" $C_x$ be the event that all the walks from level $[[x-R]]_R$ first hit level $[[x]]_R$ at $x$. On the event $C_x$, we say a coalescence occurs at $x$.

\begin{lem}\label{prop:coalescences}
Let $\mathcal{E}_1$ be the event that all the ${\bf X}^a$ are transient to the right, that all steps to the left and right are bounded by $L$ and $R$, respectively, and that infinitely many coalescences occur to the left and to the right of 0. Then $\Pa(\mathcal{E}_1)=1$.
\end{lem}

\begin{proof}
Boundedness of steps has probability 1 by assumption (\hyperlink{cond:C3}{C3}), and by assumption (\hyperlink{cond:C4}{C4}) all the walks ${\bf X}^a$ are transient to the right with probability 1. Now for $k\geq 2$ and $x\in\Z$, let $C_{x,k}$ be the event that all the walks from level $[[x-R]]_R$ first hit level $[[x]]_R$ at $x$ without ever having reached $[[x-kR]]_R$. Choose $k$ large enough that $\Pa(C_{0,k})>0$; then under the law $\Pa$, the events $\{C_{nkR,k}\}_{n\in\Z}$ are all independent and have equal, positive probability. Thus, infinitely many of them will occur in both directions, $\Pa$--a.s. By definition, $C_{x,k}\subset C_x$, and so infinitely many of the events $C_x$ occur in both directions, $\Pa$--a.s.
\end{proof}

Assume the environment and cascade are in the event $\mathcal{E}_1$. Let $(x_k)_{k\in\Z}$ be the locations of coalescence events (with $x_0$ the smallest non-negative $x$ such that $C_x$ occurs). By definition of the $x_k$, for every $k$ and for every $a$ to the left of $[[x_k]]_R$, $H_{[[x_k]]_R}({\bf X}^a)=H_{x_k}({\bf X}^a)<\infty$. Now for $j<k$, it necessarily holds that $x_j$ is to the left of $[[x_k]]_R$, since there can be only one $x_k$ per level. Define $\nu(j,k):=H_{x_k}({\bf X}^{x_j})$. By definition of the walks ${\bf X}^a$, we have for $j<k$, $n\geq0$, 
\begin{equation}\label{eqn:996}
X_{n+\nu(j,k)}^{x_j}=X_n^{x_k}.
\end{equation}
From this one can easily check that the $\nu(j,k)$ are additive; that is, for $j<k<\ell$, we have $\nu(j,\ell)=\nu(j,k)+\nu(k,\ell)$. 
Because all the ${\bf X}^{x_k}$ agree with each other in the sense of \eqref{eqn:996}, we may define a single, bi-infinite walk $\BIX=(\biX_n)_{n\in\Z}$ that agrees with all of the ${\bf X}^{x_k}$. For $n\geq0$, let $\biX_n=X_n^{x_0}$. For $n<0$, choose $j<0$ such that $\nu(j,0)>|n|$, and let $X_n=X_{\nu(j,0)-|n|}^{x_j}$. This definition is independent of the choice of $j$, because if $j<k<0$ with $v(k,0)>|n|$, then by \eqref{eqn:996} and the additivity of the $\nu(j,k)$, we have
\begin{equation*}
    X_{\nu(j,0)-|n|}^{x_j}=X_{\nu(j,k)+\nu(k,0)-|n|}^{x_j}=X_{\nu(k,0)-|n|}^{x_k}.
\end{equation*}
We may then define $\overline{N}_{x}:=\#\{n\in\Z:X_n=x\}$ to be the amount of time the walk $\BIX$ spends at $x$. Thus, $\overline{N}_x=\lim_{a\to-\infty}N_x({\bf X}^a)$.

\begin{lem}\label{prop:Ergodic}
Both of the sequences $({\bf X}^a)_{a\in\Z}$ and $(\overline{N}_x)_{x\in\Z}$ are stationary and ergodic. 
\end{lem}

\begin{proof}
For a given environment, the cascade that defines $\BIX$ may be generated by a (countable) family ${\bf U}=\left(U_n^{a}\right)_{n\in\N,a\in\Z}$ of i.i.d. uniform random variables on $[0,1]$. For such a collection, and an $a\in\Z$, let ${\bf U}^a$ be the projection $\left(U_n^{a}\right)_{n\in\N}$. Given an environment $\omega$, the finite walk from $a$ to level $[[a+R]]_R$ may be generated using the first several $U_n^{a}$. (One of the $U_n^{a}$ is used for each step. Once the walk terminates, the rest of the $U_n^{a}$ are not needed, but one does not know in advance how many will be needed.) Let $\hat{\omega}^x=(\omega^x,{\bf U}^x)$, and $\hat{\omega}=(\hat{\omega}^x)_{x\in\Z}$. Define the left shift $\hat{\theta}$ by $\hat{\theta}(\hat{\omega}):=(\hat{\omega}^{x+1})_{x\in\Z}$.
Then $(\hat{\omega}^x)_{x\in\Z}$ is an i.i.d. sequence. We have ${\bf X}^0={\bf X}^0(\hat{\omega})$ and ${\bf X}^a={\bf X}^0(\hat{\theta}^a\hat{\omega})$. Similarly, $\overline{N}_0=\overline{N}_0(\hat{\omega})$ and $N_x=\overline{N}_0(\hat{\theta}^x\hat{\omega})$. So it suffices to show that ${\bf X}^0$ and $\overline{N}_0$ are measurable.
The measurability of ${\bf X}^0$ is obvious. For $\overline{N}_0$, let $A_{k,\ell,B,r}$ be the event that:
\begin{enumerate}[(a)]
    \item for some $x<0$, a coalescence event $C_{x,k}$ (as defined in the proof of Lemma \ref{prop:coalescences}) occurs with $-B\leq x-kR<x<0$, so that $\BIX$ agrees with ${\bf X}^x$ to the right of $x$;
    \item $N_0^{[-B,B]}({\bf X}^x)\geq \ell$, where $N_0^{[-B,B]}$ is the amount of time the walk spends at $x$ before exiting $[-B,B]$; and
    \item none of the walks from sites $a\in[-B,B]$ uses more than $r$ of the random variables  $U_r^{a}$.
\end{enumerate}
On this event, $\overline{N}_0$ is seen to be at least $\ell$ by looking only within $[-B,B]$ and only at the first $r$ uniform random variables at each site. The event $A_{k,\ell,B,r}$ is measurable, because it is a measurable function of finitely many random variables, and the event $\{\overline{N}_0>\ell\}$ is, up to a null set, simply the union over all $r$, then over all $B$, and then over all $k$ of these events. Thus, $\overline{N}_0$ is measurable.
\end{proof}

We now give the connection between $\overline{N}_0$ and the limiting velocity $v$. 

\begin{lem}\label{lem:LimitingVelocity}
 $v=\frac{1}{\Ea[\overline{N}_0]}$. Consequently, the walk is ballistic if and only if $\Ea[\overline{N}_0]<\infty$.
 \end{lem}

We note that a similar formula for the limiting speed in the ballistic case can be obtained from \cite[Theorem 6.12]{Dolgopyat&Goldsheid2019} for discrete-time RWRE on a strip, although the probabilistic interpretation is less explicit, and an ellipticity assumption that does not hold for Dirichlet RWRE is required. 

\begin{proof}
By Lemma \ref{prop:Ergodic} and Birkhoff's Ergodic theorem, for any $c\in\Z$ we have

\begin{equation*}\lim_{n\to\infty}\frac1n\sum_{k=c}^n\overline{N}_{k}=\Ea[\overline{N}_0],~\Pa\text{--a.s.}\end{equation*}

Fix $a\in\Z$. For large enough $k$, $N_k({\bf X}^a)=\overline{N}_{k}$. We therefore get
\begin{equation*}\lim_{n\to\infty}\frac1n\sum_{k=c}^nN_k({\bf X}^a)=\Ea[\overline{N}_0],~\Pa\text{--a.s.}\end{equation*}
It follows that
\begin{equation*}\lim_{n\to\infty}\frac1n\sum_{k=c}^nN_k({\bf X})=\Ea[\overline{N}_0],~\Pa^a\text{--a.s.}\end{equation*}
By Lemma \ref{known}, we get $v=\frac{1}{\Ea[\overline{N}_{0}]}$.
\end{proof}

Now we can see that the walk is ballistic if and only if $\Ea[\overline{N}_0]<\infty$. In order to prove Lemma \ref{lem:ballistic}, we need to compare $\Ea[\overline{N}_0]$ with $\Ea^0[N_0]$. 

\begin{lem}\label{lem:limitfinitecase}
$\Ea[\overline{N}_0]\leq\Ea^0[N_0]$. 
\end{lem}

\begin{proof}
If $\Ea^0[N_0]=\infty$, the inequality is trivial. Assume, therefore, that $\Ea^0[N_0]<\infty$. 

Note that ${\lim_{x\to\infty}N_0({\bf X}^{-x})=\overline{N}_0}$, $\Pa$--a.s. Assuming we are able to interchange a limit with an expectation, we have

\begin{align}
\notag
\Ea[\overline{N}_0] &= \Ea\left[\lim_{x\to\infty}N_0({\bf X}^{-x})\right] 
\\
\label{eqn:1074}
&= \lim_{x\to\infty}\Ea\left[N_0({\bf X}^{-x})\right] 
\\
\notag
&= \lim_{x\to\infty}\Ea^{-x}[N_0({\bf X})]. 
\end{align}
But each term $\Ea^{-x}[N_0({\bf X})]=E[E_{\omega}^{-x}[N_0]]$ is less than $\Ea^0[N_0]=E[E_{\omega}^0[N_0]]$, since $E_{\omega}^{-x}[N_0]=P_{\omega}^{-x}(H_0<\infty)E_{\omega}^{0}[N_0]$.
Therefore, we may conclude $\Ea[\overline{N}_0]\leq\Ea^0[N_0]$, provided we can justify \eqref{eqn:1074}. To do this, we will apply the dominated convergence theorem, noting that $N_0({\bf X}^{-x})\leq \max_{1-R\leq y\leq 0}N_0({\bf X}^{y})$ for all ${x>R}$. To see that the latter has finite expectation, we have

\begin{align*}
\Ea\left[\max_{1-R\leq y\leq 0}N_0({\bf X}^{y})\right]&\leq
\sum_{y=1-R}^0\Ea\left[N_0({\bf X}^{y})\right]
\\
&=\sum_{y=1-R}^0E\left[E_{\omega}[N_0({\bf X}^{y})]\right]
\\
&=\sum_{y=1-R}^0E\left[E_{\omega}^{y}[{N}_0]\right]
\\
&=\sum_{y=1-R}^0E\left[P_{\omega}^{y}(H_{0}<\infty)E_{\omega}^0[{N}_0]\right]
\\
&\leq \sum_{y=1-R}^0E\left[E_{\omega}^0[{N}_0]\right]
\\
&=R\Ea^0[{N}_0]
\\
&<\infty.
\end{align*}
This allows us to justify our use of the dominated convergence theorem, completing the proof.
\end{proof}

We must now handle the case where $\Ea^0[N_0]=\infty$. Our first step is to prove Lemma \ref{lem:BremontLemma}, which states that $v>0$ if and only if $\Ea^0[H_{\geq1}]<\infty$.

\setcounter{ExampleSave}{\value{exmp}}
\setcounter{SectionSave}{\value{section}}
\setcounter{section}{1}
\setcounter{exmp}{2}
\setcounter{claim}{0}

\begin{proof}[Proof of Lemma \ref{lem:BremontLemma}]\hypertarget{proof:BremontLemma}
Suppose $\Ea^0[H_{\geq1}({\bf X})]<\infty$. 
Then for $P_{\omega}$--almost every cascade, we have
\begin{align}
\notag
\frac{H_{\geq x}({\bf X}^0)}{x}
&=\frac{1}{x}\sum_{k=1}^x(H_{\geq k}({\bf X}^0)-H_{\geq k-1}({\bf X}^0))
\\
\label{eqn:k-1}
&\leq \frac{1}{x}\sum_{k=1}^xH_{\geq k}({\bf X}^{k-1}),
\end{align}
where the inequality comes from the fact that if ${\bf X}^0$ hits $[k-1,\infty)$ at $k-1$, then it follows the same path from there as ${\bf X}^{k-1}$, while if it hits $[k-1,\infty)$ at a point to the right of $k-1$, then $H_{\geq k}-H_{\geq k-1}=0$. By Birkhoff's Ergodic Theorem, the right side $\Pa$--a.s. approaches $\Ea[H_{\geq 1}({\bf X}^0)]=\Ea^0[H_{\geq1}]$. Now we know from Proposition \ref{prop:LimitingVelocity} that $\lim_{n\to\infty}\frac{n}{X_n}=\frac{1}{v}$, so the subsequence $\frac{H_{\geq x}({\bf X}^0)}{X_{H_{\geq x}}^0}$ must have the same limit. Since, for $x>R$, we have
\begin{equation*}
    \frac{H_{\geq x}({\bf X}^0)}{X_{H_{\geq x}}^0+R}\leq \frac{H_{\geq x}({\bf X}^0)}{x}\leq\frac{H_{\geq x}({\bf X}^0)}{X_{H_{\geq x}}^0-R},
\end{equation*}
we get $\frac{1}{v}=\lim_{x\to\infty}\frac{H_{\geq x}({\bf X}^0)}{x}$. Applying \eqref{eqn:k-1}, we get
\begin{align*}
    \frac{1}{v}&=\lim_{n\to\infty}\frac{H_{\geq x}({\bf X}^0)}{x}
    \\
    &\leq \lim_{n\to\infty}\frac{1}{x}\sum_{k=1}^xH_{\geq k}({\bf X}^{k-1})
    \\
    &=\Ea^0[H_{\geq 1}].
\end{align*}
Therefore, if $\Ea^0[H_{\geq 1}]<\infty$, then $v>0$.

On the other hand, suppose $\Ea^0[H_{\geq 1}]=\infty$. We will show that $v=0$. 
\begin{claim}
$\Ea\left[\min_{1\leq i\leq R}H_{\geq R+1}({\bf X}^i)\right]=\infty$.
\end{claim}
By assumptions (\hyperlink{cond:C1}{C1}), (\hyperlink{cond:C2}{C2}), and (\hyperlink{cond:C3}{C3}), we still have an $m_0\geq\max(L,R)$ large enough that every interval of length $m_0$ is irreducible, $P$--a.s.
Let $A$ be the event that
\begin{itemize}
    \item For each $i=1,\ldots, R-1$, the walk ${\bf X}^i$ hits $R$ before leaving $[R-m_0+1, R]$. 
    \item The walk ${\bf X}^R$ first exits $[R-m_0+1,R]$ by hitting $R-m_0$. 
\end{itemize}
Then under $P$, the quenched probability of $A$ is independent of $\omega^{\leq R-m_0}$. Now, on the event $A$, the minimum $\min_{1\leq i\leq R}H_{\geq R+1}({\bf X}^i)$ is attained for $i=R$, since all the other walks take time to get to $R$ and then simply follow ${\bf X}^R$. 
Now on $A$, $H_{\geq R+1}({\bf X}^R)$ is greater than the amount of time it takes for the walk ${\bf X}^R$ to cross back to $[R-m_0+1,\infty)$ after first hitting $R-m_0$. The quenched expectation of this time, conditioned on $A$, is $E_{\omega}^{R-m_0}[H_{\geq R-m_0+1}]$ by the strong Markov property, and this depends only on $\omega^{\leq R-m_0}$. Hence
\begin{align}
\notag
    \Ea\left[\min_{1\leq i\leq R}H_{\geq R+1}({\bf X}^i)\right]
    &\geq E\left[P_{\omega}(A)E_{\omega}[H_{\geq R+1}({\bf X}^R)|A]\right]
    \\
\notag
    &\geq E\left[P_{\omega}(A)E_{\omega}^{R-m_0}[H_{\geq R-m_0+1}({\bf X})]\right]
    \\
\notag
    &= \Pa(A)\Ea^{R-m_0}[H_{\geq R-m_0+1}({\bf X})]
    \\
\notag    
    &=\Pa(A)\Ea^{0}[H_{\geq 1}({\bf X})]
    \\
\notag
    &=\infty.
\end{align}
This proves our claim. 
Now for $x\geq1$, 
\begin{equation*}
    H_{\geq xR+1}({\bf X}^0)\geq H_{\geq1}({\bf X}^0)+\sum_{k=1}^x\min_{1\leq i \leq R}H_{\geq kR+1}({\bf X}^{(k-1)R+i}).
\end{equation*}
Dividing by $xR$ and taking limits as $x\to\infty$, we get $\lim_{x\to\infty}\frac{H_{\geq xR+1}({\bf X}^0)}{xR}=\infty$, $\Pa$--a.s. by Birkhoff's ergodic theorem. Hence $\lim_{x\to\infty}\frac{H_{\geq xR+1}({\bf X})}{xR+1}=\infty$, $\Pa^0$--a.s. It follows that $v=0$. 
\end{proof}
\setcounter{exmp}{\value{ExampleSave}}
\setcounter{section}{\value{SectionSave}}

Now we can handle the case $\Ea^0[N_0]=\infty$.

\begin{prop}\label{prop:4.8}
If $\Ea^0[N_0]=\infty$, then $v=0$. \end{prop}
\begin{proof}
Suppose $\Ea^0[N_0]=\infty$. We want to show that $v=0$. By Lemma \ref{lem:BremontLemma}, it suffices to show that $\Ea^0[H_{\geq1}]=\infty$. 

Now $N_0$ is the total number of visits the walk makes to 0. These visits may be sorted based on the farthest point to the right that the walk has hit in the past at the time of each visit. 
As in the proof of Lemma \ref{known}, we use $N_{k}^{(-\infty,x)}$ to denote the amount of time the walk spends at $k$ before $H_{\geq x}$.
Thus, for a walk started at 0 we get
\begin{equation}\label{eqn:1207}
    N_0=\sum_{x=0}^{\infty}\left(N_0^{(-\infty,x+1)}-N_0^{(-\infty,x)}\right).
\end{equation}
Taking expectations on both sides, we get
\begin{equation}\label{eqn:1213}
    \Ea^0[N_0]
    =\sum_{x=0}^{\infty}E\left[E_{\omega}^0\left[N_0^{(-\infty,x+1)}-N_0^{(-\infty,x)}\right]\right]
\end{equation}
Now $N_0^{(-\infty,x)}$ and $N_0^{(-\infty,x+1)}$ can only differ if the walk hits $[x,\infty)$ at $x$. Conditioned on this event, the distribution under $P_{\omega}^0$ of the walk $(X_{n+H_{\geq x}})_{n=0}^{\infty}$ is the distribution of ${\bf X}$ under $P_{\omega}^x$. Thus, 
\begin{equation}\label{eqn:1218}
    E_{\omega}^0\left[N_0^{(-\infty,x+1)}-N_0^{(-\infty,x)}\right]
    =
    P_{\omega}^0(X_{H_{\geq x}}=x)E_{\omega}^x\left[N_0^{(-\infty,x+1)}\right].
\end{equation}
Combining \eqref{eqn:1213} and \eqref{eqn:1218}, we get
\begin{align*}
    \Ea^0[N_0]
    &=\sum_{x=0}^{\infty}E\left[P_{\omega}^0(X_{H_{\geq x}}=x)E_{\omega}^x\left[N_0^{(-\infty,x+1)}\right]\right]
    \\
    &\leq \sum_{x=0}^{\infty}E\left[E_{\omega}^x\left[N_0^{(-\infty,x+1)}\right]\right]
    \\
    &=\sum_{x=0}^{\infty}\Ea^x\left[N_0^{(-\infty,x+1)}\right].
\end{align*}
By stationarity, 
\begin{align*}
    \Ea^0[N_0]
    &\leq\sum_{x=0}^{\infty}\Ea^0\left[N_{-x}^{(-\infty,1)}\right]
    \\
    &=\Ea^0\left[\sum_{x=0}^{\infty}N_{-x}^{(-\infty,1)}\right]
    \\
    &=\Ea^0[H_{\geq1}].
\end{align*}
If $\Ea^0[N_0]=\infty$, it follows that $\Ea^0[H_{\geq1}]$, and by Lemma \ref{lem:BremontLemma}, $v=0$. 
\end{proof}

We can now complete the proof of our main lemma.

\begin{proof}[Proof of Lemma \ref{lem:ballistic}]\hypertarget{proof:ballistic}
Assume the walk is transient to the right. 
If $\Ea^0[{N}_0]=\infty$, the conclusion is that of Proposition \ref{prop:4.8}. Otherwise, combining Lemmas \ref{lem:LimitingVelocity} and \ref{lem:limitfinitecase} gives $v>0$.
The left-transient case follows by symmetry. By the 0-1 law of \cite{Key1984}, the remaining case is where the walk is recurrent. This implies that $N_0=\infty$, $\Pa^0$--a.s., so that $\Ea^0[N_0]=\infty$. Because $v=0$ in the recurrent case, the lemma is true.
\end{proof}

\section{Ballistic parameters}\label{sec:ballisticity}

We now return to the Dirichlet model. In this section, we will characterize ballisticity in terms of $L$, $R$, and the parameters $(\alpha_i)_{i=-L}^R$. 
Lemma \ref{lem:ballistic} tells us that the walk is ballistic precisely when the quantity $\Ea_{\G}^0[N_0]=E_{\G}[E_{\omega}^0[N_0]]$ is finite. Although we cannot usually calculate this expectation, we are able to characterize when it is finite in terms of our Dirichlet parameters. We assume throughout this section that $\kappa_1>0$, so that the walk is transient to the right, and we examine the integrability of $E_{\omega}^0[N_0]$ under $P_{\G}$. In fact, we generalize the question, examining when $E_{\G}[E_{\omega}^0[N_0]^s]<\infty$ for $s>0$. The goal of this section is to prove Theorems \ref{thm:gamma}, \ref{thm:TFAE2}, and \ref{thm:MainMomentTheorem}. Theorem
\ref{thm:MainBallisticTheorem} then easily follows from Theorem \ref{thm:MainMomentTheorem} and Lemma \ref{lem:ballistic}.

\subsection{Finite traps: the parameter \texorpdfstring{$\kappa_0$}{kappa0}}\label{subsec:FiniteTraps}

In this subsection, we use Tournier's lemma to study the existence of finite traps---finite sets in which the walk is expected to spend an infinite amount of time before exiting under the annealed measure.

We note that although Tournier's lemma is stated for a finite graph, it can readily be applied to walks that are killed upon exiting a finite subset of an infinite graph. In particular, for any positive integer $M$, the annealed expected number of visits to $0$ before exiting $[-M,M]$ is infinite if and only if there is a strongly connected subset $S$ of $[-M,M]$ containing 0 such that $\beta_S\leq1$. 

This discussion motivates us to define, for our graph $\mathcal{G}$, the quantity 
\begin{equation}\label{eqn:kappa0def}
    \kappa_0:=\inf\{\beta_S:S\subset \Z\text{ finite, strongly connected}\}.
\end{equation}
Recall that Theorem \ref{thm:gamma} states the equivalence of the following statements (where $N_0^S$ is the amount of time the walk spends at 0 before first exiting $S$):
\begin{enumerate}[(a)]
    \item $\kappa_0\leq s$.
    \item For all sufficiently large $M$, $E_{\G}\left[E_{\omega}^0\left[N_0^{[-M,0]}\right]^s\right]=\infty$.
    \item For some $M\geq0$, $E_{\G}\left[E_{\omega}^0\left[N_0^{[-M,M]}\right]^s\right]=\infty$.
\end{enumerate}
The proof is essentially a straightforward application of Tournier's lemma as we have just discussed; however, in order to handle the boundary case $s=\kappa_0$, we first need to show that the infimum in the definition of $\kappa_0$ is actually a minimum. For example, in the case $\kappa_0=s=1$, showing that $\kappa_0$ is a minimum means showing that there is actually a finite set containing 0 that the walk is expected to get stuck in for an infinite amount of time. We also give an algorithm to compute $\kappa_0$.

\begin{prop}\label{prop:GammaAttained}
The infimum $\kappa_0$ for the graph $\G$ is actually a minimum attained by a set $S\subset \Z$. Moreover, there is an integer $M$, which may be calculated from $L$, $R$, and the weight assignments $(\alpha_{i})_{-L\leq i\leq R}$, such that the infimum is attained on a subset $S$ with diameter at most $M$. Hence $\kappa_0$ can be calculated directly.
\end{prop}

By translation invariance of the graph $\G$, this implies that it is possible to compute $\kappa_0$ by looking at strongly connected subsets $S$ with diameter no more than $M$ and with leftmost point $0$.

\setcounter{ExampleSave}{\value{exmp}}
\setcounter{exmp}{\value{thm}}
\begin{proof}

We prove this in a series of claims. Recall that $m_0\geq\max(L,R)$ is large enough that every interval of length $m_0$ is strongly connected.

\begin{claim}\label{claim:GammaClaim1}
 $\kappa_0\leq d^{\totheright}+d^{\totheleft}$. 
\end{claim}

To prove this claim, it suffices to exhibit a finite, strongly connected set $S\subset\Z$ with $\beta_S=d^{\totheright}+d^{\totheleft}$. Let $S=[0,m_0-1]$. Then $S$ is strongly connected. Now $\beta_S$ is the total weight of edges from $[0,m_0-1]$ to other vertices. Since $m\geq\max(L,R)$, it is easy to check that this is exactly $d^{\totheright}+d^{\totheleft}$.

\begin{claim}\label{claim:GammaClaim3}
 Let $S\subset\Z$ be a finite, strongly connected set of vertices. If $x$ is a vertex to the left or to the right of $S$, then $\beta_{S\cup\{x\}}\geq\beta_{S}$. 
\end{claim}

The quantity $\beta_S$ is the sum of all weights from vertices in $S$ to vertices not in $S$. The quantity $\beta_{S\cup\{x\}}$ counts all same weights, except for weights of edges from $S$ to $x$, and it also counts weights of edges from $x$ to vertices not in $S\cup\{x\}$. 
If $x$ is to the right of $S$, then the total weight of edges from $S$ to $x$ cannot be more than $c^{\totheright}$, because $c^{\totheright}$ is the total weight into $x$ from all vertices to the left of $x$. On the other hand, $c^{\totheright}$ is also the total weight from $x$ to all vertices to the right of $x$, which are necessarily not in $S\cup\{x\}$. Thus, the additional weight from $x$ to the right at least makes up for any weight into $x$ from $S$. This proves the claim in the case that $x$ is to the right of $A$, and a similar argument proves the symmetric case.
\begin{rem}
Note the importance of the assumption that $x$ is to the left or right of $S$. If $x$ is in between some of the vertices of $S$, then it is certainly possible that  $\beta_{S\cup\{x\}}<\beta_{S}$. See Examples \ref{kappa0example9} and \ref{kappa0example6}. 
\end{rem}
\begin{claim}\label{claim:GammaClaim4}
Let $S$ be a finite, strongly connected subset of $\Z$. Say that a vertex $x\in S$ is {\em insulated} if every site reachable in one step from $x$ is also in $S$. Then if $x<y$ are consecutive non-insulated vertices in $S$ with $y-x>m_0$, it must be the case that all vertices between $x$ and $y$ are in $S$.
\end{claim}

Suppose there are two consecutive non-insulated vertices $x,y\in S$ with $y-x>m_0$. Because $m_0\geq\max(L,R)$, there must be other vertices from $S$ strictly between $x$ and $y$ in order for $S$ to contain a path from $x$ to $y$ and $y$ to $x$. By assumption, all such vertices are insulated. Therefore, if there is an edge from a vertex in $(x,y)\cap S$ to another vertex in $(x,y)$, then the latter vertex must also be in $S$, and since it is strictly between $x$ and $y$, it must also be insulated. Applying this fact repeatedly, we see that any two vertices that communicate within $(x,y)$ are either both insulated in $S$ or both in $S^c$. Since the length of $(x,y)$ is at least $m_0$, all sites in the interval communicate, and so all are in $S$.

\begin{claim}\label{claim:GammaClaim5}
The infimum $\kappa_0$ is attained as a minimum; $\kappa_0=\beta_{S_0}$ for some $S_0$. Moreover, there is an algorithm to find it.
\end{claim}
Let $\varepsilon$ be the smallest weight any edge in $\G$ has, let $N$ be an integer such that $N\varepsilon\geq d^{\totheright}+d^{\totheleft}$, and let $M=(N-1)(m_0)$ (note this implies $M\geq m_0$).
Then if $S$ is a set of vertices with diameter greater than $M$, it must either have at least $N$ non-insulated vertices or have consecutive non-insulated vertices that differ by more than $m_0$. If there are at least $N$ non-insulated vertices, then there is an edge from each of these to at least one vertex outside of $S$, which means $\beta_{S}\geq N\varepsilon\geq d^{\totheright}+d^{\totheleft}\geq\kappa_0$, the last inequality coming from Claim \ref{claim:GammaClaim1}. On the other hand, if two consecutive non-insulated vertices $x<y$ differ by more than $m_0$, then $[x,y]\subseteq S$ by Claim \ref{claim:GammaClaim4}. Now $\beta_{[x,y]}=d^{\totheright}+d^{\totheleft}$, and vertices to the left and to the right of $[x,y]$ can only increase $\beta_{S}$ by Claim \ref{claim:GammaClaim3}. Thus, $\beta_{S}\geq d^{\totheright}+d^{\totheleft}\geq\kappa_0$. Therefore, one can compute $\kappa_0$ by looking only at $\beta_{S}$ for subsets $S$ of $\Z$ with diameter no larger than $M$ (note that this includes $[0,m_0-1]$, which has $\beta_{[0,m_0-1]}=d^{\totheright}+d^{\totheleft}$). By shift invariance, one can in fact look only at subsets of $[0,M]$. Since there are only finitely many such sets, the infimum in the definition of $\kappa_0$ is actually a minimum. Since a suitable $M$ can be easily calculated from $L$, $R$, and the $\alpha_i$, finding such an $M$ and then examining all strongly connected subsets of $\Z$ with leftmost point 0 and diameter $\leq M$ gives us an algorithm to find $\kappa_0$. 
\end{proof}
\setcounter{exmp}{\value{ExampleSave}}

We are now able to prove Theorem \ref{thm:gamma}.
\begin{proof}[Proof of Theorem \ref{thm:gamma}]\hypertarget{proof:gamma}
(a) $\Rightarrow$ (b)
Suppose $\kappa_0\leq s$. By Proposition \ref{prop:GammaAttained}, this means there is a finite, strongly connected set $S$ of vertices in $\Z$ such that $\beta_S\leq s$. By translation invariance, we may assume $0$ is the rightmost point of $S$. Let $-M$ be large enough that $S\subseteq[-M,0]$. By collapsing all vertices not in $[-M,0]$ into a sink and then arguing as in Claim \ref{claim:ComparingGraphs}, we may apply Tournier's lemma along with the amalgamation property to see that $E_{\G}\left[E_{\omega}^0\left[N_0^{[-M,0]}\right]^s\right]=\infty$.

(b) $\Rightarrow$ (c)
Immediate, since $N_0^{[-M,0]}\leq N_0^{[-M,M]}$.

(c) $\Rightarrow$ (a)
Suppose $E_{\G}\left[E_{\omega}^0\left[N_0^{[-M,M]}\right]^s\right]=\infty$. Again, collapsing all vertices not in $[-M,M]$ into a single sink, we may apply Tournier's lemma to see that there is a strongly connected set $S\subseteq[-M,M]$ such that $\beta_S\leq s$. Hence $\kappa_0\leq s$.
\end{proof}

The method for finding $\kappa_0$ given in the proof of Proposition \ref{prop:GammaAttained} requires knowledge of the $\alpha_i$, and the number of sets to examine grows exponentially in the smallest positive $\alpha_i$. We prove in Appendix \ref{append:kappa} that given only $L$, $R$, and the $i$ for which $\alpha_i>0$, $\kappa_0$ can be expressed as a minimum of finitely many positive integer combinations of the $\alpha_i$. If one has this formula, then one may easily compute $\kappa_0$ for any specific values of the $\alpha_i$.

\begin{prop}\label{prop:GammaAttainedweakstrong}
Given $L$, $R$, and the $i$ for which $\alpha_i>0$, $\kappa_0$ is an elementary function (a minimum of finitely many positive integer combinations) of the $\alpha_i$.
\end{prop}

Notice that Proposition \ref{prop:GammaAttainedweakstrong} would have sufficed in place of Proposition \ref{prop:GammaAttained} to show that the infimum in the definition of $\kappa_0$ is a minimum, which is enough to prove Theorem \ref{thm:gamma}. Proposition \ref{prop:GammaAttainedweakstrong} is stronger than Proposition \ref{prop:GammaAttained} in that it shows, given the structure of the graph $\G$, that there is an elementary formula for $\kappa_0$ that holds for all possible choices of $\alpha_i$, whereas Proposition \ref{prop:GammaAttained} gives $\kappa_0$ as a minimum of $\beta_S$ for $S$ in a set that is finite, but whose size depends on the size of the $\alpha_i$. On the other hand, Proposition \ref{prop:GammaAttainedweakstrong} is weaker than Proposition \ref{prop:GammaAttained} in that it does not provide an explicit algorithm for finding $\kappa_0$. This is because we do not know of a general way to explicitly find the finite set $\allowable^*$ given in the proof. Nevertheless, we do give examples in Appendix \ref{append:kappa} where we are able to find this set and thus give $\kappa_0$ as a minimum of finitely many sums. We leave it as an open question to find a general algorithm to do this.

\subsection{Large-scale backtracking: The parameter \texorpdfstring{$\kappa_1$}{kappa1}}\label{subsec:backtracking}

We have seen that the parameter $\kappa_0$ controls moments of the quenched expected time a walk spends at 0 before exiting a finite set. We will now show that in a similar way, $\kappa_1$ controls moments of backward traversals of arbitrarily large stretches of the graph. 

In our proof of Theorem \ref{thm:transiencecriterion}, we used the graphs $\G_M$, finite graphs that looked like $\G$ except near endpoints. Here, we consider these along with a ``limiting graph" that is half infinite. Let $\G_+$ be a graph with vertex set $[0,\infty)$. The graph $\G_+$ contains all the same edges between vertices to the right of 0 with the same weights as $\G$. For vertices $1\leq i\leq L$, there is an edge from $i$ to 0 with weight $\sum_{j={1-L}}^0\alpha_{j-i}$. And to each vertex $1\leq j\leq R$ is added an edge from 0 with weight $\sum_{i=j}^R\alpha_i$.

\begin{figure}
    \centering
    \includegraphics[width=6.5in]{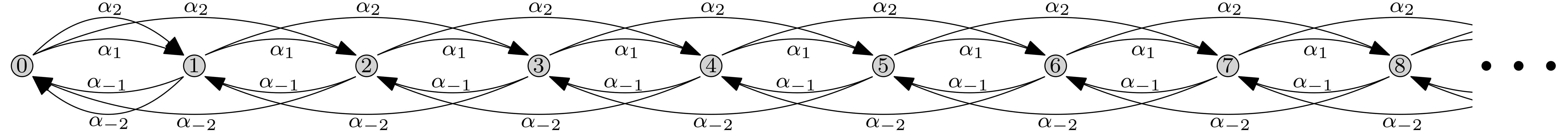}
    \caption{The graph $\G_+$}
    \label{fig:graphG+}
\end{figure}

The graph $\G_+$ has zero divergence at all sites except 0, where the divergence is $d^{\totheright}-d^{\totheleft}=\kappa_1$. Thus, in a sense there is a ``net flow" of strength $\kappa_1$ from 0 to infinity, and to motivate the following lemma, the reader may imagine an edge ``from $\infty$ to 0" with weight $\kappa_1$. In some sense, the following lemma extends Corollary \ref{cor:loopreversal} (1) to this infinite graph. One can prove it using a comparison between $\G$ and $\G_M$.

\begin{lem}[{\cite[Theorem 2]{Tournier2015}}]\label{lem:InfiniteTimeReversal}
Under $P_{\G_+}$, $P_{\omega}^0(\tilde{H}_0=\infty)\sim \text{\em{Beta}}(\kappa_1,d^{\totheleft})$.
\end{lem}



We will use Lemma \ref{lem:InfiniteTimeReversal} to prove Theorem \ref{thm:TFAE2} in two separate propositions. Recall that for a given walk ${\bf X}$ and integers $x<y$, the quantity $N_{x,y}=N_{x,y}({\bf X})$ is defined as the number of trips from $y$ to $x$.

\begin{prop}\label{prop:InfiniteTraps}
Suppose $s\geq \kappa_1$. Then the following hold:
\begin{enumerate}
    \item $E_{\G}[E_{\omega}^0[N_0]^s]=\infty$.
    \item For all $x<y\in\Z$, $E_{\G}[E_{\omega}^0[N_{x,y}]^s]=\infty$.
\end{enumerate}
\end{prop}

Like Theorem \ref{thm:gamma}, this proposition not only gives a sufficient condition for zero speed, but gives the reason for zero speed. If $\kappa_0\leq1$, the speed is zero because it takes a long time for the walk to exit small traps; on the other hand, if $\kappa_1\leq1$, the speed is zero because the walk traverses large regions of $\Z$ many times.

\setcounter{ExampleSave}{\value{exmp}}
\setcounter{exmp}{\value{thm}}
\setcounter{claim}{0}
\newpage
\begin{proof}\hspace{0pt}

\begin{center}{\bf Outline}\end{center}
The philosophy of the proof
is that the components of a Dirichlet random vector become more and more independent as their values become small.
If $(Y_1,Y_2,Y_3,Y_4)$ is a Dirichlet vector with parameters $(a,b,c,d)$, and $(X_1,X_2)$ and $(X_3,X_4)$ are independent Dirichlet vectors with parameters $(a,y)$ and $(c,z)$, respectively, then $P(X_1<\tepsilon,X_3<\tdelta)\asymp P(Y_1<\tepsilon,Y_3<\tdelta)$ as $(\tepsilon,\tdelta)\to(0,0)$. In other words, although these probabilities do not necessarily become approximately equal (even if $y=b$ and $z=c$), they are bounded by constant multiples of each other. For each $i=1,\ldots,R$, the Dirichlet weight entering $i$ from $[1-R,0]$ in $\G$ is the same as the Dirichlet weight entering $i$ from $0$ in $\G_+$. The goal of this proof is to exploit the comparability of small-value probabilities and perform a coupling between $(\omega(0,i))_{i=1}^R$ under $P_{\G_+}$ and $\left(\sum_{j=1-R}^0\omega(j,i)\right)_{i=1}^R$ under $P_{\G}$.

Actually, we will couple vectors that distinguish between different edges to the same vertex. Let $E_0:=\{(i,j):1-R\leq i\leq0,1\leq j\leq i+R,\alpha_{j-i}>0\}$ be the set of right-oriented edges in $\G$ that originate from or cross 0, and for every $e=(i,j)\in E_0$, let $\alpha_e=\alpha_{j-i}$. We consider random vectors ${\bf Z}=\left(Z_e\right)_{e\in E_0}$ and ${\bf Y}=\left(Y_e\right)_{e\in E_0}$, and a measure $P'$ such that the distribution of $\left(Z_e\right)_{e\in E_0}$ under $P'$ is the distribution of $\left(\omega(i,j)\right)_{(i,j)\in E_0}$ under $P_{\G}$, and such that ${\bf Y}$ is a Dirichlet random vector with parameters $(\alpha_e)_{e\in E_0}$. The amalgamation property implies that $\left(\sum_{i=1-R}^0Y_{(i,j)}\right)_{j=1}^R$ is distributed like $(\omega(0,j)_{j=1}^R$ under $P_{\G_+}$. The idea is to define a coupling event $K$, independent of ${\bf Y}$ and with positive probability, on which $Z_e\leq Y_e$ for all $e$. We do not quite accomplish this, but we come close enough that we are able to use ${\bf Z}$ and ${\bf Y}$ to construct random environments $\omega_1$ and $\omega_2$, drawn respectively according to $P_{\G}$ and $P_{\G_+}$, such that on $K$, $\sum_{i=1-R}^0P_{\omega_1}^i(\tilde{H}_{[1-R,0]}=\infty)$ is bounded above by a constant multiple of $P_{\omega_2}^0(\tilde{H}_0=\infty)$. From Lemma \ref{lem:InfiniteTimeReversal} we get that $E_{\G_+}\left[\frac{1}{P_{\omega}^0(\tilde{H}_0=\infty)}\right]=\infty$, so our coupling gives us $E'\left[\frac{1}{\sum_{i=1-R}^0P_{\omega_1}^i(\tilde{H}_{[1-R,0]}=\infty)}\right]=\infty$. This is enough to give us $\Ea_{\G}^0[N_0]=\infty$, and more careful analysis yields $\Ea_{\G}^0[N_{x,y}]=\infty$.

\begin{center}{\bf Groundwork for the coupling}\end{center}
Suppose the edges in $E_0$ are enumerated as $e_1,\ldots,e_k$ in some way. (In fact, we will enumerate them in a random way, yet to be described, but for now assume the enumeration is fixed.) We have said how ${\bf Z}$ and ${\bf Y}$ will be distributed under $P'$. By the amalgamation property, $Z_{e_1}$ and $Y_{e_1}$ are both beta random variables whose first parameter is $\alpha_{e_1}$. Their second parameters may differ, but by \eqref{eqn:moments} we nonetheless have $P'(Z_{e_1}<\tepsilon)\asymp \tepsilon^{\alpha_{e_1}} \asymp P'(Z_{e_1}<\tepsilon)$, where $f(\tepsilon)\asymp g(\tepsilon)$ means there exist positive constants $c,C$ such that $cg(\tepsilon)\leq f(\tepsilon)\leq Cg(\tepsilon)$ for all $\tepsilon\in[0,1]$.

Note that for $1\leq i\leq k-1$, $Y_{e_i}':=\frac{Y_{e_i}}{1-\sum_{j=1}^{i-1}Y_{e_j}}=\frac{Y_{e_i}}{\sum_{j=i}^kY_{e_j}}$ is a beta random variable, independent of $Y_{e_1},\ldots,Y_{e_{i-1}}$, and with first parameter $\alpha_{e_i}$ (this comes from the restriction property, along with the amalgamation property). Let $Y_{e_k}':=\frac{Y_{e_k}}{1-\sum_{j=1}^{k-1}Y_{e_j}}=1$. Likewise, for $1\leq i\leq k$, $Z_{e_i}':=\frac{Z_{e_i}}{1-\sum_{j=1}^{i-1}Z_{e_j}\one_{\{\underline{e_j}=\underline{e_i}\}}}$ is a beta random variable, independent of $Z_{e_1},\ldots,Z_{e_{i-1}}$, and with first parameter $\alpha_{e_i}$. (We do not have $Z_{e_k}=1$ a.s., because $Z_{e_k}$ corresponds to an edge from a vertex $i\leq0$ to a vertex $j>0$, and  there are still edges from $i$ to vertices to the left of 0.)  By \eqref{eqn:moments}, $P'(Z_{e_i}'<\tepsilon)\asymp \tepsilon^{\alpha_i}\asymp P'(Y_{e_i}'<\tepsilon)$, for $0\leq i\leq k-1$. Thus, there exists a constant $c$ such that for all $\tepsilon\in[0,1]$, and for all $i=1,\ldots,k-1$, we have $P'(Z_{e_k}'<\tepsilon)\geq cP'(Y_{e_k}'<\tepsilon)$. This $c$ may depend on the chosen permutation $(e_1,\ldots,e_k)$ of $E_0$, but there are only finitely many permutations, so we may assume $c$ is small enough to work for any of them. For each $1\leq i\leq k$, let $F_{Z_{e_i}'}$ be the cdf for $Z_{e_i}'$, and let $Q_{Z_{e_i}'}$ be the associated quantile function (since $F_{Z_{e_i}'}$ is continuous and strictly increasing on $[0,1]$, $Q_{Z_{e_i}'}$ is simply the inverse of $F_{Z_{e_i}'}$ restricted to the interval $[0,1]$). Similarly, for $1\leq i\leq k$, let $F_{Y_{e_i}'}$ and $Q_{Y_{e_i}'}$ be the cdf and quantile function for $Y_{e_i}'$.\footnote{Since $Y_{e_k}'$ is identically 1, its cdf is 0 to the left of 1 and 1 at 1, and its quantile function is identically 1. All other $Y_{e_i}'$ have continuous cdfs and quantile functions.} Choose $\ell$ to be an integer large enough that $\frac{1}{\ell}<c$. Then
\begin{equation}\label{eqn:1515}
    F_{Z_{e_i}'}\geq \frac{1}{\ell}F_{Y_{e_i}'},\quad 1\leq i\leq k.
\end{equation}

Now given ${\bf Y}'=(Y_{e_1}',\ldots,Y_{e_k}')$ and ${\bf Z}'=(Z_{e_1}',\ldots,Z_{e_k}')$, we can recover ${\bf Y}$ and ${\bf Z}$. First, $Y_{e_1}=Y_{e_1}'$ and $Z_{e_1}=Z_{e_1}'$, and then if $Y_j$ and $Z_j$ are known for $1\leq j\leq i$, then the formulas for $Y_{e_i}'$ and $Z_{e_i}'$ can be used to find $Y_{e_i}$ and $Z_{e_i}$. Therefore, one way to generate the vector ${\bf Y}$ is to generate independent beta random variables $Y_{e_1}',Y_{e_2}',\ldots,Y_{e_{k-1}}'$ with appropriate parameters (which depend on the permutation $(e_1,\ldots,e_k)$) and then use these to recover ${\bf Y}$. Similarly, ${\bf Z}$ can be generated by means of independent beta random variables $Z_{e_1}',\ldots,Z_{e_k}'$. Under this method, the chosen permutation $(e_1,\ldots,e_k)$ affects the parameters for the $Y_{e_i}'$ and $Z_{e_i}'$, as well as the order in which they are put together, but the distributions of ${\bf Y}$ and ${\bf Z}$, are the same regardless of the chosen permutation.

\begin{center}{\bf The coupling}\end{center}
 Our probability space is $[0,1]^{k}\times\{0,\ldots,\ell-1\}^k\times\Omega_{\Z}$. Let $P'$ be the product measure whose marginals on $\Omega_{\Z}$ are equal to $P_{\G}$, and whose marginals on $[0,1]^{k}\times\{0,\ldots,\ell-1\}^k$ are uniform. An element of our probability space will be of the form
\begin{equation*}\left(U_{1},\ldots,U_{k},V_1,\ldots,V_k,\omega\right),\end{equation*}
where the $U_i$ take values on $[0,1]$, the $V_i$ take integer values from 0 to $\ell-1$, and $\omega$ can be any environment on $\Z$. Define the function $W_i:=\frac{V_i+U_i}{\ell}$. Then the $W_i$ are i.i.d. uniform $[0,1]$ under $P'$.

Recall that $\omega^{>0}$ is the environment $\omega$ to the right of 0; that is, if $\omega=(\omega(a,b))_{a,b\in\Z}$, then $\omega^{>0}=(\omega(a,b))_{a,b\in\Z,a>0}$. The values of $P_{\omega}^i(H_0=\infty)$ are determined by $\omega^{>0}$ for $1\leq i\leq R$. For a given $\omega^{>0}$, let $(e_1,\ldots,e_k)$ be a permutation of $E_0$ such that

\begin{equation}\label{eqn:1532}
i<j \quad \Rightarrow \quad P_{\omega}^{\overline{e_j}}(H_0=\infty)\leq P_{\omega}^{\overline{e_i}}(H_0=\infty).
\end{equation}
To get such an arrangement, sort vertices $1\leq j\leq R$ in order of $P_{\omega}^j(H_0=\infty)$, then sort the edges $(i,j)\in E_0$ primarily according to the rank of $j$ and secondarily according to the value of $i$.

We can now use uniform random variables along with quantile functions to get
${\bf Z}'=(Z_{e_1}',\ldots,Z_{e_k}')$ and ${\bf Y}'=(Y_{e_1}',\ldots,Y_{e_k}')$. 
Letting $Z_{e_i}'=Q_{Z_{e_i}'}(W_1)$ gives us the desired distribution for each $Z_{e_i}'$ under $P'$, and letting $Y_{e_i}'=Q_{Y_{e_i}'}(U_1)$ gives us the desired distribution for each $Y_{e_i}'$. This gives us ${\bf Y}'$ and ${\bf Z}'$, from which we may then recover ${\bf Z}$ and ${\bf Y}$.
Even though the specific permutation $(e_1,\ldots, e_k)$ is used along with the $U_i$ in defining ${\bf Y}$, the distribution of ${\bf Y}$ is the same for any fixed permutation, and so ${\bf Y}$ is independent of $\omega$. Similarly, ${\bf Z}$ is also independent of $\omega$.

Define the {\em coupling event} $K$ to be the event that $V_1=\ldots=V_k=0$, the walk is transient to the right, $P_{\omega}$--a.s., and $\omega(i,j)>0$ iff $\alpha_{j-i}>0$ for all $i,j\in\Z$. Because these last two conditions each have $P'$ probability 1, $K$ is independent of $\omega$ as well as ${\bf Y}$, and has positive probability $P'(K)=\left(\frac{1}{\ell}\right)^k$. On $K$, $W_i=\frac{1}{\ell}U_i$ for $1\leq i\leq R$. Let $t\in[0,1]$. Then $Q_{Z_{e_i}'}(\frac{1}{\ell}t)$ is the unique $x$ such that $F_{Z_{e_i}'}(x)=\frac1{\ell}t$.
For this $x$, applying \eqref{eqn:1515} gives us
$\frac{1}{\ell}t\geq\frac{1}{\ell}F_{Y_{e_i}'}(x)$, or $t\geq F_{Y_{e_i}'}(x)$. Applying the increasing function $Q_{Y_{e_i}'}$ to both sides, we get  $Q_{Y_{e_i}'}(t)\geq x=Q_{Z_{e_i}'}(\frac{1}{\ell}t)$. This is true for all $t\in[0,1]$, so on the event $K$, we have
\begin{equation}\label{eqn:1544}
Z_{e_i}'\leq Y_{e_i}',\quad\quad i=1,\ldots,k.
\end{equation}

We now describe how to use ${\bf Z}$, ${\bf Y}$, and $\omega$ to create environments $\omega_1$ and $\omega_2$, drawn according to $P_{\G}$ and $P_{\G_+}$, respectively. For $i,j>0$, let $\omega_1(i,j)=\omega_2(i,j)=\omega(i,j)$. We also let $\omega_1(i,j)=\omega(i,j)$ for $i>0$ and any $j\leq0$. But for $\omega_2$, transition probabilities from positive to negative vertices are ``collapsed" to 0. That is, for $i>0$, we let $\omega_2(i,0)=\sum_{j\leq0}\omega(i,j)$. By the amalgamation property, one can check that transition probabilities at sites greater than $0$ are drawn according to $P_{\G_+}$. Moreover, for $i\geq1$, we have

\begin{equation}\label{eqn:1547}
    P_{\omega_1}^i(H_{\leq0}=\infty)=P_{\omega_2}^i(H_0=\infty)
\end{equation}

Then let $\omega_2(0,j)=\sum_{i=1-R}^0Y_{(i,j)}$, $j=1,\ldots,R$. By the distribution of ${\bf Y}$ and the fact that it is independent of $\omega$, transition probabilities at sites greater than or equal to 0 of $\omega_2$ are drawn according to $P_{\G_+}$ (sites less than 0 don't matter, so for example we can let $\omega_2(i,j)=\one_{\{j=i\}}$ whenever $i<0$).

For $i\notin[1-R,0]$, and for all $j$, let $\omega_1(i,j)=\omega(i,j)$. For $i\in[1-R,0]$ and for all $j>0$ with $(i,j)\in E_0$, let $\omega_1(i,j)=Z_{(i,j)}$. For all $i\in[1-R,0]$ and $j\leq0$, we keep $\omega_1(i,j)$ the same as $\omega(i,j)$, but scaled to ensure that $\sum_j\omega_1(i,j)=1$. That is,

\begin{equation*}\omega_1(i,j)=\left(1-\sum_{r>0}Z_{(i,r)}\right)\frac{\omega(i,j)}{\sum_{r\leq0}\omega(i,r)}\end{equation*}

Notice that under this definition, $\left(\frac{\omega_1(i,j)}{\sum_{r\leq0}\omega_1(i,r)}\right)_{r\leq0}=\left(\frac{\omega(i,j)}{\sum_{r\leq0}\omega(i,r)}\right)_{r\leq0}$. By the restriction property, $\left(\frac{\omega(i,j)}{\sum_{r\leq0}\omega(i,r)}\right)_{r\leq0}$ is independent of $\left(\omega(i,r)\right)_{r>0}$ under $P_{\G}$, and by the way we have defined $\omega_1$, $\left(\frac{\omega_1(i,j)}{\sum_{r\leq0}\omega_1(i,r)}\right)_{r\leq0}$ is independent of $(Z_{(i,r)})_{r>0}$ under $P'$. Therefore, since  $(Z_{(i,r)})_{r>0}\overset{(L)}{=}\left(\omega(i,r)\right)_{r>0}$ by construction, we have $\left(\omega_1(i,j)\right)_j\overset{(L)}{=}\left(\omega(i,j)\right)_j$ for each $i$, and the transition probability vectors are all independent. Hence the law of $\omega_1$ is $P_{\G}$.
\newpage
\begin{center}{\bf Comparing sums of probabilities}\end{center}

The goal for this part of the proof is to show that on the event $K$, 
\begin{equation}\label{eqn:1557}
\sum_{i=1}^{k}Z_{e_i}P_{\omega_1}^{\overline{e_i}}(H_{\leq0}=\infty)
\leq
C\sum_{i=1}^{k}Y_{e_i}P_{\omega_2}^{\overline{e_i}}(H_0=\infty)
\end{equation}
for some deterministic constant $C$. Note that the sum on the right side of \eqref{eqn:1557} is equal to $P_{\omega_2}^0(\tilde{H}_0=\infty)$. 

By \eqref{eqn:1547}, we could get \eqref{eqn:1557} by showing that $Z_{e_i}\leq CY_{e_i}$ for all $i$. However, achieving this precisely would require a more elaborate coupling. The difficulty is that if, for example, $Y_{e_1}$ is very close to 1, all other $Y_{e_i}$ are forced to be very small, whereas some of the $Z_{e_i}$ are independent of $Z_{e_1}$. Our specific ordering of the edges $e_1,\ldots,e_k$ allows us to get around this difficulty. 

Let $r$ be the smallest integer in $\{1,\ldots,R\}$ such that $\sum_{i=1}^{r}Y_{e_i}>\frac12$ ($r$ is random). For $1\leq i\leq r$, on the event $K$ we have
\begin{equation*}
Z_{e_i}<\frac{Z_{e_i}}{1-\sum_{j=1}^{i-1}Z_{e_j}\one_{\{\underline{e_j}=\underline{e}_i\}}}
\leq
\frac{Y_{e_i}}{1-\sum_{j=1}^{i-1}Y_{e_j}}
\leq
2Y_{e_i}.
\end{equation*}
(The middle terms are the definitions of $Z_{e_i}'$ and $Y_{e_i}'$, respectively.) We now have
\begin{align}
\notag
    \sum_{i=1}^kZ_{e_i}P_{\omega_1}^{\overline{e}_i}(H_{\leq0}=\infty)
    &=\sum_{i=1}^rZ_{e_i}P_{\omega_1}^{\overline{e}_i}(H_{\leq0}=\infty) + \sum_{i=r+1}^kZ_{e_i}P_{\omega_1}^{\overline{e}_i}(H_{\leq0}=\infty)
    \\
    \notag
    &\leq \sum_{i=1}^r2Y_{e_i}P_{\omega_2}^{\overline{e}_i}(H_0=\infty) + \sum_{i=r+1}^kP_{\omega_2}^{\overline{e}_i}(H_0=\infty)
    \\
    \label{eqn:31}
    &\leq \sum_{i=1}^r2Y_{e_i}P_{\omega_2}^{\overline{e}_i}(H_0=\infty) + kP_{\omega_2}^{\overline{e}_r}(H_0=\infty),
\end{align}
where, for the last line, we used the fact that $P_{\omega}^{\overline{e}_i}(H_0=\infty)$ is non-increasing in $i$ by \eqref{eqn:1532}. We want to combine the two terms from \eqref{eqn:31} into one. To do this, we note
\begin{align*}
    \sum_{i=1}^r2Y_{e_i}P_{\omega_2}^{\overline{e}_i}(H_0=\infty)
    &\geq \sum_{i=1}^r2Y_{e_i}P_{\omega_2}^{\overline{e}_r}(H_0=\infty)
    \\
    &=2P_{\omega_2}^{\overline{e}_r}(H_0=\infty)\sum_{i=1}^rY_{e_i}
    \\
    &\geq P_{\omega_2}^{\overline{e}_r}(H_0=\infty),
\end{align*}
where we used the same non-increasing property for the first line, and the definition of $r$ in the last line. Applying this to \eqref{eqn:31} gives us
\begin{align}
    \notag
     \sum_{i=1}^kZ_{e_i}P_{\omega_1}^{\overline{e}_i}(H_{\leq0}=\infty)
     &\leq \sum_{i=1}^r2Y_{e_i}P_{\omega_2}^{\overline{e}_i}(H_0=\infty) + k \sum_{i=1}^r2Y_{e_i}P_{\omega_2}^{\overline{e}_i}(H_0=\infty)
     \\
     \notag
     &= 2(k+1)\sum_{i=1}^rY_{e_i}P_{\omega_2}^{\overline{e}_i}(H_0=\infty)
     \\
     \label{eqn:MainInequalityforTransienceProof}
     &\leq 2(k+1)\sum_{i=1}^kY_{e_i}P_{\omega_2}^{\overline{e}_i}(H_0=\infty).
\end{align}
This is exactly \eqref{eqn:1557}. 

\begin{center}{\bf Comparing expectations}\end{center}

We consider the probability in $\omega_1$, starting from a point $a$ in $(-\infty,0]$, of never hitting the set $[1-R,0]$ at a positive time. If $\omega_1$ is transient to the right and jumps to the right are bounded by $R$, then the only way for this to occur is for $a$ to be in $[1-R,0]$, for the first step to be to the right of 0, and then for the walk to never again hit a site to the left of 0. Thus,
on the coupling event $K$,
\begin{align}
\notag
\max_{a\leq0}P_{\omega_1}^a(\tilde{H}_{[1-R,0]}=\infty)
&\leq \sum_{i=1-R}^0P_{\omega_1}^i(\tilde{H}_{[1-R,0]}=\infty)
\\
\notag
&=\sum_{i=1-R}^0\sum_{j=1}^{R}Z_{(i,j)}P_{\omega_1}^j(H_{\leq0}=\infty)
\\
\notag
&\leq 2(k+1)\sum_{i=1-R}^0\sum_{j=1}^{R}Y_{(i,j)}P_{\omega_2}^j(H_0=\infty)
\\
\label{eqn:1682}
&=2(k+1)P_{\omega_2}^0(\tilde{H}_0=\infty).
\end{align}

It is straightforward to check by induction that for all $n\geq1$, $a\leq0$, \begin{equation}\label{eqn:1685}P_{\omega_1}^a(N_{[1-R,0]}\geq n)\geq \min_{1-R\leq i\leq0}P_{\omega_1}^i(\tilde{H}_{[1-R,0]}<\infty)^{n-1}.
\end{equation}
Summing over all $n\geq1$ in \eqref{eqn:1685} and applying \eqref{eqn:1682}, we get on the coupling event $K$,

\begin{align}
    \notag
    E_{\omega_1}^a[N_{[1-R,0]}]
    &=\sum_{n=1}^{\infty}P_{\omega_1}^a(N_{[1-R,0]}\geq n)
    \\
    \notag
    &\geq \sum_{n=1}^{\infty}\min_{1-R\leq i\leq0}P_{\omega_1}^i(\tilde{H}_{[1-R,0]}<\infty)^{n-1}
    \\
    \notag
    &=\frac{1}{\max_{1-R\leq i\leq0}P_{\omega_1}^i(\tilde{H}_{[1-R,0]}=\infty)}
    \\
    \notag
    &\geq\frac{1}{2(k+1)P_{\omega_2}^0(\tilde{H}_0=\infty)}
    \\
    \label{eqn:Comparison}
    &=\frac{1}{2(k+1)}E_{\omega_2}^0[N_0].
\end{align}

Since the event $K$ is independent of $\omega_2$, we conclude that 
 \begin{align*}
 E'\left[E_{\omega_1}^0[N_{[1-R,0]}]^s\right]
 &\geq  E'\left[E_{\omega_1}^0[N_{[1-R,0]}]^s\one_K\right]
 \\
 &\geq E'\left[\frac{1}{2^s(k+1)^s}E_{\omega_2}^0[N_0]^s\one_K\right]
 \\
 &=\frac{1}{2^s(k+1)^s}E'\left[E_{\omega_2}^0[N_0]^s\right]P'(K)
 \\
 &=\frac{1}{\ell^k}\frac{1}{2^s(k+1)^s}E'\left[E_{\omega_2}^0[N_0]^s\right].
 \end{align*}
 By the way $\omega_1$ and $\omega_2$ are distributed, this means
 \begin{align}
 \notag
     E_{\G}\left[E_{\omega}^0[N_{[1-R,0]}]^s\right]
     &\geq \frac{1}{\ell^k}\frac{1}{2^s(k+1)^s}E_{\G_+}\left[E_{\omega}^0[N_0]^s\right]
     \\
     \label{eqn:1695}
     &=\frac{1}{\ell^k}\frac{1}{2^s(k+1)^s}E_{\G_+}\left[\frac{1}{P_{\omega}^0(\tilde{H}_0=\infty)^s}\right]
 \end{align} 
If $s\geq\kappa_1$, then the right side of \eqref{eqn:1695} is infinite by Proposition \ref{lem:InfiniteTimeReversal} and \eqref{eqn:moments}, so we have
\begin{align*}
    \infty&=\E_{\G}[E_{\omega}^0[N_{[1-R,0]}]^s]
    \\
    &=E_{\G}\left[\left(\sum_{i=1-R}^0E_{\omega}^0[N_i]\right)^s\right]
    \\
    &\leq E_{\G}\left[\left(\sum_{i=1-R}^0E_{\omega}^i[N_i]\right)^s\right]
    \\
    &\leq R^sE_{\G}\left[\left(\max_{1-R\leq i\leq 0}E_{\omega}^i[N_i]\right)^s\right]
    \\
    &=R^s\sum_{i=1-R}^0E_{\G}\left[\left(E_{\omega}^i[N_i]\right)^s\right]
    \\
    &=R^{s+1}E_{\G}\left[\left(E_{\omega}^0[N_0]\right)^s\right].
\end{align*}
This proves the first part of the proposition. At this point, the reader interested only in characterizing ballisticity may skip the remainder of the proof, and may also skip Proposition \ref{prop:FewOscillations}, going straight to Section \ref{subsec:ballisticity}. However, the next part of this proof and Proposition \ref{prop:FewOscillations} together provide an important insight into the behavior of the walk: namely, that the walk is expected to oscillate back and forth between any two points infinitely many times precisely when $\kappa_1\leq1$.

\begin{center}{\bf Arbitrarily large backtracking}\end{center}

We now want to prove the second part of the proposition, which strengthens our result to show that the expected number of oscillations between any two points is infinite. We do this via the following claim.

\begin{claim}\label{claim:backtrackingclaim}
 For any $a\leq 0$ and $x<y\leq0$ we have
\begin{equation*}
    E_{\G}\left[E_{\omega}^a[N_{x,y}]^s\given \omega^{\leq -R} \right]=\infty,~P_{\G}\text{--a.s.}
\end{equation*}
\end{claim}

Assume for now that the claim is true. Taking expectations on both sides gives us 
\begin{equation}\label{eqn:1735}
    E_{\G}\left[E_{\omega}^a[N_{x,y}]^s\right]=\infty.
\end{equation}
Let $x<y\in\Z$. If $y\leq0$, then letting $a=0$ in \eqref{eqn:1735} gives us $E_{\G}\left[E_{\omega}^0[N_{x,y}]^s\right]=\infty$, which is exactly what we needed to show for the second part of the proposition. If $y>0$, then \eqref{eqn:1735} gives us $E_{\G}\left[E_{\omega}^{-y}[N_{x-y,0}]^s\right]=\infty$, and then the translation invariance of $\G$ gives us $E_{\G}\left[E_{\omega}^0[N_{x,y}]^s\right]=\infty$.

It remains, then, to prove the claim. Under $P'$, $\omega_1$ is drawn according to $P_{\G}$. Since $\omega_1$ agrees with $\omega$ on $(-\infty,-R]$, our claim is equivalent to the statement $E'\left[E_{\omega_1}^a[N_{x,y}]^s\given \omega^{\leq -R} \right]=\infty$, $P'$--a.s. And since $\sigma(\omega^{\leq -R})$ is coarser than $\sigma(\omega^{\leq0})$, it suffices to show that 
\begin{equation}\label{eqn:1734}
E'\left[E_{\omega_1}^a[N_{x,y}]^s\given\omega^{\leq0} \right]=\infty,\quad P'\text{--a.s.}
\end{equation}

We first show that for $i\in[1-R,0]$ and $j\leq0$, there is a constant $C>0$ such that on the event $K$,
\begin{equation}\label{eqn:1741}
\frac{\omega_1(i,j)}{\omega(i,j)}\geq C.
\end{equation}
Recall that for $i\in[1-R,0]$ and $j\leq0$, we have defined
\begin{equation*}
\omega_1(i,j)=\left(1-\sum_{r>0}Z_{(i,r)}\right)\frac{\omega(i,j)}{\sum_{r\leq0}\omega(i,r)}.
\end{equation*}
Now for $i\in[1-R,0]$, we know that $\sum_{r>0}Z_{(i,r)}$ can be arbitrarily close to 1. However, we assert it is bounded away from 1 on the coupling event $K$, where necessarily $W_1,W_2,\ldots,W_k<\frac{1}{\ell}$. In fact, we assert that on this event, $\sum_{r>0}Z_{(i,r)}$ is maximized when $W_1=\cdots=W_k=\frac{1}{\ell}$. One can check by induction that for a given $1-R\leq i\leq 0$,
\begin{equation}\label{eqn:1743}
\sum_{r>0}Z_{(i,r)}=1-\prod_{r>0}(1-Z_{(i,r)}').
\end{equation}
Since all the $Z_{(i,r)}'$ are independent, the right side of \eqref{eqn:1743} is maximized when they are all as large as possible. On $K$, the largest they can get is when all the $W_i$ are equal to $\frac1{\ell}$, and this yields a value less than 1, proving our assertion and giving us \eqref{eqn:1741}.

Let $-M<0$. For $1-R\leq i\leq 0$, consider the ($i,x,y$) {\em excursion event} $\mathcal{E}_{i,x,y}$ where:
\begin{itemize}
    \item $X_0=i$
    \item If $i\neq y$, the walk hits $y$ before returning to $i$ or leaving $(-\infty,0]$. 
    \item After $H_y$, the walk hits $x$ without hitting $i$ more than once in between and without leaving $(-\infty,0]$.\footnote{It is necessary to allow for the walk to hit $i$ once on the way to $x$ in case $b>i$ and the only way to reach $x$ from $y$ is through $i$. For all other cases, we could require that the walk avoid $i$ between hitting $y$ and hitting $x$, but to avoid treating this case separately, we allow one visit to $i$ in all cases, at the cost of a factor of 2.}
\end{itemize}
We say an $(i,x,y)$ excursion event starts at time $n$ if $(X_n,X_{n+1},\ldots)\in\mathcal{E}_{i,x,y}$. Then the number of such excursion events for any $i$ is no more than $2N_{x,y}$. This is because each trip from $y$ to $x$ can count toward at most two excursion events due to the requirement that there be only one visit to $i$ in between visiting $y$ and $x$.\footnote{For example, if the walk goes from $i$ to $y$ to $i$ again to $y$ again and then to $x$, then an excursion event started at each of the times the walk was at $i$, but there was only one trip from to $y$ to $x$.}

Fix $\omega^{\leq0}$. For any $i\leq0$, the probability under $P_{\omega}^i$ of any finite path that stays within $(-\infty,0]$ is fixed. On the event $K$ (which is independent of $\omega^{\leq0}$ and therefore still has probability $\frac{1}{\ell^k}$ conditioned on $\omega^{\leq0}$), the probability under $P_{\omega_1}^i$ of such a path is bounded from below due to \eqref{eqn:1741}. 
Therefore, on the event $K$, there exists a positive constant $c=c(\omega^{\leq0})$ such that on $K$, 
\begin{equation}\label{eqn:1764}
 \min_{1-R\leq i\leq0}P_{\omega_1}^i(\mathcal{E}_{i,x,y})>c.   
\end{equation}
(For each $i$ consider a particular finite path that achieves $\mathcal{E}_{i,x,y}$, take $c_i$ to be a lower bound for the probability under $P_{\omega_1}^i$ of taking that path, then take $c$ to be the minimum of the $c_i$.)

We have from \eqref{eqn:Comparison} that on $K$, $E_{\omega_1}^a[N_{[1-R,0]}]\geq\frac{1}{2(k+1)}E_{\omega_2}^0[N_0]$. Taking conditional expectations, we almost surely have 
\begin{align*}
E'\left[E_{\omega_1}^a[N_{[1-R,0]}]^s\given\omega^{\leq0} \right]
&\geq
E'\left[E_{\omega_1}^a[N_{[1-R,0]}]^s\one_K\given\omega^{\leq0} \right]
\\
&\geq
\frac{1}{2^s(k+1)^s} E'\left[E_{\omega_2}^0[N_0]^s\one_K\given\omega^{\leq0} \right]
\\
&=\frac{1}{2^s(k+1)^s}E'\left[E_{\omega_2}^0[N_0]^s \right]P'(K)
\\
&=\infty,
\end{align*}
where the first equality comes from the fact that $\omega^{\leq0}$, $K$, and $\omega_2$ are all independent. Now with probability 1,
\begin{align*}
\infty&=E'\left[E_{\omega_1}^a[N_{[1-R,0]}]^s\given\omega^{\leq0} \right]
\\
&=
E'\left[\left(\sum_{n=1}^{\infty}P_{\omega_1}^a(X_n\in[1-R,0])\right)^s\given\omega^{\leq0}\right]
\\
&=
E'\left[\left(\sum_{i=1-R}^0\sum_{n=1}^{\infty}P_{\omega_1}^a(X_n=i)\right)^s\given\omega^{\leq0}\right]
\end{align*}

Multiplying by $c^s$, where $c$ is the constant from \eqref{eqn:1764}, will not change the fact that the expression is infinite. And since $c^s$ depends only on $\omega^{\leq0}$, it may be pulled inside of an expectation conditioned on $\omega^{\leq0}$. Therefore we almost surely have
\begin{align*}
\infty&=E'\left[\left(\sum_{i=1-R}^0\sum_{n=1}^{\infty}cP_{\omega}^a(X_n=i)\right)^s\given\omega^{\leq0}\right]
\\
&\leq E'\left[\left(\sum_{i=1-R}^0\sum_{n=1}^{\infty}P_{\omega_1}^a(X_n=i)P_{\omega_1}^i(\mathcal{E}_{i,x,y})\right)^s\given\omega^{\leq0}\right]
\\
&= E'\left[\left(\sum_{i=1-R}^0E_{\omega_1}^a\left[\#\left\{n\in\N_0:(X_k)_{k=n}^{\infty}\in\mathcal{E}_{i,x,y}\right\}\right]\right)^s\given\omega^{\leq0}\right]
\\
&\leq E'\left[\left(\sum_{i=1-R}^0E_{\omega_1}^a[2N_{x,y}]\right)^s\given\omega^{\leq0}\right]
\\
&=
(2R)^sE'\left[E_{\omega_1}^a[N_{x,y}]^s\given\omega^{\leq0}\right]
\end{align*}
Thus, $E'\left[E_{\omega_1}^a[N_{x,y}]^s\given\omega^{\leq0}\right]=\infty$ with probability 1. This is \eqref{eqn:1734}, which suffices to prove our claim, and with it our theorem.
\end{proof}

\setcounter{exmp}{\value{ExampleSave}}

We now prove a slightly strengthened converse to Proposition \ref{prop:InfiniteTraps}. 

\begin{prop}\label{prop:FewOscillations}
If $s<\kappa_1$, then there is an $M\geq0$ such that for all $x,y\in\Z$ with $y-x\geq M$,  $E_{\G}\left[E_{\omega}^0[N_{x,y}']^s\right]<\infty$.
\end{prop}

\begin{proof}
Let $s<\kappa_1$. It suffices to find $M$ such that for any $a\in\Z$,
\begin{equation}\label{eqn:1850}
E_{\G}\left[E_{\omega}^a[N_{0,M}']^s\right]<\infty.
\end{equation}
This is because for any $x,y$ with $y-x=M$, \eqref{eqn:1850} gives us 
\begin{equation*}
E_{\G}\left[E_{\omega}^{-x}[N_{0,M}']^s\right]<\infty,
\end{equation*}
and then the shift-invariance of $\G$ gives us  $E_{\G}\left[E_{\omega}^0[N_{x,y}']^s\right]<\infty$.
For $y-x>M$, we have $N_{x,y}'\leq N_{x,x+M}'$, which finishes the proof. 

Let $M\geq1$ and and fix $a\in\Z$. 
Suppose $\omega$ is such that jumps to the left and right are bounded by $L$ and $R$, respectively, and the Markov chain is irreducible with almost-sure transience to the right (this is all true with $P_{\G}$-probability 1). Then after the walk traverses $[0,M]$ from right to left $n-1$ times, it will almost surely return to $[M,\infty)$ by transience to the right, and then traversing the interval an $n$th time will require visiting one of the sites $M,\ldots,M+L-1$ and then backtracking past 0. Hence 
\ifSHOWEXTRA
\else
$P_{\omega}^a(N_{0,M}'\geq n\given N_{0,M}'\geq n-1)\leq\max_{0\leq i\leq L-1}P_{\omega}^{M+i}(H_{\leq0}<\infty)$. 
\fi
It is therefore straightforward to check by induction that
\begin{equation*}    
P_{\omega}^a(N_{0,M}'\geq n) \leq \left(\max_{0\leq i\leq L-1}P_{\omega}^{M+i}(H_{\leq0}<\infty)\right)^{n-1}.
\end{equation*}
Summing over all $n\geq1$ we get
\begin{equation}\label{eqn:243}
E_{\omega}^a[N_{0,M}']\leq \frac{1}{\min_{0\leq i\leq L-1}P_{\omega}^{M+i}(H_{\leq0}=\infty)}.
\end{equation}
It therefore suffices to show that the right hand side has finite $s$th moment. 

Notice that by our assumptions on $\omega$, for any $x\geq0$ and $y\geq x$, there is a $z$ in $[y+1,y+R]$ such that $P_{\omega}^z(H_{\leq0}=\infty)\geq P_{\omega}^{x}(H_{\leq0}=\infty)$. This is because if a walk from $x$ is to avoid backtracking to 0, it must enter $[y+1,y+R]$ before backtracking to 0, and then continue to avoid backtracking to 0. Thus, in every interval of length $R$ to the right of 0, there is at least one site $z$ (which depends on $\omega$) such that $P_{\omega}^z(H_{\leq0}=\infty)\geq\max_{1\leq i\leq R}P_{\omega}^i(H_{\leq0}=\infty)$. Call such a site an {\em escape site}, and let 
\begin{equation}
    \mathcal{E}=\mathcal{E}(\omega):=\{x\in[1,M]:x\text{ is an escape site}\}.
\end{equation}
Then for $x>0$,

\begin{align*}
    P_{\omega}^x(H_{\leq0}=\infty)
    &\geq P_{\omega}^x(H_{\mathcal{E}}<H_{\leq0})\max_{1\leq i\leq R}P_{\omega}^i(H_{\leq0}=\infty)+P_{\omega}^x(H_{\mathcal{E}}=H_{\leq0}=\infty)
    \\
    &\geq
     P_{\omega}^x(H_{\mathcal{E}}\leq H_{\leq0})\max_{1\leq i\leq R}P_{\omega}^i(H_{\leq0}=\infty).
\end{align*}

Substituting this into \eqref{eqn:243} we get

\begin{align}
\notag
E_{\omega}^a[N_{0,M}']
&\leq 
\frac{1}{\min_{0\leq i\leq L-1}P_{\omega}^{M+i}(H_{\mathcal{E}}\leq H_{\leq0})}
\cdot
\frac{1}{\max_{1\leq i\leq R}P_{\omega}^i(H_{\leq0}=\infty)}
\\
\label{eqn:343}
&=\frac{1}{1-\max_{0\leq i\leq L-1}P_{\omega}^{M+i}( H_{\leq0}<H_{\mathcal{E}})}
\cdot
\frac{1}{\max_{1\leq i\leq R}P_{\omega}^i(H_{\leq0}=\infty)}
\end{align}

To show that this has finite $s$th moment for large enough $M$ we will use H\"{o}lder's inequality. Choose $s'$ such that $s<s'<\kappa_1$, and let $t=\frac{ss'}{s'-s}$. Thus, for any random variables $X$ and $Y$, where $E[X^t]<\infty$ and $E[Y^{s'}]<\infty$, we will have $E[(XY)^s]<\infty$. We will show that the second term of the right hand side of \eqref{eqn:343} has finite $s'$th moment, and that the first term can have arbitrarily high finite moments for $M$ sufficiently large, so that for large enough $M$, the first term has finite $t$th moment. 

By arguing along the lines of Claim \ref{claim:ComparingGraphs}, we see that the distribution of $\max_{1\leq i\leq R} P_{\omega}^i(H_0=\infty)$ under $P_{\G_+}$ is the distribution of $\max_{1\leq i\leq R} P_{\omega}^i(H_{\leq 0}=\infty)$ under $P_{\G}$.
Recall from Lemma \ref{lem:InfiniteTimeReversal} that under $P_{\G_+}$, $P_{\omega}^0(\tilde{H}_0=\infty)\sim \text{Beta}(\kappa_1,d^{\totheleft})$. Now for $P_{\G_+}$--a.e. environment $\omega$, $P_{\omega}^0(\tilde{H}_0=\infty)\leq \max_{1\leq i\leq R} P_{\omega}^i(H_0=\infty)$ by the Markov property, because $X_1\in[1,R]$ $\Pa_{\G_+}^0$--a.s. We may conclude that
\begin{equation}\label{eqn:1982}
E_{\G}\left[\left(\frac{1}{\max_{1\leq i\leq R} P_{\omega}^i(H_{\leq0}=\infty)}\right)^{s'}\right]\leq E_{\G_+}\left[\left(\frac{1}{P_{\omega}^0(\tilde{H}_{0}=\infty)}\right)^{s'}\right]<\infty.
\end{equation}

We now show that the first term of \eqref{eqn:343} has finite $t$th moment, provided $M$ is chosen large enough. Let $A\in[0,L-1]$ be the maximizer in the denominator this term. That is, $M+A$ is the site within $[M,M+L-1]$ from which there is the highest   probability of hitting $0$ before hitting an escape site between $0$ and $M$. For $i\in[0,L-1]$, let $\omega_i$ be the environment $\omega$, modified at sites other than $M+i$ in $[M,\infty)$ so that the walk jumps from these sites to $M+i$ with probability 1 under $\omega_i$. That is, for $y\geq M$, $\omega_i(y,M+i)=1$.

Now under $\omega$, a walk from any site to the right of $[M,M+L-1]$ must enter $[M,M+L-1]$ before hitting 0. By the strong Markov property, the site in $[M,\infty)$ with the best probability (under $\omega$) of hitting 0 strictly before $\mathcal{E}$ is therefore $A$. Forcing the walk to jump from other sites in $[M,\infty)$ to site $A$ can only increase the probability that the walk hits 0 before $\mathcal{E}$, by Lemma \ref{lem:HarmonicLemma}. Therefore,

\begin{equation*}
\max_{0\leq i\leq L-1}P_{\omega}^{M+i}(H_{\leq0}<H_{\mathcal{E}})
=
P_{\omega}^{M+A}(H_{\leq0}<H_{\mathcal{E}})
\leq
P_{\omega_A}^{M+A}(H_{\leq0}<H_{\mathcal{E}}).
\end{equation*}
From this we get
 \begin{equation}\label{eqn:2003}
 \frac{1}{1-\max_{0\leq i\leq L-1}P_{\omega}^{M+i}(H_{\leq0}< H_{\mathcal{E}})}
\leq
\frac{1}{P_{\omega_A}^{M+A}(H_{\mathcal{E}}\leq H_{\leq0})}
\leq
\sum_{i=0}^{L-1}\frac{1}{P_{\omega_i}^{M+i}(H_{\mathcal{E}}\leq H_{\leq0})},
\end{equation}
and it suffices to show that each term in the sum has finite $t$th moment for large enough $M$. 

Say that a set $W\subseteq[1,M]$ is an {\em escape-type set} if $W$ contains at least one element of every interval of length $R$ contained in $[1,M]$. Then $\mathcal{E}$ is an escape-type set. Now for each escape-type set $W$, consider an environment $\omega_{i,W}$ such that:

\begin{enumerate}
    \item All sites $y\in W(\omega)$ are sinks: for all $y\in W,~z\in\Z$, $\omega_{i, W}(y,z)=\one_{\{z=y\}}$;
    \item For $y=1,\ldots L-1$, $\omega_{i, W}(y,0)=\sum_{z\leq0}\omega_i(y,z)$ and for $z<0$, $\omega_{i, W}(y,z)=0$;
    \item For all $z$, $\omega_{i, W}(0,z)=\one_{\{z=M+i\}}$;
    \item All other transition probabilities are the same in $\omega_{i, W}$ as in $\omega_i$.
\end{enumerate}

By construction, $P_{\omega_i}^{M+i}(H_{\mathcal{E}}\leq H_{\leq0})=P_{\omega_{i,\mathcal{E}}}^0(\tilde{H}_0=\infty)$. Hence
\begin{equation}\label{eqn:2018}
\frac{1}{P_{\omega_i}^{M+i}(H_{\mathcal{E}}\leq H_{\leq0})}=E_{\omega_{i,\mathcal{E}}}^0[N_0].
\end{equation}
We wish to show, with Tournier's lemma, that this quantity has finite $t$th moment for sufficiently large $M$. Since $\mathcal{E}$ is random, $\omega_{i,\mathcal{E}}$ is not a Dirichlet environment, because the set $\mathcal{E}$ of sink sites is random. Nevertheless, for a fixed $M$, there are finitely many possible escape sets. 
Hence it suffices to show that for large enough $M$, $E_{\omega_{i,W}}^0[N_0]$ has finite $t$th moment for every escape-type set $W$. Note that sites outside of the set $[0,M+i+R]$ are unreachable from sites inside the set under $\omega_{i,W}$. 
By the amalgamation property, the restriction of $\omega_{i,W}$ to $[0,M+i+R]^2$ is distributed as a Dirichlet environment on a graph $\G_{M,i,W}$ with vertex set $[0,M+i+R]$ that looks like $\G$ on these vertices except that:
\begin{enumerate}
    \item Directed edges from sites $w\in W$ are removed and replaced with one self-loop at each such site;
    \item For each $y=1,\ldots L-1$, all directed edges from $y$ to sites less than or equal to 0 are replaced with one directed edge to $0$ with the sum of their weights (in our illustration we use multiple edges for visual clarity);
    \item All directed edges from 0 are replaced with one directed edge to $M+i$;
    \item All directed edges from each $y\in [M+L,M+i+R]$ are replaced with one edge to $M+i$.
\end{enumerate}
When there is only one edge from a vertex, its weight does not matter---weight 1 is as good as any. 
Figure \ref{fig:graphGMiWexample} illustrates an example of a graph $\G_{M,i,W}$. 
\begin{figure}
    \centering
    \includegraphics[width=6.5in]{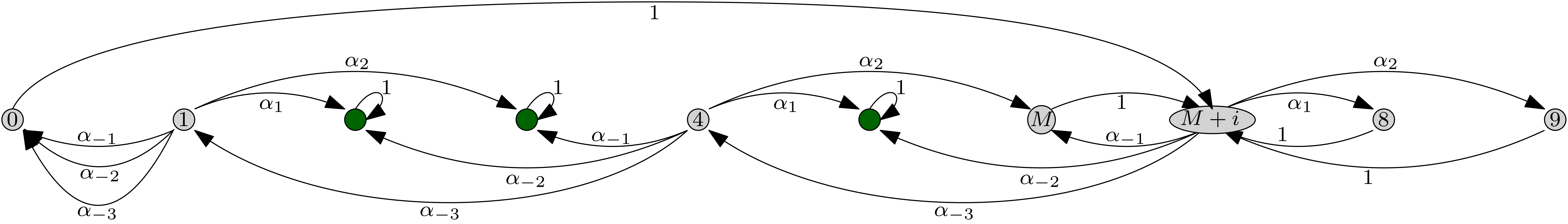}
    \caption{An example of the graph $\G_{M,i,W}$ with $L=3$, $R=2$, $M=6$, $i=1$, and $W=\{2,3,5\}$}
    \label{fig:graphGMiWexample}
\end{figure}
Any strongly connected set $S$ of vertices containing 0 must contain $M+i$ and a path from $M+i$ to $0$, but cannot contain any vertices in $W$. Dividing $[0,M]$ into consecutive intervals of length $m_0$ (where $m_0\geq\max(L,R)$ is large enough that every interval of length $m_0$ is strongly connected in $\G$), we see that every such interval contains a vertex in $S$ (because a path from $M+i$ to $0$ in $\G_{i,W}$ cannot jump over $L$ vertices), and every such interval contains a vertex in $W$ (because such a vertex exists in every interval of length $R$). Hence, by the definition of $m_0$, every such interval must contain an edge from $S$ to a vertex not in $S$. If $M\geq qm_0$, so that $[0,M]$ contains at least $q$ disjoint intervals of length $m_0$, and if $\varepsilon$ is the smallest weight any edge in $\G$ has, then $\beta_S\geq q\varepsilon$. Taking $M$ sufficiently large raises this lower bound above $t$. Fix this large $M$. Then Tournier's lemma ensures that $E_{\G_{M,i,W}}\left[E_{\omega}^0[N_0]^t\right]<\infty$ for all escape-type sets $W$, which implies that $E_{\G}\left[E_{\omega_{i,W}}^0[N_0]^t\right]<\infty$. Since $M$ is fixed, there are finitely many such $W$, so
\begin{equation*}
E_{\G}\left[E_{\omega_{i,\mathcal{E}}}^0[N_0]^t\right]\leq E_{\G}\left[\sum_{W}E_{\omega_{i,W}}^0[N_0]^t\right]<\infty,     
\end{equation*}
where the sum is taken over all escape-type sets $W$. 
Now, by \eqref{eqn:2018}, each term of the sum on the right hand side of \eqref{eqn:2003} has finite $t$th moment, giving finite $t$th moment to the left hand side. 
This is what we needed to complete the proof, as we may now apply H\"older's inequality to \eqref{eqn:343} to see that $E_{\G}\left[E_{\omega}^a[N_{0,M}']^s\right]<\infty$, which is precisely \eqref{eqn:1850}.
\end{proof}

We have now proven each part of Theorem \ref{thm:TFAE2}, which states that if $\kappa_1>0$ and $s>0$, then the following are equivalent:

\begin{enumerate}[(a)]
    \item $\kappa_1>s$.
    \item There is an $M\geq0$ such that for all $x,y\in\Z$ with $y-x\geq M$, $E_{\G}\left[E_{\omega}^0[N_{x,y}']^s\right]<\infty$.
    \item There exist $x<y\in\Z$ such that $E_{\G}\left[E_{\omega}^0[N_{x,y}]^s\right]<\infty$.
\end{enumerate}

\begin{proof}[Proof of Theorem \ref{thm:TFAE2}]\hypertarget{proof:TFAE2}
$(a) \Rightarrow (b)$ This is Proposition \ref{prop:FewOscillations}.

$(b) \Rightarrow (c)$ This follows from the fact that $N_{-M,0}'\geq N_{-M,0}$.

$(c) \Rightarrow (a)$ This is the contrapositive of Proposition \ref{prop:InfiniteTraps} (2).
\end{proof}

\subsection{Using \texorpdfstring{$\kappa_0$}{kappa0} and \texorpdfstring{$\kappa_1$}{kappa1} to characterize ballisticity}\label{subsec:ballisticity}

From Theorems \ref{thm:gamma} and \ref{thm:TFAE2}, we can conclude that if $s\geq \kappa_0$, then $E_{\G}[E_{\omega}^0[N_0]^s]=\infty$ due to finite trapping effects, and that if $s\geq \kappa_1$, then $E_{\G}[E_{\omega}^0[N_0]^s]=\infty$ due to large-scale backtracking effects. On the other hand, if $s<\min(\kappa_0,\kappa_1)$, then neither effect, on its own, is enough to cause $E_{\G}[E_{\omega}^0[N_0]^s]$ to be infinite. However, we must consider the possibility that the two effects could ``conspire together", since the quenched probability of backtracking and hitting 0 is not completely independent of the quenched probability of hitting 0 a large number of times before exiting a small region. Nevertheless, we are able to show that there is enough independence that this is not an issue; $E_{\omega}^0[N_0]$ indeed has finite moments up to the minimum of $\kappa_0$ and $\kappa_1$. In this subsection, we prove Theorem \ref{thm:MainMomentTheorem}, which states that $E_{\G}[E_{\omega}^0[N_0]^s]<\infty$ if and only if $s<\min(\kappa_0,\kappa_1)$. 
Before proving it, we begin with the following lemma. 

\begin{lem}\label{lem:2030}
For any $z>0$, 
\begin{equation*} 
    E_{\G}\left[\frac{1}{P_{\omega}^0(H_{\geq z}<\tilde{H}_0)^s}\right]<\infty\quad\Leftrightarrow\quad s<\min(\kappa_0,d^{\totheright}).
\end{equation*}
\end{lem}
Equivalently, $E_{\G}\left[E_{\omega}^0\left[N_0^{(-\infty,z)}\right]^s\right]=\infty$ if and only if $s<\min(\kappa_0,d^{\totheright})$.
\begin{proof}
First, assume $s<\min(\kappa_0,d^{\totheright})$. 
We will consider an integer $-M<-R$. For an environment $\omega$, let $j(\omega)$ be the vertex $j$ in $[1-M,R-M]$ that maximizes $P_{\omega}^j(H_0<H_{\geq z})$.
Consider an environment $\omega'$ such that:
\begin{enumerate}
    \item for $i\in[z-R,z-1]$, $\omega'(i,z)=\sum_{j\geq z}\omega(i,j)$ and for $j>z$, $\omega'(i,j)=0$;
    \item $z$ is a sink---for all $x$, $\omega'(z,x)=\one_{\{x=z\}}$;
    \item all other transition probabilities are the same in $\omega'$ as in $\omega$.
\end{enumerate}

By construction, $P_{\omega'}^x(\tilde{H}_0<\infty)=P_{\omega'}^x(\tilde{H}_0<H_z)=P_{\omega}^x(\tilde{H}_0<H_{\geq z})$ for all $x$. We now modify the environment further. For each $1\leq \ell\leq R$, let $\omega_{\ell}$ be the environment such that
\begin{enumerate}
    \item for $i\in [-M-L+1,-M]$, $\omega'(i,j)=\one_{\{j=\ell\}}$ for all $j$;
    \item All other transition probabilities are the same in $\omega_{\ell}$ as in $\omega'$. 
\end{enumerate}
In particular, we are interested in $\omega_{j(\omega)}$. Because $j(\omega)$ maximizes $P_{\omega}^j(H_0<H_{\geq z})$ on $[-M-L+1,-M]$,  modifying sites in this set to send the walk directly to $j(\omega)$ can only increase the probability, from any starting point, that $H_0<H_z$ (see Lemma \ref{lem:HarmonicLemma} for details). Thus, 
\begin{equation*}
    P_{\omega_{j(\omega)}}^0(\tilde{H}_0<H_{\geq z})\leq P_{\omega'}^0(\tilde{H}_0<H_{\geq z})=P_{\omega}^0(\tilde{H}_0<H_{\geq z}).
\end{equation*} 
Moreover, for $P_{\G}$--a.e. $\omega$, we have $E_{\omega_{j(\omega)}}^0[N_0]=\frac{1}{P_{\omega_{j(\omega)}}^0(\tilde{H}_0=\infty)}$. It follows that
\begin{equation*}
    \frac{1}{P_{\omega}^0(H_{\geq z}<\tilde{H}_0)}\leq E_{\omega_{j(\omega)}}^0[N_0].
\end{equation*}

It suffices, therefore, to show that $E_{\G}[E_{\omega_{j(\omega})}[N_0]^s]<\infty$. We would like to use Tournier's lemma. Although $\omega_{j(\omega)}$ is not distributed according to a Dirichlet distribution (since $j(\omega)$ is random), the amalgamation property implies that each $\omega_{\ell}$ is distributed as a Dirichlet environment. In particular, the restriction of $\omega_{\ell}$ to $[-M-L+1,z]^2$ is distributed as a Dirichlet environment on a graph $\G_{\ell}$ with vertices $[-M-L+1,z]$ and the following properties:
\begin{enumerate}
    \item each vertex in $[-M-L+1,-M]$ has one edge to $\ell$ with arbitrary weight, say 1;
    \item vertices in $[1-M,z-R-1]$ have the same edges with the same weights as in $\G$;
    \item vertices in $[z-R,z-1]$ have the same edges with the same weights as in $\G$, except that edges that would terminate to the right of $z$ terminate at $z$;\footnote{As usual, if this would result in multiple edges, they are collapsed into one edge with the sum of their weights, but we leave multiple edges in our illustration for clarity.}
    \item $z$ has one self-loop with arbitrary weight, say 1.
\end{enumerate}
With probability 1, no vertices to the left of $-M-L+1$ or to the right of $z$ are reachable from $[-M,z]$, so what happens at these vertices does not really matter. Figure \ref{fig:graphGLexample} shows an example of the graph $\G_{\ell}$.

\begin{figure}
    \centering
    \includegraphics[width=6.5in]{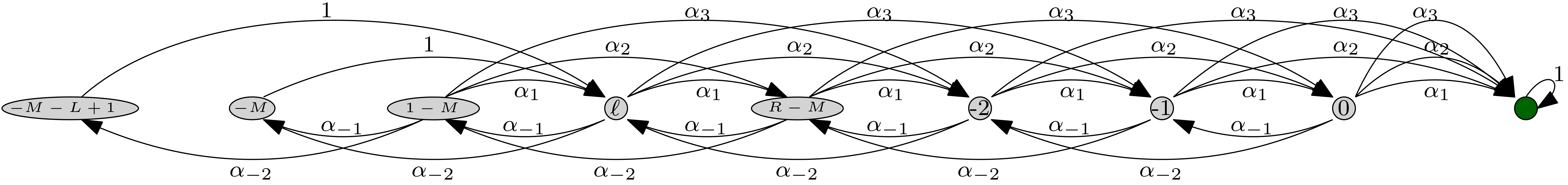}
    \caption{An example of the graph $\G_{\ell}$, where $L=2$, $R=3$, $z=1$, $M=6$, and $\ell=2-M$}
    \label{fig:graphGLexample}
\end{figure}

Moreover, since $1-M\leq j(\omega)\leq R-M$, we have $E_{\omega_{j(\omega)}}^0[N_0]\leq\sum_{\ell=1-M}^{R-M}E_{\omega_{\ell}}^0[N_0]$. Thus, to show that $E_{\G}\left[\frac{1}{P_{\omega}^0(H_{\geq z}<\tilde{H}_0)^s}\right]<\infty$, it suffices to show that $E_{\G}[E_{\omega_{\ell}}^0[N_0]^s]=E_{\G_{\ell}}[E_{\omega}^0[N_0]^s]$ is finite for each $1-M\leq \ell\leq R-M$.

\setcounter{ExampleSave}{\value{exmp}}
\setcounter{exmp}{\value{thm}}
\setcounter{claim}{0}

We make the following claim.
\begin{claim}
If $M$ is chosen large enough, $\G_{\ell}$ will have the property that every finite, strongly connected set $S$ of vertices with $0\in S$ has $\beta_S\geq\min(d^{\totheright},\kappa_0)$.
\end{claim}

\setcounter{exmp}{\value{ExampleSave}}

The proof of this claim is very similar to the proof of Proposition \ref{prop:GammaAttained}, so we only sketch it here, making special note of important similarities to and differences from the earlier proof. Like $\G$, the graph $\G_{\ell}$ has the property that for any finite, strongly connected set set $S$ of vertices, if $x<y$ are consecutive ``non-insulated" vertices to the right of $-M$ that differ by more than $m_0$, then $(x,y)\subset S$. The graph $\G_{\ell}$ also has the property that every strongly connected set $S'$ that contains $m_0$ or more consecutive vertices has total weight at least $d^{\totheright}$ exiting $S'$ to the right. However, because of how the graph is modified, there need not be any weight exiting $S'$ to the left. We may therefore choose $M$ large enough that every finite, strongly connected set $S$ containing 0 either (a) does not reach $-M$, in which case $S$ looks like a subset of $\G$ and $\beta_{S}\geq\kappa_0$, or (b) contains at least $m_0$ consecutive vertices in a row, in which case $\beta_{S}\geq d^{\totheright}$, or (c) contains enough non-insulated vertices that $\beta_{S}\geq \min(d^{\totheright},\kappa_0)$. This completes the proof of the claim.

Since we have assumed $s<\min(\kappa_0,d^{\totheright})$, our claim gives us $s<\beta_S$ for all strongly connected sets $S$ containing $0$, and therefore Tournier's Lemma tells us that $E_{\G_{\ell}}[E_{\omega}^0[N_0]^s]<\infty$. 

This finishes one direction of the lemma, and the only direction that is needed for the rest of the paper. For the other direction, assume $s\geq\min(\kappa_0,d^{\totheright})$. 
If $s\geq\kappa_0$, then Theorem \ref{thm:gamma} implies $E_{\G}\left[E_{\omega}^0\left[N_0^{[-M,0]}\right]^s\right]=\infty$, which implies $E_{\G}\left[E_{\omega}^0\left[N_0^{(-\infty,z)}\right]^s\right]=\infty$. If $s\geq d^{\totheright}$, then one can check that $E_{\G}\left[\left(\sum_{i<z,j\geq z}\omega(i,j)\right)^{-s}\right]=\infty$, and then arguments along the lines of the last part of the proof of Proposition \ref{prop:InfiniteTraps} give the desired result. 
\end{proof}

We are now ready to prove Theorem \ref{thm:MainMomentTheorem}, which gives us the final piece for our characterization of ballisticity. 

\begin{proof}[Proof of Theorem \ref{thm:MainMomentTheorem}]\hypertarget{proof:MainMomentTheorem}
The forward direction is implied by Theorem \ref{thm:gamma} and Theorem \ref{thm:TFAE2}. For the reverse direction, let $s<s'<\min(\kappa_1,\kappa_0)$. We want to show that $E_{\G}\left[E_{\omega}^0[N_0]^s\right]<\infty$. Since, for all $\omega$, $E_{\omega}^0[N_0]=\frac{1}{P_{\omega}^0(\tilde{H}_0=\infty)}$, we must examine the quantity $P_{\omega}^0(\tilde{H}_0=\infty)$.

Let $M=qm_0$ be a positive multiple of $m_0$ (we will later take $q$ to be large enough to satisfy a given condition, but we will not take $q$ to infinity). For each $1\leq j\leq q$, let $B_j=((j-1)m_0,jm_0]$, and define
\begin{equation*}
f_{\omega}(B_j):=\min_{x\in B_j, 1\leq i \leq R}P_{\omega}^x(X_{H_{B_j^c}}=jm_0+i)
\end{equation*}
Since $f_{\omega}(B_j)$ only depends on $\omega^{B_j}$, the $f_{\omega}(B_j)$ are i.i.d. 

We will show that for some $a,C>0$, $P_{\G}(f_{\omega}(B_1)<\tepsilon)\leq C\tepsilon^a$. For a given vertex $x$, let $\underline{\omega^x}=\min_{\alpha_{y-x}>0}\omega(x,y)$.
Since intervals of length $m_0$ are strongly connected, 
there is a path in $\G$ from every $x\in B_1$ to every $m_0+i$, $1\leq i\leq R$, that uses only vertices in $B_1$, and each at most once (since the existence of a path with loops implies the existence of a path without loops).
Therefore, 
\begin{equation*}
    f_{\omega}(B_1)\geq\left(\prod_{x\in B_1}\underline{\omega^x}\right).
\end{equation*}
Now for each $x$,
\begin{align*}
P_{\G}\left(\underline{\omega^x}<\tepsilon\right)
&\leq(R+L)\max_{\alpha_{y-x}>0}P_{\G}(\omega(x,y)<\tepsilon)
\\
&\leq C\tepsilon^{a}
\end{align*}
for some constant $C>0$, where $a=\min_{\alpha_i\neq0}\alpha_i$. Now by \cite[Lemma 9]{Tournier2009}, there is some $r$ such that $P_{\G}(\prod_{x\in B_1}\underline{\omega^x}\leq \tepsilon)\leq C'\tepsilon^a(-\ln \tepsilon)^r$ for all $\tepsilon$ sufficiently small. Let $a'<a$; then for large enough $C''$, we have $P_{\G}(\prod_{x\in B_1}\underline{\omega^x}<\tepsilon)\leq C''\tepsilon^{a'}$, and therefore $P_{\G}(f_{\omega}(B_1)<\tepsilon)\leq C''\tepsilon^{a'}$

Let $M=m_0q$ be large enough that $a'q>t:=\frac{ss'}{s'-s}$. Since all the $f_{\omega}(B_j)$ are i.i.d., we have
\begin{equation}\label{eqn:jbest}
P_{\G}\left(\max_{0\leq j\leq q-1}f_{\omega}(B_j)<\tepsilon\right)\leq (C''\tepsilon^{a'})^{q}.
\end{equation}
Now let $j_{\best}$ be the maximizer of $f_{\omega}(B_j)$ over $0\leq j\leq q-1$. By \eqref{eqn:jbest} and the choice of $M$, we have
\begin{equation}\label{eqn:jbestmoments}
    E_{\G}\left[\left(\frac{1}{f_{\omega}(B_{j_{\best}})}\right)^t\right]<\infty.
\end{equation}
As in the proof of Proposition \ref{prop:FewOscillations}, we note that for any random variables $X$ and $Y$, where $E[X^t]<\infty$ and $E[Y^{s'}]<\infty$, we will have $E[(XY)^s]<\infty$ by H\"older's inequality. 

For $P_{\G}$--a.e. environment $\omega$ we have
\begin{align*}
    P_{\omega}^0(\tilde{H}_0=\infty)
    &\geq 
    P_{\omega}^0(H_{B_{j_{\best}}}<\tilde{H}_0)
    f_{\omega}(B_{j_{\best}})
    \max_{1\leq i\leq R}P_{\omega}^{j_{\best}m_0+i}(H_{\leq j_{\best}m_0}=\infty).
\end{align*}
This is because one way for the walk to never return to 0 is for it to hit $B_{j_{\best}}$ before returning to 0, then once it is in $B_{j_{\best}}$, to make its way to the vertex $j_{\best}m_0+i$ just to the right of $B_{j_{\best}}$ that maximizes the probability of never backtracking to $j_{\best}m_0$, and then to avoid backtracking to $j_{\best}m_0$. Thus

\begin{align}
\notag
    \frac{1}{P_{\omega}(\tilde{H}_0=\infty)}
    &\leq\frac{1}{P_{\omega}^0(H_{B_{j_{\best}}}<\tilde{H}_0)f_{\omega}(B_{j_{\best}})\max_{1\leq i\leq R}P_{\omega}^{j_{\best}m_0+i}(H_{\leq j_{\best}m_0}=\infty)}
    \\
    \notag
    &=\sum_{j=1}^{q}\frac{\one_{\{j=j_{\best}\}}}{P_{\omega}^0(H_{B_j}<\tilde{H}_0)f_{\omega}(B_{j_{\best}})\max_{1\leq i\leq R}P_{\omega}^{jm_0+i}(H_{\leq jm_0}=\infty)}
    \\
    \label{73}
    &\leq \frac{1}{f_{\omega}(B_{j_{\best}})}\sum_{j=1}^{q}\frac{1}{P_{\omega}^0(H_{B_j}<\tilde{H}_0)}\cdot\frac{1}{\max_{1\leq i\leq R}P_{\omega}^{jm_0+i}(H_{\leq jm_0}=\infty)}
\end{align}
Now for each fixed $j$, $P_{\omega}^0(H_{B_j}<\tilde{H}_0)$ and $\max_{1\leq i\leq R}P_{\omega}^{jm_0+i}(H_{\leq jm_0}=\infty)$ are independent. Under $P_{\G}$, the reciprocal of the former has finite $s'$th moment by Lemma \ref{lem:2030}, since $s'\leq\min(\kappa_0,\kappa_1)\leq\min(\kappa_0,d^{\totheright})$.

We show that the reciprocal of $\max_{1\leq i\leq R}P_{\omega}^{jm_0+i}(H_{\leq jm_0}=\infty)$ has finite $s'$th moment. Under $P_{\G}$, its distribution is the same as that of $\max_{1\leq i\leq R}P_{\omega}^{i}(H_{\leq 0}=\infty)$, which (arguing as in Claim \ref{claim:ComparingGraphs}) is also the distribution of $\max_{1\leq i\leq R}P_{\omega}^{i}(H_0=\infty)$ under $P_{\G_+}$. Now for $P_{\G_+}$--a.e. environment $\omega$, $P_{\omega}^0(\tilde{H}_0=\infty)\leq \max_{1\leq i\leq R}P_{\omega}^{i}(H_0=\infty)$. Hence
\begin{equation}\label{eqn:2164}
    E_{\G}\left[\frac{1}{\max_{1\leq i\leq R}P_{\omega}^{jm_0+i}(H_{\leq jm_0}=\infty)^{s'}}\right]
    \leq
    E_{\G_+}\left[\frac{1}{P_{\omega}^0(\tilde{H}_0=\infty)^{s'}}\right].
\end{equation}
Lemma \ref{lem:InfiniteTimeReversal} tells us that under $P_{\G_+}$, $P_{\omega}^0(\tilde{H}_0=\infty)\sim \text{Beta}(\kappa_1,d^{\totheleft})$. Since $s'\leq\kappa_1$, the right side of \eqref{eqn:2164} is finite by \eqref{eqn:moments}. 

Returning now to \eqref{73}, we have seen that the two fractions inside the sum are independent and each have finite $s'$th moment. Thus, the entire sum has finite $s'$th moment. By this and by \eqref{eqn:jbestmoments}, we may apply H\"older's inequality to get

\ifSHOWEXTRA
\else
\begin{equation*}
     E_{\G}[E_{\omega}^0[N_0]^s]
    =
    E_{\G}\left[
    \frac{1}{P_{\omega}(\tilde{H}_0=\infty)^s}\right]
    <\infty.
\end{equation*}
\fi
\end{proof}


We are now able to give a complete characterization of ballisticity. Assume $\kappa_1>0$. Theorem \ref{thm:MainBallisticTheorem} states that the walk is ballistic if and only if $\min(\kappa_1,\kappa_0)>1$. 

\begin{proof}[Proof of Theorem \ref{thm:MainBallisticTheorem}]\hypertarget{proof:MainBallisticTheorem}
Let $s=1$ in Theorem \ref{thm:MainMomentTheorem}, and then apply Lemma \ref{lem:ballistic}.
\end{proof}

\appendix

\section{Auxiliary Proofs}\label{append:Leftoverproofs}

We outline a proof of Proposition \ref{prop:LimitingVelocity}.

Recall that we have defined $\tau_0:=0$, and for $k\geq1$, $\tau_k:=\min\{n>\tau_{k-1}:X_n>X_j\text{ for all }j<n,X_n\leq X_j\text{ for all }j>n\}$.

We are to prove that for $P$ satisfying (\hyperlink{cond:C1}{C1}), (\hyperlink{cond:C2}{C2}), (\hyperlink{cond:C3}{C3}), and (\hyperlink{cond:C4}{C4}), the following hold:
\begin{enumerate}
    \item There is a $\Pa^0$--almost sure limiting velocity 
\begin{equation}\label{eqn:2034}
    v:=\lim_{n\to\infty}\frac{X_n}{n}=\frac{\Ea^0[X_{\tau_2}-X_{\tau_1}]}{\Ea[\tau_{2}-\tau_{1}]},
\end{equation}
where the fraction is understood to be 0 if the denominator is infinite, and the numerator is always finite.
    \item $\lim_{x\to\infty}\frac{H_{\geq x}}{x} = \frac{1}{v}$, $\Pa^0$--a.s., where $\frac{1}{v}$ is understood to be $\infty$ if $v=0$.
\end{enumerate}

It is standard (see, for example, \cite{Kesten1977}, \cite{Sznitman&Zerner1999}) to prove a LLN like \eqref{eqn:2034} by the following steps:
\begin{enumerate}[(a)]
    \item Show that $\frac{X_{\tau_k}}{k}$ approaches $\Ea[X_{\tau_2}-X_{\tau_1}]$
    \item Show that $\frac{\tau_k}{k}$ approaches $\Ea[\tau_2-\tau_1]$
    \item Show that $\Ea[X_{\tau_2}-X_{\tau_1}]<\infty$
    \item Conclude that the limit \eqref{eqn:2034} holds for the subsequence $\left(\frac{X_{\tau_k}}{\tau_k}\right)_k$
    \item Use straightforward bounds that come from the definitions of the $\tau_k$ to get the limit for the entire sequence $\left(\frac{X_n}{n}\right)_n$. 
\end{enumerate}
Part (2) then follows from a comparison of $\frac{x}{H_{\geq x}}$ with a subsequence of $\frac{X_n}{n}$. 

The definition of the regeneration times is precisely set up so that both the sequences $(\tau_k-\tau_{k-1})_{k\geq2}$ and $(X_{\tau_k}-X_{\tau_{k-1}})_{k\geq2}$ are i.i.d. sequences, so proving the limits (a) and (b) is a matter of tracing how the i.i.d. feature follows from the definitions and applying the strong law of large numbers. In fact, arguing as in \cite[Lemma 1]{Kesten1977}, one can show that the triples
\begin{equation*}
    \xi_n:=\left(\tau_n-\tau_{n-1}~,~(X_{\tau_{n-1}+i}-X_{\tau_{n-1}})_{i=1}^{\tau_n-\tau_{n-1}}~,~(\omega^x)_{x=X_{\tau_{n-1}}}^{X_{\tau_n}-1}\right)
\end{equation*}
are i.i.d. under $\Pa^0=P\times P_{\omega}^0$ for $n\geq2$.

The finiteness in (c) can be shown using arguments along the lines of those in \cite[Lemma 3.2.5]{Zeitouni2004}. Because we have assumed almost-sure transience to the right, the measure $\Q$ introduced there is unnecessary. Another difference is that in our model, transience to the right does not imply that every vertex to the right of the origin is hit. So instead of studying the probability, for a given $x$, that a regeneration occurs at site $x$, we focus on the probability that the regeneration occurs on a given interval of length $R$. For $z\geq0$, let $B_z$ be the event that for some $k$, $X_{\tau_k}\in[zR,(z+1)R)$. Then
\begin{align}
    \notag
    \Pa^0(B_z)&=E\left[P_{\omega}^0(B_z)\right]
    \\
    \notag
    &\geq
    E\left[\sum_{i=0}^{R-1}P_{\omega}^0(H_{[zR,(z+1)R)}=zR+i)P_{\omega}^{zR+i}(H_{<zR+i}=\infty)\right]
    \\
    \notag
    &= 
    \sum_{i=0}^{R-1}E\left[P_{\omega}^0(H_{[zR,(z+1)R)}=zR+i)P_{\omega}^{zR+i}(H_{<zR+i}=\infty)\right]
    \\
\notag
    &=\sum_{i=0}^{R-1}\Pa^0(H_{[zR,(z+1)R)}=zR+i)\Pa^{zR+i}(H_{<zR+i}=\infty)
    \\
    \label{eqn:1935}
    &=\Pa^{0}(H_{<0}=\infty),
\end{align}
where the second to last equality comes from the fact that $\omega^{<zR}$ is independent of $\omega^{\geq zR+i}$, and the last comes from translation invariance and the fact that $H_{[zR,(z+1)R)}<\infty$ $\Pa^0$--a.s. The rest of the argument from \cite{Zeitouni2004} goes through to prove (c), and (d) and (e) easily follow.

We now prove a lemma that is used to justify arguments where we replace transition probabilities at a particular site $x$ with an almost-sure jump to another site $y$.

\begin{lem}\label{lem:HarmonicLemma}
Consider an environment $\omega$ on a finite or countable vertex set $V$. Assume the Markov chain is irreducible under $\omega$. Let $A$ and $B$ be disjoint finite subsets of $V$, with $A$ nonempty. Now let $x\in V-A\cup B$ and $y\in V-B$, with $P_{\omega}^y(H_A<H_B)\geq P_{\omega}^x(H_A<H_B)$.
Suppose $\omega'$ is an environment that agrees with $\omega$ at all sites other than $x$ but has $\omega'(x,y)=1$. Then for all $z\in V$, $P_{\omega'}^z(H_A<H_B)\geq P_{\omega}^z(H_A<H_B)$ 
\end{lem}

\begin{proof}
We first assume that $V$ is finite. For $z\in V$, let $f_0(z)=P_{\omega}^z(H_A<H_B)$. Then, for $i\geq1$, define

\begin{equation*}
    f_i(z):=\begin{cases}
    \sum_{w\in V}\omega'(z,w)f_{i-1}(w)~~~&z\notin A\cup B
    \\0 & z\in B
    \\1 & z\in A
    \end{cases}.
\end{equation*}

Now for $z\neq x$, $f_1(z)=f_0(z)$, since $f$ is harmonic with respect to $\omega$. On the other hand, $f_1(x)=f_0(y)\geq f_0(x)$. Thus, $f_1(z)\geq f_0(z)$ for all $z$. A straightforward induction now shows that for all $z\in V$, the sequence $(f_i(z))$ is increasing and bounded. Hence $f(z):=\lim_{i\to\infty}f_i(z)$ exists, and $f(z)\geq f_0(z)$ for all $z$.

We will have completed the proof if we can show that $f(z)=P_{\omega'}^z(H_A<H_B)$. One can easily check that $f$ is $\omega'$-harmonic on $V-A\cup B$, identically 1 on $A$, and identically 0 on $B$. The function $z\to P_{\omega'}^z(H_A<H_B)$ has these same properties, and by the maximum principle there is only one such function.

Now suppose $V$ is infinite, and let $z_0\in V$. Take any finite subset $S\subset V$ containing $A$, $B$, $x$, $y$, and $z_0$, and take one ``sink" vertex $\partial$ not in $S$. For $z,w\in S$, define $\omega_{\ast}(z,w):=P_{\omega}^z(\tilde{H}_S<\infty, X_{\tilde{H}_S}=w)$, and $\omega_{\ast}(z,\partial):=P_{\omega}^z(\tilde{H}_S=\infty)$. Define $(\omega')_{\ast}$ similarly. Then for all $z\in S$,
\begin{equation}\label{eqn:1724}
P_{\omega}^z(H_A<H_B)=P_{\omega_{\ast}}^z(H_A<H_B).
\end{equation} 
To see this, it is straightforward to check that $z\to P_{\omega}^z(H_A<H_B)$ is harmonic with respect to $\omega_{\ast}$ on $S-\{A\cup B\}$, identically 1 on $A$, and identically $0$ on $B\cup\{\partial\}$. The function $z\to P_{\omega_{\ast}}^z(H_A<H_B)$ has these same properties, and there is only one such function. Similarly, \begin{equation}\label{eqn:1728}
P_{\omega'}^z(H_A<H_B)=P_{(\omega')_{\ast}}^z(H_A<H_B).
\end{equation}

We want to show that $(\omega')_{\ast}$ has the same properties relative to $\omega_{\ast}$ that $\omega'$ has relative to $\omega$, so that we may apply the conclusion from the finite case with $(\omega_{\ast})':=(\omega')_{\ast}$. That is, we need to show that 

(a) For $z\neq x$, and $w\in S\cup\{\partial\}$,  we have $(\omega')_{\ast}(z,w)=\omega_{\ast}(z,w)$.

(b) $(\omega')_{\ast}(x,y)=1$.

The statement (b) comes from the fact that $\omega'(x,y)=1$. The statement (a) says that $P_{\omega'}^z(\tilde{H}_S<\infty,X_{\tilde{H}_S}=w)=P_{\omega}^z(\tilde{H}_S<\infty,X_{\tilde{H}_S}=w)$. This is true because these probabilities do not depend on the environment at $x$, and $\omega\equiv\omega'$ everywhere else.

We may now apply the finite case to the measures $\omega_{\ast}$ and $(\omega')_{\ast}$, concluding that 
$P_{(\omega')_{\ast}}^{z_0}(H_A<H_B)\geq P_{\omega_{\ast}}^{z_0}(H_A<H_B)$. By \eqref{eqn:1724} and \eqref{eqn:1728}, we get 
$P_{\omega'}^{z_0}(H_A<H_B)\geq P_{\omega}^{z_0}(H_A<H_B)$. Since $z_0$ was arbitrary, this is enough to conclude the argument.
\end{proof}

\section{Calculating \texorpdfstring{$\kappa_0$}{kappa0}}\label{append:kappa}

In the body of the paper, we gave an algorithm to compute $\kappa_0$, but this algorithm grows in complexity with the smallest positive $\alpha_i$. Here, we prove Proposition \ref{prop:GammaAttainedweakstrong}, which asserts that given an underlying directed graph, $\kappa_0$ is a minimum of finitely many positive integer combinations of the $\alpha_i$. Although we do not have a general algorithm to find the formula, we give several examples where we are able to do so. These examples exhibit various important features that sets $S$ attaining $\beta_S=\kappa_0$ may have. We also prove Proposition \ref{prop:orderedpair}, showing that $\kappa_0$ and $\kappa_1$ are unrestricted by each other. We begin with the following lemma.
\begin{lem}\label{lem:SetTheory}
Let $\allowable\subseteq\N^{k}$ be a set of ordered $k$-tuples of positive integers. Let $\sqsubseteq$ be the natural partial ordering on $\N^k$, $(n_1,\ldots,n_k)\sqsubseteq(n_1',\ldots,n_k')$ if $n_i\leq n_i'$ for all $i$. Then there is a finite subset $\allowable^*\subseteq\allowable$ such that for all ${\bf x}\in\allowable$, there is an ${\bf x}^*\in\allowable^*$ such that ${\bf x}^*\sqsubseteq{\bf x}$.
\end{lem}

\begin{proof}
We prove this by induction on $k$. The base case $k=1$ is trivial; a 1-tuple is simply a positive integer, and we can let $\allowable^*:=\{\min\allowable\}$.
Now suppose the result is true for all subsets of $\N^k$, and let $\allowable\subseteq \N^{k+1}$. Now let $\underline{\allowable}:=\{(n_1,\ldots,n_k):(n_1,\ldots,n_k,n_{k+1}\in\allowable\text{ for some }n_{k+1}\in\N\}$ be the projection of $\allowable$ onto $\N^k$, and for $n\in\N$ define $\underline{\allowable}(n):=\{(n_1,\ldots,n_k):(n_1,\ldots,n_k,n)\in\allowable\}$. Thus, $\underline{\allowable}=\bigcup_{n=1}^{\infty}\underline{\allowable}(n)$.  
Now by the inductive hypothesis, there is a finite set $\underline{\allowable}^*\subseteq \underline{\allowable}$ such that every element of $\underline{\allowable}$ is greater than some element of $\underline{\allowable}^*$. Since $\underline{\allowable}^*$ is finite, $\underline{\allowable}^*\subseteq\bigcup_{n=1}^{N}\underline{\allowable}(n)$ for some $N$. Applying the inductive hypothesis to each $\underline{\allowable}(n)$ gives us sets $\underline{\allowable}^*(n)$ such that for all ${\bf x}\in\underline{\allowable}(n)$, there is a ${\bf x}^*\in\underline{\allowable}^*(n)$ such that ${\bf x}^*\sqsubseteq{\bf x}$.

Now define $\allowable^*:=\bigcup_{n=1}^N\{(n_1,\ldots,n_k,n):(n_1,\ldots,n_k)\in\underline{\allowable}^*(n)\}$. It is easy to see that this is a finite subset of $\allowable$. Now suppose $(n_1,\ldots,n_k,n_{k+1})\in\allowable$. 
Suppose $n_{k+1}< N$. Then $(n_1,\ldots,n_k)\in\underline{\allowable}(n_{k+1})$, so there exists $(n_1^*,\ldots,n_k^*)\in\underline{\allowable}^*(n_{k+1})$ with $n_i^*\leq n_i$ for $1\leq i\leq k$. Since $(n_1^*,\ldots,n_k^*,n_{k+1})\in\allowable^*$, we are done. 

On the other hand, suppose $n_{k+1}\geq N$. Then since $(n_1,\ldots,n_k)\in\underline{\allowable}$, there exists $(n_1^*,\ldots,n_k^*)\in\underline{\allowable}^*$ such that $n_i^*\leq n_i$ for all $i=1,\ldots,k$. Now $(n_1^*,\ldots,n_k^*)\in\underline{\allowable}(n)$ for some $n\leq N\leq n_{k+1}$. Hence there exists $(n_1^{**},\ldots,n_k^{**})\in\underline{\allowable}^*(n)$ such that $n_i^{**}\leq n_i^*$ for all $1\leq i\leq k$. Thus $(n_1^{**},\ldots,n_k^{**},n)\in\allowable^*$ with $n_i^{**}\leq n_i$ for $1\leq i\leq k$, and $n\leq n_{k+1}$. \end{proof}

We are now able to prove Proposition \ref{prop:GammaAttainedweakstrong}.

\begin{proof}[Proof of Proposition \ref{prop:GammaAttainedweakstrong}]
For any finite set $S\subset \Z$, $\beta_S$ is a sum of weights of edges exiting $\beta$. The weight of each edge is $\alpha_i$ for some $i$, and each $\alpha_i$ must be included at least once, as the weight of an edge exiting either the rightmost or leftmost point of $S$. Thus,
\begin{equation*}
    \beta_S=\sum_{i=-L}^Rx_i\alpha_i
\end{equation*}
where $x_i=x_i(S):=\#\{z\in S:z+i\notin 
S\}\geq1$.

Now let $\allowable\subset\N^{R+L}$ be the set of ordered tuples $(y_{-L},\ldots,y_R)$ such that there is some finite set $S$ with $x_i(S)=y_i$ for all $-L\leq i\leq R$. Thus,
\begin{equation*}
    \kappa_0=\inf\left\{\sum_{i=-L}^Ry_i\alpha_i:(y_{-L},\ldots,y_R)\in\allowable\right\}.
\end{equation*}
 Applying Lemma \ref{lem:SetTheory}, we get a finite set $\allowable^*\subseteq\allowable$ such that for any $S$, there is a $(y_{-L},\ldots,y_R)\in\allowable^*$ with $y_i\leq x_i(S)$ for all $-L\leq i\leq R$. Thus,
\begin{equation}\label{kappa0elementaryfunction}
    \kappa_0=\min\left\{\sum_{i=-L}^Ry_i\alpha_i:(y_{-L},\ldots,y_R)\in\allowable^*\right\}.
\end{equation}
This is a minimum of finitely many positive integer combinations of the $\alpha_i$.
\end{proof}

We now give examples where we can find the formula for $\kappa_0$. Recall that $\kappa_0:=\inf\{\beta_S:S\subset \Z\text{ finite, strongly connected}\}$, where $\beta_S$ is the sum of edge weights leaving the set $S$. By shift invariance of the graph $\G$, it suffices to consider sets $S$ whose leftmost point is 0. 

We already showed in Claim \ref{claim:GammaClaim1} that $\kappa_0\leq d_{\totheright}+d^{\totheleft}$. We can also give a general lower bound: $\kappa_0\geq c^{\totheright}+c^{\totheleft}$. This is because any strongly connected set will have weight at least $c^{\totheright}$ exiting from the rightmost point and weight at least $c^{\totheleft}$ exiting from the leftmost point. Therefore, every strongly connected set $S$ has $\beta_S\geq c^{\totheright}+c^{\totheleft}$, and taking the infimum preserves the inequality. So we have the bounds
\begin{equation}\label{eqn:boundsforkappa0}
    c^{\totheright}+c^{\totheleft}\leq \kappa_0\leq d^{\totheright}+d^{\totheleft}.
\end{equation}

\begin{exmp}\label{kappa0example1}
$L=R=1$.

In this case, $d^{\totheright}=c^{\totheright}=\alpha_1$, and $d^{\totheleft}=c^{\totheleft}=\alpha_{-1}$, so \eqref{eqn:boundsforkappa0} immediately implies $\kappa_0=\alpha_1+\alpha_{-1}$. This can also be seen by noting that the only strongly connected sets are intervals, which all have the same exit weight.
\end{exmp}

\begin{exmp}\label{kappa0example7}
$\alpha_0>0$. 
 
In this case, $\{0\}$ is already a strongly connected set, so \eqref{eqn:boundsforkappa0} gives $\kappa_0=\beta_{\{0\}}=c^{\totheright}+c^{\totheleft}$. 
\end{exmp}

\begin{exmp}\label{kappa0example2}
$L=2$, $R=3$, $\alpha_i=0$ for $i=-1,\ldots,2$. 

In this case, we also have $\kappa_0=d^{\totheright}+d^{\totheleft}$. Let $S$ be a strongly connected finite set of vertices with left endpoint 0. Then $S$ contains 3, since there must be a vertex reachable from 0 in one step, and by assumption there are no vertices to the left of 0. Also, $S$ contains 2, since 0 must be reachable in one step from a vertex in $S$. Likewise, 2 must then also be reachable, and since $-1\notin S$, $S$ must contain 4 as well. Now since a vertex must be reachable from 3, $S$ must contain either 1 or 6. Suppose $S$ contains 1. Then $S$ contains $[0,4]$, which has exit weight $d^{\totheright}+d^{\totheleft}$, and by Claim \ref{claim:GammaClaim3}, $\beta_S\geq \beta_{[0,4]}=d^{\totheright}+d^{\totheleft}$. On the other hand, suppose $S$ does not contain 1. Then it contains 6. If $S$ also contains 5, then $S$ contains the interval $[2,6]$, which is shift-equivalent to $[0,4]$. If $S$ contains neither 1 nor 5, then it contains exactly the set $\{0,2,3,4,6\}$ and possibly vertices to the right and/or left of this set, so by Claim \ref{claim:GammaClaim3}, $\beta_S\geq\beta_{\{0,2,3,4,6\}}$. One can easily check that in this case, $\beta_{\{0,2,3,4,6\}}=d^{\totheright}+d^{\totheleft}$. 

We have calculated $\kappa_0$ without even showing that either $[0,4]$ or $\{0,2,3,4,6\}$ is strongly connected, but in fact they both are. Consider the path $0\to3\to1\to4\to2\to0$ in $[0,4]$ and the path $0\to3\to6\to4\to2\to0$ in $\{0,2,3,4,6\}$. Thus, even when $\kappa_0=d^{\totheright}+d^{\totheleft}$, a minimizing set $S$ for $\beta_S$ need not be an interval (although a large enough interval will always be a minimizing set). 
\end{exmp}

\begin{exmp}\label{kappa0example3}
$L=1$, $R\geq2$, $\alpha_0=0$ $\alpha_1>0$, $\alpha_i=0$ for $i=2,\ldots, R-1$.

In this case, we show that $\kappa_0=2\alpha_R+\alpha_1+\alpha_{-1}$; thus if $R>2$, then $\kappa_0<d^{\totheright}+d^{\totheleft}$. Let $S$ be a strongly connected set with left endpoint 0. Since 0 must be reachable from another point in $S$, we have $1\in S$. Now by Claim \ref{claim:GammaClaim3}, this implies $\beta_S\geq\beta_{\{0,1\}}=2\alpha_R+\alpha_1+\alpha_{-1}$. The set $\{0,1\}$ is strongly connected, and hence $\kappa_0=\beta_{\{0,1\}}=2\alpha_R+\alpha_1+\alpha_{-1}$. 
\end{exmp}

This example lets us show that $\kappa_0$ and $\kappa_1$ are independent in the sense that no information about either may be inferred from the other. 

\begin{prop}\label{prop:orderedpair}
The ordered pair $(\kappa_0,\kappa_1)$ may take on any value in the first quadrant of $\R^2$.
\end{prop}

\begin{proof}
Let $a,b>0$. We will show that $L$, $R$, and the $\alpha_i$ may be chosen such that $\kappa_0=a$ and $\kappa_0=b$. Let $L=1$, and let $R\geq2$ be large enough that $\frac{a}{2}>\frac{b}{R}$. Then let $\alpha_1=\alpha_{-1}=\frac{a}{2}-\frac{b}{R}$, $\alpha_R=\frac{b}{R}$, and all other $\alpha_i=0$. Then $\kappa_1=-\alpha_{-1}+\alpha_1+R\alpha_R=b$, and by Example \ref{kappa0example3}, $\kappa_0=a$. 
\end{proof}

\begin{exmp}\label{kappa0example4}
$L=R=2$, $\alpha_{-1},\alpha_1>0$, $\alpha_0=0$.

There are two possibilities. Let $S$ be a strongly connected set with leftmost point 0. If $1\in S$, then by Claim \ref{claim:GammaClaim3}, $\beta_S\geq\beta_{\{0,1\}}=d^{\totheright}+d^{\totheleft}$. On the other hand, if $1\notin S$, then $2\in S$, and by Claim \ref{claim:GammaClaim3}, $\beta_S\geq\beta_{\{0,2\}}=\alpha_{-2}+2\alpha_{-1}+2\alpha_1+\alpha_2$. Since both $\{0,1\}$ and $\{0,2\}$ are strongly connected, $\kappa_0$ may be either $d^{\totheright}+d^{\totheleft}$ or $\alpha_{-2}+2\alpha_{-1}+2\alpha_1+\alpha_2$, depending on whether $\alpha_{-1}+\alpha_1$ or $\alpha_{-2}+\alpha_2$ is smaller. That is, 
\begin{align*}
        \kappa_0&=\min(2\alpha_{-2}+\alpha_{-1}+\alpha_1+\alpha_2,\alpha_{-2}+2\alpha_{-1}+2\alpha_1+\alpha_2)
        \\
        &=\alpha_{-2}+\alpha_{-1}+\alpha_1+\alpha_2+\min(\alpha_{-1}+\alpha_1,\alpha_{-2}+\alpha_2)
\\
&=\min(\beta_{\{0,1\}},\beta_{\{0,2\}})
. 
\end{align*}
\end{exmp}

\begin{exmp}\label{kappa0example5}
$L=6$, $R=3$, $\alpha_2>0$, $\alpha_i=0$ for $i=-5,\ldots,1$.  

If $S$ is a finite, strongly connected set with 0 the leftmost vertex, then $6\in S$, since $0$ must be reachable from the right. We consider possible sets $S\cap[0,6]$. There are 32 subsets of $[0,6]$ that contain 0 and 6; however, $S$ must contain either 2 or 3, since there must be edges from 0 to other sets in $S$ and nothing to the left of 0 is allowed. Similarly, if $S$ contains 1, then it must contain either 3 or 4, and if $S$ contains 2, then it must contain either 4 or 5. This eliminates 12 of the 32 possibilities, leaving 20 possibilities for $S\cap[0,6] $. Of these, we first consider two candidates, $\{0,3,6\}$ and $\{0,2,4,6\}$. Both of these are strongly connected, and so $\beta_{\{0,3,6\}}=2\alpha_{-6}+3\alpha_2+\alpha_3$ and $\beta_{\{0,2,4,6\}}=3\alpha_{-6}+\alpha_2+4\alpha_3$ both provide upper bounds for $\kappa_0$. Depending on the values of the $\alpha_i$, either can be lower than the other. The set $\{0,2,3,4,6\}$ is also strongly connected, but has $\beta_{\{0,2,3,4,6\}}=4\alpha_{-6}+2\alpha_2+3\alpha_3$. Thus, if $\alpha_2\geq\alpha_3$, then $\beta_{\{0,2,4,6\}}<\beta_{\{0,2,3,4,6\}}$, and if $\alpha_2\leq\alpha_3$, then $\beta_{\{0,3,6\}}<\beta_{\{0,2,3,4,6\}}$. One can simply check that other 17 of the possible sets $D=S\cap[0,6]$ either have $\beta_D>\beta_{\{0,3,6\}}$ for all possible values of the $\alpha_i$, 
$\beta_D>\beta_{\{0,2,4,6\}}$ for all possible values of the $\alpha_i$,
or
$\beta_D>\beta_{\{0,2,3,4,6\}}$ for all possible values of the $\alpha_i$. By Claim \ref{claim:GammaClaim3}, this implies that $\beta_S\geq\min(\beta_{\{0,3,6\}},\beta_{\{0,2,4,6\}})$. Therefore, 
\begin{align*}
    \kappa_0
    &=\min(2\alpha_{-6}+3\alpha_2+\alpha_3,3\alpha_{-6}+\alpha_2+4\alpha_3)
    \\
    &=\min(\beta_{\{0,3,6\}},\beta_{\{0,2,4,6\}}).
\end{align*}
\end{exmp}

In all five of the above examples, there is always a set $S$ minimizing $\beta_S$ that represents a single, simple loop. The exit time from $S$ is the first time the walk stops repeating this loop. Thus, if $\kappa_0\leq1$, then there is a single loop that the walk is expected to repeat infinitely many times before deviating from it. 

In the nearest-neighbor case, treated in Example \ref{kappa0example1}, $\kappa_0\leq1$ means the walk is expected to repeat the loop $0\to1\to0$ infinitely many times before ever taking a different step (and, likewise, the walk is expected to repeat the loop $0\to-1\to0$ infinitely many times before ever stepping to $1$). This does not mean the only finite traps are sets of the form $\{x,x+1\}$. For example, it is also the case that $\beta_{[0,5]}\leq1$, so that the walk is expected to spend an infinite amount of time in $[0,5]$ before leaving it, regardless of the precise path (and even if transition probabilities at sites 1,2,3, and 4 are conditioned to be moderate). But there are no finite traps ``worse" (in the sense of finite moments of quenched expected exit time) than the set $\{0,1\}$. 

In fact, for nearest-neighbor RWDE on $\Z^d$, pairs of adjacent vertices are always the worst finite traps, and if $\kappa_0\leq1$, then the walk is expected to bounce back and forth between 0 and one other vertex infinitely many times before doing anything else \cite{Sabot&Tournier2016}.

Our other examples so far match this trend in a sense; although the worst traps are not necessarily pairs of vertices, the worst traps are loops, and $\kappa_0\leq1$ means there is a loop that the walk is expected to iterate infinitely many times before doing anything else. 
\begin{itemize}
    \item In Example \ref{kappa0example7}, one such loop is $0\to0$. 
    \item In Example \ref{kappa0example2}, one such loop is $0\to3\to6\to4\to2\to0$. 
    \item In Example \ref{kappa0example3}, one such loop is $0\to1\to0$. 
    \item In Example \ref{kappa0example4}, one such loop is $0\to1\to0$ (if $\beta_{\{0,1\}}\leq1$) or $0\to2\to0$ (if $\beta_{\{0,2\}}\leq1$).
    \item In Example \ref{kappa0example5}, one such loop is $0\to3\to6\to0$ (if $\beta_{\{0,3,6\}}<1$) or $0\to2\to4\to6$ (if $\beta_{\{0,2,4,6\}}<1$). 
\end{itemize}

Our next example shows that unlike in the nearest-neighbor case on $\Z^d$, there are parameters where the strongest finite traps never represent just one loop. In particular, one can find cases where $\kappa_0\leq1$, so there are finite traps in which the walk is expected to be stuck for an infinite amount of time, but there is no single loop that the walk is expected to iterate infinitely many times before deviating from it.

\begin{exmp}\label{kappa0example9}
$L=R=2$, $\alpha_{-1}=\alpha_0=0$, $\alpha_1>0$.
\end{exmp}

A finite, strongly connected set with 0 as its leftmost point will necessarily contain 2, since 0 must be reachable from the right. Thus, by claim \ref{claim:GammaClaim3}, for any finite, strongly connected $S$, $\beta_S\geq\min(\beta_{\{0,1,2\}},\beta_{\{0,2\}})$. Now $\{0,2\}$ and $\{0,1,2\}$ are already strongly connected, so $\kappa_0=\min(\beta_{\{0,1,2\}},\beta_{\{0,2\}})$. The minimum may be achieved on either set, depending on the $\alpha_i$. 

\begin{figure}
    \centering
    \includegraphics[width=6.5in]{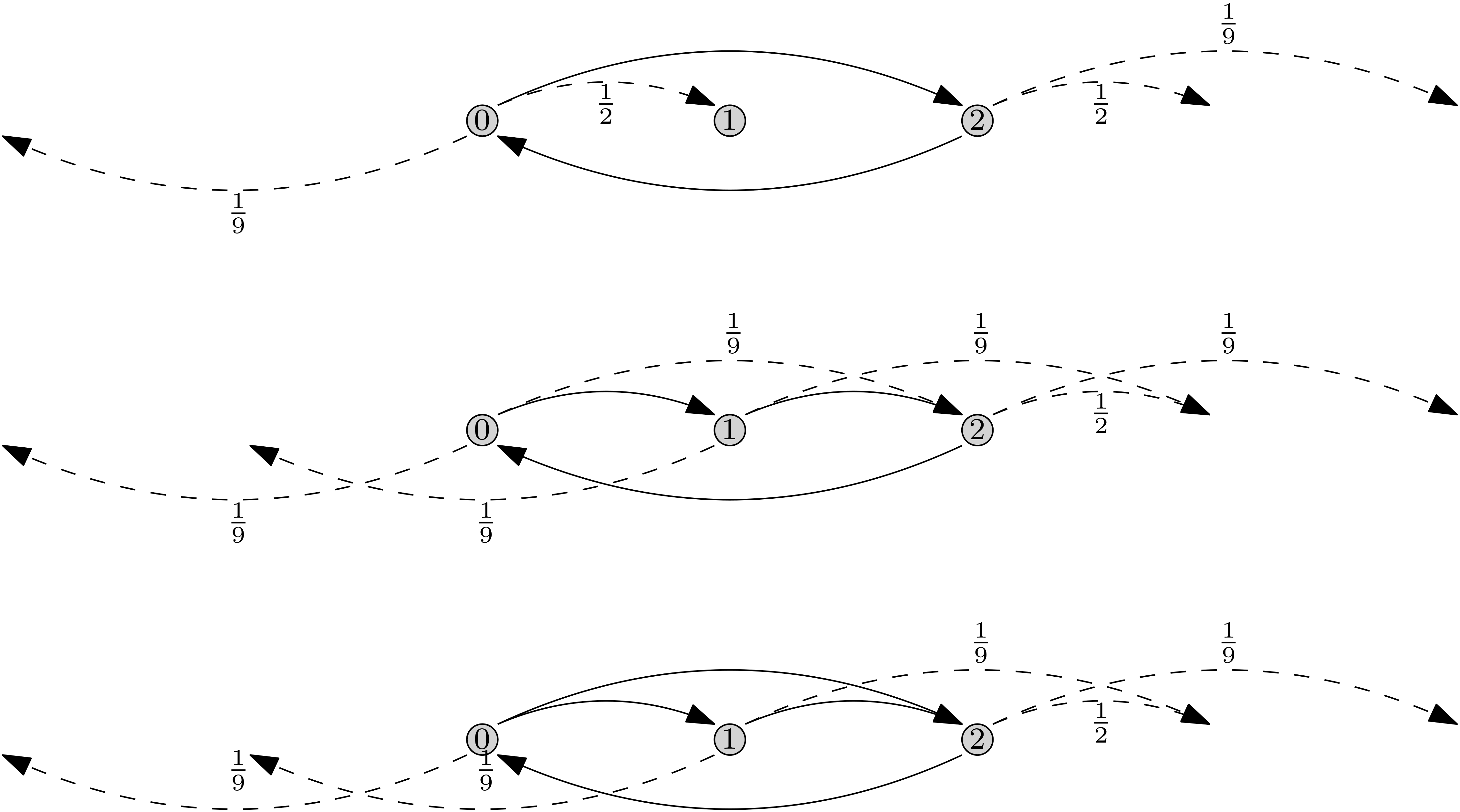}
    \caption{The top shows the weights exiting the loop $0\to2\to0$. The middle shows the weights exiting the loop $0\to1\to2\to0$. The bottom shows the weights exiting the union of these two loops, or the set $\{0,1,2\}$.}
    \label{fig:kappa0figure}
\end{figure}

We now examine a case where $\kappa_0\leq1$, but there are no loops that the walk is expected to iterate infinitely many times before doing anything else. Suppose $\alpha_{-2}=\alpha_{2}=\frac19$, and $\alpha_1=\frac12$. Then $\beta_{\{0,2\}}=\frac{11}{9}>1$, and $\beta_{\{0,1,2\}}=\frac{17}{18}<1$. Thus, $\kappa_0=\frac{17}{18}$, and a walk started from 0 is expected to spend an infinite amount of time in $\{0,1,2\}$ before exiting. However, because $\beta_{\{0,2\}}=\frac{11}{9}>1$, the expected exit time from $\{0,2\}$ is finite, so the walk is not expected to iterate the loop $0\to2\to0$ infinitely many times before deviating from it. Moreover, the walk is not expected to iterate the loop $0\to1\to2\to0$ infinitely many times before deviating it, but to see this, we must use the original formulation of Tournier's lemma from \cite{Tournier2009}. The formulation there is in terms of sets of edges rather than vertices. The edges that are not in the loop $0\to1\to2$ but have tails in the vertex set touched by this loop have weights that add up to $\frac{19}{18}>1$. Hence \cite[Theorem 1]{Tournier2009} implies that the expected time to deviate from this set of edges (and thus from the loop $0\to1\to2\to0$) is finite. Nevertheless, the weight exiting $\{0,1,2\}$ is $\frac{17}{18}<\frac12$, so the walk is expected to stick to the vertex set $\{0,1,2\}$, and thus to the pair of loops $0\to2\to0$ and $0\to1\to2\to0$, for an infinite amount of time before doing anything else. See Figure \ref{fig:kappa0figure}.

Note that in this example, $\beta_{\{0,1,2\}}<\beta_{\{0,2\}}$. This shows that in Claim \ref{claim:GammaClaim3}, the assumption that $x$ is to the right or left of $x$ is really needed.
 
Our next example presents a similar phenomenon: the walk is not expected to get stuck in any one loop for an infinite amount of time, but the walk is expected to spend an infinite amount of time in a set of vertices. In the previous example, the vertex set $S$ minimizing $\beta_S$ can have all of its vertices hit by one loop, the loop $0\to1\to2\to0$, but this loop alone does not have a trapping effect as strong as the whole set $S$. In the next example, there is no single loop that can hit all the vertices in the minimizing $S$, so our formulation of Tournier's lemma in terms of vertices only is enough to see that there is no single loop the walk is expected to traverse infinitely many times before straying from it. The next example also presents a calculation of $\kappa_0$ for a situation where it is less straightforward than the others we've examined.

\begin{exmp}\label{kappa0example6}
$L=16$, $R=5$, $\alpha_{-16},\alpha_2,\alpha_5>0$, all other $\alpha_i=0$. 

In this case, there are four possible values for $\kappa_0$, three of which can be attained by sets of vertices representing single loops, but one of which cannot. We will show that $\kappa_0$ is attained by one of the following four sets:

\begin{itemize}
    \item ${S_1}=\{0,2,4,6,8,10,12,14,16\}$. This set represents a loop that steps up by 2s from 0 to 16 and then jumps back to 0. $\beta_{S_1}=8\alpha_{-16}+   \alpha_2+9\alpha_5$.
    \item ${S_2}=\{0,5,10,15,16,20,25,30,32\}$.  The set ${S_2}$ represents a loop that steps up by 5s from 0 to 30, then steps to 32 and jumps back to 16 and then to 0. This is one of 28 loops that all step up by 5 six times, up by 2 once, and down by 16 twice, having vertex set $S$ with leftmost point 0. All such loops have the same associated $\beta_{S}=\beta_{S_2}=7\alpha_{-16}+   8\alpha_2+     3\alpha_5$. 
    \item $S_3=\{0,5,10,12,14,16\}$. The set ${S_3}$ represents a loop that steps up by 5s from 0 to 10, then by 2s from 0 to 16, then jumps back to 0. This is one of 10 loops that all have vertex set $S\subset[0,16]$, all of which have $\beta_{S}=\beta_{S_3}=    5\alpha_{-16}+   3\alpha_2+     4\alpha_5$.
    \item ${S_4}=\{0,2,4,5,6,7,8,9,10,11,12,14,16\}$. This set does {\em not} represent one single loop; in fact, it represents all 10 loops that stay within $[0,16]$.  $\beta_{S_4}=    12\alpha_{-16}+   2\alpha_2+     5\alpha_5$.
\end{itemize}
One can check that any of $\beta_{S_1}$, $\beta_{S_2}$, $\beta_{S_3}$, and $\beta_{S_4}$ can be the smallest, depending on the $\alpha_i$. We will show
\begin{equation*}
    \kappa_0=\min(\beta_{S_1},\beta_{S_2},\beta_{S_3},\beta_{S_4}). 
\end{equation*}

To confirm that these are the possible values for $\kappa_0$, let $S$ be a finite, strongly connected set with leftmost point 0, and we will show that $\beta_S$ is at least as large as one of these values. First, we note
\begin{equation*}
    \beta_S=x_{-16}\alpha_{-16}+x_2\alpha_2+x_5\alpha_5,
\end{equation*}
where $x_i=x_i(S):=\#\{z\in S:z+i\notin S\}\geq1$. Since 0 must be reachable in one step from a vertex to its right, $16\in S$.

\begin{claim}\label{claim:16proofclaim1}
$x_{-16}\geq 5$.
\end{claim}

$S$ must contain $0$ and $16$, so there must be a path $\sigma$ from $0$ to $[16,\infty)$ that does not leave $S$. Since the only step to the left is down 16, all steps in this path must be to the right (and must therefore be of length 2 or 5). If this path includes two or more steps of length 2, then $S\cap[0,15]$ must have at least 5 elements. But for each $z\in S\cap[0,15]$, $z-16\notin S$, so each element of $S\cap[0,15]$ contributes 1 to $x_{-16}$. Hence $x_{-16}\geq5$. On the other hand, if the path $\sigma$ includes no steps or one step of length 2, then the path $\sigma$ includes four vertices in $S\cap[0,15]$, does not include 1 or 4, and lands on either 17 or 20. Since 1 and 4 are not in $\sigma$, either 17 or 20 contributes 1 to $x_{-16}$, in which case we then have $x_{-16}\geq5$, or else 1 or 4 is in $S$, in addition to the four vertices from $\sigma$ that are in $[0,15]$, so that $|S\cap[0,15]|\geq5$, and so $x_{-16}\geq5$. This proves our claim. 

\begin{claim}\label{claim:16proofclaim2}
If $x_2=1$, then $x_{-16}\geq8$ and $x_5\geq9$.
\end{claim}

To see this, note that if $x_2=1$, then $S$ includes only even vertices; otherwise, the rightmost odd vertex and the rightmost even vertex would each contribute 1 to $x_2$, giving $x_2\geq2$. Moreover, since $S$ contains $0$ and $16$, it must contain every even vertex between, in order to prevent any vertex other than the rightmost from contributing to $x_2$. The 8 even vertices $z=0$ through $z=14$ each have $z-16\notin S$, so $x_{-16}\geq8$, and the 9 even vertices $z=0$ through $z=16$ each have $z+5\notin S$, so $x_5\geq9$. 

\begin{claim}\label{claim:16proofclaim3}
$x_5\geq3$. 
\end{claim}

To see this, note that the rightmost vertex from each equivalence class $\pmod 5$ will contribute 1 to $x_5$. We already have $0,16\in S$, so the equivalence classes $0$ and $1$ are represented. But the equivalence class 1 is only reachable from equivalence class 2 (via a downward step of length 16) and from equivalence class 4 (via an upward step of length 2). Hence $S$ must contain an element from one of the equivalence classes 2 or 4 $\pmod 5$, and therefore at least three equivalence classes are represented, so $x_5\geq3$. 

\begin{claim}\label{claim:16proofclaim4}
If $x_5=3$, then $x_{-16}\geq7$, and $x_2\geq8$. 
\end{claim}

We first note that since each equivalence class contributes only 1 to $x_5$, all elements in each equivalence class must form an unbroken arithmetic progression from the lowest to the highest. That is, letting $z_i^{\text{least}}$ and $z_i^{\text{greatest}}$ be, respectively, the least and greatest $z$ such that $z\equiv i \pmod 5$ and $z\in S$, we have $\{z_i^{\text{least}},z_i^{\text{least}}+5,z_i^{\text{least}}+10,\ldots,z_i^{\text{greatest}}\}\subset S$. We now examine two separate cases.

{\em Case 1: $S$ contains elements from equivalence classes 0,1, and 2 $\pmod5$.}

Since $S$ contains no elements from equivalence class 4, equivalence $z_1^{\text{least}}$ can only be reached from equivalence class 2, which occurs via a leftward step of length $16$. Thus $z_1^{\text{least}}+16\in S$.
On the other hand, since $S$ contains no elements from equivalence class 3, equivalence class 2 can only be reached from equivalence class 0, via a rightward step of length 2. Thus $z_2^{\text{least}}-2\in S$.

Therefore, the path
\begin{equation*}
    0\to5\to10\to\cdots\to (z_2^{\text{least}}-2)\to z_2^{\text{least}}
    \to\cdots\to (z_1^{\text{least}}+16)\to z_1^{\text{least}}\to\cdots\to 16\to0
\end{equation*}
is in $S$. All steps marked out by ellipses are upward steps of length 5. It follows that this is a path of length 9, since any other number of steps would result in ending at a point other than 0. 

All the vertices in equivalence class 0 contribute 1 to $x_{-16}$, since stepping down by 16 would reach a vertex in equivalence class 4, which cannot be in $S$. Vertices from the path that are in equivalence class 2, other than $z_1^{\text{least}}+16$, are less than $z_1^{\text{least}}+16$, and so stepping down by 16 reaches a vertex that is in equivalence class 1 but not in $S$. And all vertices from the path that are in equivalence class 1, other than 16, are less than 16, so stepping down by 16 reaches a vertex not in $S$. Thus all but two of the vertices from the path shown will contribute 1 to $x_{-16}$, and therefore $x_{-16}\geq7$.

Now all the vertices in equivalence classes 1 or 2 contribute 1 to $x_2$, since stepping to the right by 2 reaches a vertex in equivalence class 3 or 4. And vertices in equivalence class 0 that are less than $z_2^{\text{least}}-2$ also contribute to $x_2$, since stepping to the right by 2 reaches a vertex in equivalence class 2 but less than $z_2^{\text{least}}$. Thus, all but one of the vertices shown in this path contribute to $x_2$, so $x_2\geq8$.

{\em Case 2: $S$ contains elements from equivalence classes 0,1, and 4 $\pmod5$.}

By a similar argument to that given in Case 1, the path
\begin{equation*}
    0\to5\to10\to\cdots\to (z_4^{\text{least}}+16)\to z_4^{\text{least}}
    \to\cdots\to (z_1^{\text{least}}-2)\to z_1^{\text{least}}\to\cdots\to 16\to0
\end{equation*}
is in $S$. All steps marked out by ellipses are upward steps of length 5. It follows that this is a path of length 9, since any other number of steps would result in ending at a point other than 0. Now, for $z=0,5,10,z_4^{\text{least}}+11$, we have $z-16\equiv 4\pmod 5$, but $z-16<z_4^{\text{least}}$, so $z-16\notin S$ and $z$ contributes 1 to $x_{-16}$. Moreover, $z_4^{\text{least}}$ and every subsequent vertex are all less than 16 (except, of course, for 16 itself), so they all contribute 1 to $x_{-16}$. Thus, $x_{-16}\geq7$.

Moreover, all vertices in equivalence class 0 or 1 contribute 1 to $x_2$, since $S$ has no vertices in equivalence class 2 or 4. And all but one of the vertices $z$ in equivalence class 4 are strictly less than $z_1^{\text{least}}-2$, so that $z+2\notin S$. Hence all but one of the vertices in the loop contribute to $x_2$, so $x_2\geq8$.

\begin{claim}\label{claim:16proofclaim5}
If $x_2=2$, then $x_{5}\geq5$. 
\end{claim}
If $x_2=2$, then $S$ contains even and odd elements (because the only even upward jumps are of length 2, a strongly connected $S$ with only even elements would have $x_2=1$). It therefore must contain every even number, from its least even number to its greatest even number. In particular, it must contain $S_1=[0,16]$. This is enough to include at least one representative from every equivalence class $\pmod 5$. The greatest element of $S$ in each of these equivalence classes contributes 1 to $x_5$, so $x_5\geq5$. 

\begin{claim}\label{claim:16proofclaim7}
If $x_2=2$, then $x_5+x_{-16}\geq 17$ and $x_{-16}\geq9$. 
\end{claim}

We have already established that if $x_2=2$, then $S$ contains $S_1=[0,16]$ and at least one odd number.
Now $S_1$ has $x_{-16}=8$ and $x_5=9$. The odd number will also contribute 1 to $x_{-16}$, giving the bound $x_{-16}\geq9$. 
The set $S_1$ includes 8, which is in equivalence class 3 $\pmod 5$, and two elements of each of the equivalence classes 0,1,2, and 4. Each of the equivalence classes must contribute at least 1 to $x_5$, and for any of the classes 0,1,2, or 4 to avoid contributing 2, the odd number in between the two even numbers from that equivalence class must be contained in $S$. This saves 1 from $x_5$ but adds 1 to $x_{-16}$, thus keeping $x_5+x_{-16}\geq17$.

\begin{claim}\label{claim:16proofclaim8}
$\beta_S\geq\min(\beta_{S_1},\beta_{S_2},\beta_{S_3},\beta_{S_4})$. 
\end{claim}

We know $\beta_S$ must have $x_5\geq3$ by Claim \ref{claim:16proofclaim3}. By Claim \ref{claim:16proofclaim4}, if $x_5=3$, then $\beta_S\geq7\alpha_{-16}+8\alpha_2+3\alpha_5=\beta_{S_2}$. Now suppose $x_5\geq 4$. If $x_2=1$, then $\beta_S\geq 8\alpha_{-16}+\alpha_2+9\alpha_5=\beta_{S_1}$ by Claim \ref{claim:16proofclaim2}. Now consider the case $x_2=2$. Then $x_5+x_{-16}\geq17$ by Claim \ref{claim:16proofclaim7}. If $\alpha_5>\alpha_{-16}$, then since $x_5\geq5$ by Claim \ref{claim:16proofclaim5}, we have $\beta_S\geq 12\alpha_{-16}+2\alpha_2+5\alpha_5=\beta_{S_4}$. On the other hand, if $\alpha_{-16}>\alpha_5$, then by Claim \ref{claim:16proofclaim7}, $\beta_S\geq 9\alpha_{-16}+2\alpha_2+8\alpha_5>8\alpha_{-16}+\alpha_2+9\alpha_5=\beta_{S_1}$. Now, if $x_2\geq3$, then by the assumption that $x_5\geq4$ and by Claim \ref{claim:16proofclaim1}, we have $\beta_S\geq5\alpha_{-16}+3\alpha_2+4\alpha_5=\beta_{S_3}$. This proves our final claim. 

Now suppose the weights are $\alpha_{-16}=\frac{1}{67}$, $\alpha_2=\frac{15}{67}$, and $\alpha_5=\frac{5}{67}$. We can check that $\kappa_0=12\alpha_{-16}+2\alpha_2+5\alpha_5=1$, achieved on the set $S_4=\{0,2,4,5,6,7,8,9,10,11,12,14,16\}$, and that this is strictly less than $\beta_{S_1}$, $\beta_{S_2}$, and $\beta_{S_3}$. By the proof of Claim \ref{claim:16proofclaim7}, any set $S$ with $x_{-16}=12,x_2=2,x_5=5$ must contain a translation of $\{0,2,4,5,6,7,8,9,10,11,12,14,16\}$, and a  so there is no possibility that a set $S$ which we did not consider, and which represents a single loop, {\em also} achieves $\beta_S=1$. This means that the walk is expected to spend an infinite amount of time in the set $\{0,2,4,5,6,7,8,9,10,11,12,14,16\}$ before ever leaving it, but there is no single loop that the walk is expected to take infinitely many times before deviating from it. 
\end{exmp}

\section{Notation}\label{append:notation}

Here we collect notation that is used throughout the paper as a convenient reference.

{\bf General}

\begin{itemize}
    \item $\N=\{1,2,3,\ldots\}$. $\N_0=\{0,1,2,\ldots\}$. $\R^{\geq0}=\{x\in\R:x\geq0\}$. $\R^{>0}=\{x\in\R:x>0\}$.
    \item RWRE stands for random walk(s) in random environment(s). RWDE stands for random walk(s) in Dirichlet environment(s). 
    \item An {\em environment} $\omega$ on a countable vertex set $V$, the set $\Omega_V$ of environments on $V$, and the measurable space $(\Omega_V,\mathcal{F}_V)$, are defined in Section \ref{subsec:model}. 
    \item $\omega^x=(\omega(x,x+y))_{y\in\Z}$ is the environment $\omega$ viewed at site $x$ only. For a set $S\subseteq\Z$, $\omega^S=(\omega^x)_{x\in S}$. \item An $\omega$ with a subscript (e.g., $\omega_1$) or $\omega'$ is usually used to denote a specific environment when comparing multiple environments. 
    \item Conditions (\hyperlink{cond:C1}{C1}), (\hyperlink{cond:C2}{C2}) (\hyperlink{cond:C3}{C3}), and (\hyperlink{cond:C4}{C4}) are defined right after the statement of Theorem \ref{thm:transiencecriterion}. 
    \item $\Delta_I:=\{(x_i)_{i\in I}:\sum_{i\in I}x_i=1\}$ is the simplex of a finite set $I$.
    \item $\theta^x$ is the left shift operator on environments $\omega$, defined by $\theta^x\omega(a,b)=\omega(x+a,x+b)$.
    \item We use interval notation to denote sets of consecutive integers in the state space $\Z$, rather than subsets of $\R$. For example, $[1,\infty)$ denotes the set of integers to the right of 0. However, we make one exception, using $[0,1]$ to denote the set of all real numbers from 0 to 1.  
\end{itemize}

{\bf Graphs}

\begin{itemize}
    \item A {weighted directed graph} $\Hg=(V,E,W)$ is a vertex set $V$ with an edge set $E\subseteq V\times V$, and a weight function $w:E\to\R^{>0}$.
    \item  If $e=(x,y)\in E$, we say that $e$ is an edge from $x$ to $y$, and we say the {\em head} of $e$ is $\overline{e}=y$ and the {\em tail} of $e$ is $\underline{e}=x$.
    \item For a vertex $x\in V$, the {\em divergence} of $x$ in $\Hg$ is $\text{div}(x)=\sum_{\overline{e}=x}w(e)-\sum_{\underline{e}=x}w(e)$. If the divergence is zero for all $x$, we say the graph $\Hg$ has zero divergence.
    \item A set $S\subset V$ is {\em strongly connected} if for all $x,y\in S$, there is a path from $x$ to $y$ in $\Hg$ using only vertices in $S$.
    \item For a set $S\subseteq V$, $\beta_S$ is the sum of the weights of all edges exiting $S$. See \eqref{betaS}.
    \item $\G$ is our ``main graph." It vertex set $\Z$, edge set $\{(x,y)\in \Z\times\Z:-L\leq y-x\leq R,\alpha_{y-x}>0\}$, and weight function $(x,y)\mapsto\alpha_{y-x}$. 
    \item $\G_M$ is the finite graph with vertices $[0,M]$ that looks like $\G$ in the middle but is modified to have zero divergence at the endpoints. It is defined in the proof of Theorem \ref{thm:transiencecriterion}, and used there as well as in the proof of Lemma \ref{lem:InfiniteTimeReversal}.
    \item $\G_+$ is the half-infinite graph with vertex set $[0,\infty)$. It looks like $\G$ except near 0, where it looks like $\G_M$. It can be thought of as a limit of $\G_M$ as $M\to\infty$. It is defined at the beginning of Section \ref{subsec:backtracking}. Its crucial property is given by Lemma \ref{lem:InfiniteTimeReversal}, which is used in the proofs of Proposition \ref{prop:InfiniteTraps}, Proposition \ref{prop:FewOscillations}, and Theorem \ref{thm:MainMomentTheorem}. 
\end{itemize}

{\bf Parameters}

\begin{itemize}
    \item $(\alpha_i)_{i=-L}^R$ are Dirichlet parameters for random transition probability vectors. It is assumed that $\alpha_{-L}$ and $\alpha_R$ are positive.
    \item $d^{\totheright}=\sum_{i=1}^Ri\alpha_i$, and $d^{\totheleft}=\sum_{i=-L}^{-1}|i|\alpha_i$. 
    \item $c^{\totheright}=\sum_{i=1}^R\alpha_i$, and $c^{\totheleft}=\sum_{i=-L}^{-1}\alpha_i$. 
    \item $\kappa_1=d^{\totheright}-d^{\totheleft}=\sum_{i=-L}^Ri\alpha_i$.
    \item $\kappa_0=\inf\{\beta_S:S\subset \Z\text{ finite, strongly connected}\}$ is the minimum weight exiting a finite, strongly connected subset of $\G$. 
    \item $m_0$ is an integer large enough that every interval of length $m_0$ is strongly connected in $\G$, and also large enough that $m_0\geq\max(L,R)$.
\end{itemize}

{\bf Probability measures}

\begin{itemize}
    \item For a given weighted directed graph $\Hg=(V,E,W)$, $P_{\Hg}$ is the Dirichlet law on environments corresponding to $\Hg$; that is, the measure on $\Omega_V$ under which transition probabilities at the various vertices $x\in V$ are independent, and for each vertex $x\in V$, $(\omega(x,\overline{e}))_{\underline{e}=x}$ is distributed according to a Dirichlet distribution with parameters $(w(e))_{\underline{e}=x}$. (Or, if $V\subset\Z$, we let $P_{\Hg}$ be any measure on $\Omega_{\Z}$ whose marginals on $\Omega_V$ are as described.)
    \item $P_{\G}$ is the main measure on $\Omega_{\Z}$ that we study in this paper. That is, $P_{\G}$ is the Dirichlet measure on environments corresponding to $\G$, so that for each $x$, the transition probability vector $(\omega(x,y))_{\alpha_{y-x}>0}$ is distributed as a Dirichlet random vector with parameters $(\alpha_i)_{\alpha_i>0}$.
\end{itemize}

{\bf Functions of a walk}

For all of the functions below, we often suppress the argument of the function as is traditional with random variables. Sometimes, however, we leave the argument in when it is necessary for clarity.
\begin{itemize}
    \item $H_x({\bf X})=\inf\{n\geq0:X_n=x\}$ is the first time the walk $({\bf X})$ hits $x$. We often suppress the ${\bf X}$. 
    \item $\tilde{H}_x({\bf X})=\inf\{n>0:X_n=x\}$ is the first nonzero time ${\bf X}$ hits $x$.
    \item $N_x({\bf X})=\#\{n\in\N_0:X_n=x\}$ is the total amount of time ${\bf X}$ spends at $x$.  
    \item For a subset $S\subset\Z$, $H_S({\bf X})=\min_{x\in S}H_x({\bf X})$, and $N_S({\bf X})=\sum_{x\in S}N_x({\bf X})$.
    \item $N_x^S({\bf X})=\#\{0\leq n\leq H_{S^c}:X_n=x\}$ is the amount of time ${\bf X}$ spends at $x$ before leaving the subset $S$.
    \item $N_{x,y}({\bf X})=\#\big\{n\geq0:X_n=x,\sup\{k<n:X_k=y\}>\sup\{k<n:X_k=x\}\big\}$ is the number of times the walk hits $x$ after more recently having hit $y$, or the number of ``trips from $y$ to $x$".
    \item $N_{x,y}'({\bf X}):=\#\big\{n\in\N_0:X_n\leq x,\sup\{j<n:X_j\geq y\}>\sup\{j<n:X_j \leq x\}\big\}$ is the number of trips leftward across $[x,y]$.
    \end{itemize}

{\bf Cascade and bi-infinite walk}

\begin{itemize}
    \item A {\em cascade} is a set of finite (continuous-time) walks, one started at each point in $\Z$. The walk starting at each point terminates when it reaches or passes the next multiple of $R$.
    \item ${\bf X}^a=(X_n^a)_{n=0}^{\infty}$ is the infinite walk obtained by concatenating the finite walk starting at $a$ with the walk starting at the point where it terminates, the walk started at the point where that one terminates, and so on. 
    \item $\BIX=(\biX_n)_{n\in\Z}$ is the bi-infinite walk obtained by this process.
    \item $\overline{N}_x$, $\overline{H}_x$, $\overline{N}_{x,y}$, and so on are defined analogously to $N_x$, $H_x$, $N_{x,y}$, and so on, but with $n\geq0$ replaced with $n\in\Z$.
\end{itemize}

\bibliographystyle{plain}
\bibliography{default}







\end{document}